\documentclass[12pt, reqno]{amsart}
\usepackage{amsmath, amsthm, amscd, amsfonts, amssymb, graphicx, color}
\usepackage[bookmarksnumbered, colorlinks, plainpages]{hyperref}

\newtheorem{theorem}{Theorem}[section]
\newtheorem{lemma}[theorem]{Lemma}
\newtheorem{proposition}[theorem]{Proposition}
\newtheorem{corollary}[theorem]{Corollary}
\theoremstyle{definition}
\newtheorem{definition}[theorem]{Definition}
\newtheorem{example}[theorem]{Example}

\newtheorem{condition}[theorem]{Condition}
\theoremstyle{remark}
\newtheorem{remark}[theorem]{Remark}
\numberwithin{equation}{section}

\begin{document}
\setcounter{page}{1}

\title[The Refined Sobolev Scale, Interpolation, and Elliptic Problems]
{\large The Refined Sobolev Scale,\\ Interpolation, and Elliptic Problems}

\author[V.A. Mikhailets, A.A. Murach]
{Vladimir A. Mikhailets and Aleksandr A. Murach}

\address{Institute of Mathematics, National Academy of Sciences of Ukraine,
3, Tere\-shch\-en\-kiv\-ska Str, 01601 Kyiv-4, Ukraine.}

\email{mikhailets@imath.kiev.ua, murach@imath.kiev.ua}



\thanks{The authors were partly supported by grant no. 01-01-02 of National Academy of
Sciences of Ukraine (under the joint Ukrainian--Russian project of NAS of Ukraine
and Russian Foundation of Basic Research)}

\begin{abstract}
The paper gives a detailed survey of recent results on elliptic problems in Hilbert
spaces of generalized smoothness. The latter are the isotropic H\"ormander spaces
$H^{s,\varphi}:=B_{2,\mu}$, with
$\mu(\xi)=\langle\xi\rangle^{s}\varphi(\langle\xi\rangle)$ for
$\xi\in\mathbb{R}^{n}$. They are parametrized by both the real number $s$ and the
positive function $\varphi$ varying slowly at $+\infty$ in the Karamata sense. These
spaces form the refined Sobolev scale, which is much finer than the Sobolev scale
$\{H^{s}\}\equiv\{H^{s,1}\}$ and is closed with respect to the interpolation with a
function parameter. The Fredholm property of elliptic operators and elliptic
boundary-value problems is preserved for this new scale. Theorems of various type
about a solvability of elliptic problems are given. A~local refined smoothness is
investigated for solutions to elliptic equations. New sufficient conditions for the
solutions to have continuous derivatives are found. Some applications to the
spectral theory of elliptic operators are given.
\end{abstract}

\maketitle

\section{Introduction}\label{sec1}

\noindent In the theory of partial differential equations, the questions concerning
the existence, uniqueness, and regularity of solutions are in the focus of
investigations. Note that the regularity properties are usually formulated in terms
of the belonging of solutions to some standard classes of function spaces. Thus, the
finer a used scale of spaces is calibrated, the sharper and more informative results
will be.

In contrast to the ordinary differential equations with smooth coefficients, the
above questions are complicated enough. Indeed, some linear partial differential
equations with smooth coefficients and right-hand sides are known to have no
solutions in a neighbourhood of a given point, even in the class of distributions
\cite{Lewy57}, \cite[Sec.~6.0 and 7.3]{Hermander63}, \cite[Sec.~13.3]{Hermander83}.
Next, certain homogeneous equations (specifically, of elliptic type) with smooth but
not analytic coefficients have nontrivial solutions supported on a compact set
\cite{Plis54}, \cite[Theorem 13.6.15]{Hermander83}.  Hence, the nontrivial
null-space of this equation cannot be removed by any homogeneous boundary-value
conditions; i.e., the operator of an arbitrary boundary-value problem is not
injective. Finally, the question about the regularity of solutions is not simple
either. For example, it is known \cite[Ch. 4, Notes]{GilbargTrudinger98} that
$$
\triangle u=f\in\,C(\Omega)\nRightarrow\,u\in C^{\,2}(\Omega),
$$
with $\triangle$ being the Laplace operator, and $\Omega$ being an arbitrary
Euclidean domain.

These questions have been investigated most completely for the elliptic equations,
systems, and boundary-value problems. This was done in the 1950s and 1960s by
S.~Agmon, A.~Douglis, L.~Nirenberg, M.S.~Agranovich, A.C.~Dynin, Yu.M.~Berezansky,
S.G.~Krein, Ya.A.~Roitberg, F.~Browder, L.~Hermander, J.-L.~Lions, E.~Magenes,
M.~Schechter, L.N.~Slobodetsky, V.A.~Solonnikov, L.R.~Vo\-le\-vich and some others
(see, e.g., M.S.~Agranovich's surveys \cite{Agranovich94, Agranovich97} and the
references given therein). Note that the elliptic equations and problems have been
investigated in the classical scales of H\"older spaces (of noninteger order) and
Sobolev spaces (both of positive and negative orders).

The fundamental result of the theory of elliptic equations consists in that they
generate bounded and Fredholm operators (i.e., the operators with finite index)
between appropriate function spaces. For instance, let $Au=f$ be an elliptic linear
differential equation of order $m$ given a closed smooth manifold $\Gamma$. Then the
operator
$$
A:\,H^{s+m}(\Gamma)\rightarrow H^{s}(\Gamma),\quad s\in\mathbb{R},
$$
is bounded and Fredholm. Moreover, the finite-dimensional spaces formed by solutions
to homogeneous equations $Au=0$ and $A^{+}v=0$ both lie in $C^{\infty}(\Gamma)$.
Here $A^{+}$ is the formally adjoint operator to $A$, whereas $H^{s+m}(\Gamma)$ and
$H^{s}(\Gamma)$ are inner product Sobolev spaces over $\Gamma$ and of the orders
$s+m$ and $s$ respectively. It follows from this that the solution $u$ have an
important regularity property on the Sobolev scale, namely
\begin{equation}\label{eq1.1}
(f\in H^{s}(\Gamma)\;\;\mbox{for some}\;\;s\in\mathbb{R})\;\Rightarrow\; u\in
H^{s+m}(\Gamma).
\end{equation}

If the manifold has a boundary, then the Fredholm operator is generated by an
elliptic boundary-value problem for the equation $Au=f$, specifically, by the
Dirichlet problem.

Some of these results were extended by H.~Triebel \cite{Triebel95, Triebel83} and
the second author \cite{94UMJ12, 94Dop12} of the survey to finer scales of
function spaces, namely the Nikolsky--Besov, Zygmund, and Lizorkin--Triebel scales.

The results mentioned above have various applications in the theory of differential
equations, mathematical physics, the spectral theory of differential operators; see
M.S.~Agranovich's surveys \cite{Agranovich94, Agranovich97} and the references
therein.

As for applications, especially to the spectral theory, the case of Hilbert spaces
is of the most interest. Until recently, the Sobolev scale had been a unique
scale of Hilbert spaces in which the properties of elliptic operators were
investigated systematically. However, it turns out that this scale is not
fine enough for a number of important problems.

We will give two representative examples. The first of them concerns with the
smoothness properties of solutions to the elliptic equation $Au=f$ on the manifold
$\Gamma$. According to Sobolev's Imbedding Theorem, we have
\begin{equation}\label{eq1.2}
H^{\sigma}(\Gamma)\subset C^{r}(\Gamma)\;\;\Leftrightarrow\;\;\sigma>r+n/2,
\end{equation}
where the integer $r\geq0$ and $n:=\dim\Gamma$. This result and property
\eqref{eq1.1} allow us to investigate the classical regularity of the solution $u$.
Indeed, if $f\in H^{s}(\Gamma)$ for some $s>r-m+n/2$, then $u\in
H^{s+m}(\Gamma)\subset C^{r}(\Gamma)$. However, this is not true for $s=r-m+n/2$;
i.e., the Sobolev scale cannot be used to express unimprovable sufficient conditions
for belonging of the solution $u$ to the class $C^{r}(\Gamma)$. An analogous
situation occurs in the theory of elliptic boundary-value problems too.

The second demonstrative example is related to the spectral theory. Suppose that the
differential operator $A$ is of order $m>0$, elliptic on $\Gamma$, and self-adjoint
on the space $L_{2}(\Gamma)$. Given a function $f\in\nobreak L_{2}(\Gamma)$,
consider the spectral expansion
\begin{equation}\label{eq1.3}
f=\sum_{j=1}^{\infty}\;c_{j}(f)\,h_{j},
\end{equation}
where $(h_{j})_{j=1}^{\infty}$ is a complete orthonormal system of eigenfunctions of
$A$, and $c_{j}(f)$ is the Fourier coefficient of $f$ with respect to $h_{j}$. The
eigenfunctions are enumerated so that the absolute values of the corresponding
eigenvalues form a (nonstrictly) increasing sequence. According to the
Menshov--Rademacher theorem, which are valid for the general orthonormal series too,
the expansion \eqref{eq1.3} converges almost everywhere on $\Gamma$ provided that
\begin{equation}\label{eq1.4}
\sum_{j=1}^{\infty}\,|c_{j}(f)|^{2}\,\log^{2}(j+1)<\infty.
\end{equation}
This hypotheses cannot be reformulated in equivalent manner in terms of the
belonging of $f$ to Sobolev spaces because
$$
\|f\|_{H^{s}(\Gamma)}^{2}\,\asymp\,\sum_{j=1}^{\infty}\,|c_{j}(f)|^{2}\,j^{2s}
$$
for every $s>0$. We may state only that the condition "$f\in H^{s}(\Gamma)$ for some
$s>0$" implies convergence of the series \eqref{eq1.3} almost everywhere on
$\Gamma$. This condition does not adequately express the hypotheses \eqref{eq1.4} of
the Menshov--Rademacher theorem.

In 1963 L.~H\"ormander \cite[Sec. 2.2]{Hermander63} proposed a broad and informative
generalization of the Sobolev spaces in the category of Hilbert spaces (also see
\cite[Sec. 10.1]{Hermander83}). He introduced spaces that are parametrized by a
general enough weight function, which serves as an analog of the differentiation
order or smoothness index used for the Sobolev spaces. In particular, H\"ormander
considered the following Hilbert spaces
\begin{gather}\label{eq1.5}
B_{2,\mu}(\mathbb{R}^{n}):=
\bigl\{u:\,\mu\,\mathcal{F}u\in L_{2}(\mathbb{R}^{n})\bigr\}, \\
\|u\|_{B_{2,\mu}(\mathbb{R}^{n})}:=\|\mu\,\mathcal{F}u\|_{L_{2}(\mathbb{R}^{n})}.
\notag
\end{gather}
Here $\mathcal{F}u$ is the Fourier transform of a tempered distribution $u$ given on
$\mathbb{R}^{n}$, and $\mu$ is a weight function of $n$ arguments.

In the case where
$$
\mu(\xi)=\langle\xi\rangle^{s},\quad
\langle\xi\rangle:=(1+|\xi|^{2})^{1/2},\quad\xi\in\mathbb{R}^{n},\quad
s\in\mathbb{R},
$$
we have the Sobolev space $B_{2,\mu}(\mathbb{R}^{n})=H^{s}(\mathbb{R}^{n})$ of
differentiation order $s$.

The H\"ormander spaces occupy a central position among the spaces of generalized
smoothness, which is characterized by a function parameter, rather than a number.
These spaces are under various and profound investigations; a good deal of the work
was done in the last decades. We refer to G.A.~Kalyabin and P.I.~Lizorkin's survey
\cite{KalyabinLizorkin87}, H.~Triebel's monograph \cite[Sec.~22]{Triebel01}, the
recent papers by A.M.~Caetano and H.-G.~Leopold \cite{CaetanoLeopold06}, W.~Farkas,
N.~Jacob, and R.L.~Schilling \cite{FarkasJacobScilling01b}, W.~Farkas and
H.-G.~Leopold \cite{FarkasLeopold06}, P.~Gurka and B.~Opic \cite{GurkaOpic07},
D.D.~Haroske and S.D.~Moura \cite{HaroskeMoura04, HaroskeMoura08}, S.D.~Moura
\cite{Moura01}, B.~Opic and W.~Trebels \cite{OpicTrebels00}, and references given
therein. Various classes of spaces of generalized smoothness appear naturally in
embedding theorems for function spaces, the theory of interpolation of function
spaces, approximation theory, the theory of differential and pseudodifferential
operators, theory of stochastic processes; see the monographs by D.D.~Haroske
\cite{Haroske07}, N.~Jacob \cite{Jacob010205}, V.G.~Maz'ya and T.O.~Shaposhnikova
\cite[Sec.~16]{MazyaShaposhnikova09}, F.~Nicola and L.~Rodino \cite{NicolaRodino10},
B.P.~Paneah \cite{Paneah00}, A.I.~Stepanets \cite[Ch.~I, \S~7]{Stepanets87},
\cite[Part~I, Ch.~3, Sec. 7.1]{Stepanets05}, and also the papers by F.~Cobos and
D.L.~Fernandez \cite{CobosFernandez88}, C.~Merucci \cite{Merucci84}, M.~Schechter
\cite{Schechter67} devoted to the interpolation of function spaces, and the papers
by D.E.~Edmunds and H.~Triebel \cite{EdmundsTriebel98, EdmundsTriebel99},
V.A.~Mikhailets and V.M.~Molyboga \cite{MikhailetsMolyboga09, MikhailetsMolyboga11,
MikhailetsMolyboga12} on spectral theory of some elliptic operators appearing in
mathematical physics.

Already in 1963 L.~H\"ormander applied the spaces \eqref{eq1.5} and more general
Banach spaces $B_{p,\mu}(\mathbb{R}^{n})$, with $1\leq p\leq\infty$, to an
investigation of regularity properties of solutions to the partial differential
equations with constant coefficients and to some classes of equations with varying
coefficients. However, as distinct from the Sobolev spaces, the H\"ormander spaces
have not got a broad application to the general elliptic equations on manifolds and
to the elliptic boundary-value problems. This is due to the lack of a reasonable
definition of the H\"ormander spaces on smooth manifolds (the definition should be
independent of a choice of local charts covering the manifold) an on the absence of
analytic tools fit to use these spaces effectively.

Such a tool exists in the Sobolev spaces case; this is the interpolation of
spaces. Namely, an arbitrary fractional order Sobolev space can be obtained by the
interpolation of a certain couple of integer order Sobolev spaces. This fact
essentially facilitates both the investigation of these spaces and proofs of various
theorems of the theory of elliptic equations because the boundedness and the
Fredholm property (if the defect is invariant) preserve for linear operators under
the interpolation.

Therefore it seems reasonable to distinguish the H\"ormander spaces that are
obtained by the interpolation (with a function parameter) of couples of Sobolev
spaces; we will consider only inner product spaces. For this purpose we introduce the
following class of isotropic spaces
\begin{equation}\label{eq1.6}
H^{s,\varphi}(\mathbb{R}^{n}):=B_{2,\mu}(\mathbb{R}^{n})\quad\mbox{for}\quad
\mu(\xi)={\langle\xi\rangle^{s}\varphi(\langle\xi\rangle)}.
\end{equation}
Here the number parameter $s$ is real, whereas the positive function parameter
$\varphi$ varies slowly at $+\infty$ in the Karamata sense
\cite{BinghamGoldieTeugels89, Seneta76}. (We may assume that $\varphi$ is constant
outside of a neighbourhood of $+\infty$.) For example, $\varphi$ is admitted to be
the logarithmic function, its arbitrary iteration, their real power, and a product
of these functions.

The class of spaces \eqref{eq1.6} contains the Sobolev Hilbert scale
$\{H^{s}\}\equiv\{H^{s,1}\}$, is attached to it by the number parameter, but is
calibrated much finer than the Sobolev scale. Indeed,
$$
H^{s+\varepsilon}(\mathbb{R}^{n})\subset H^{s,\varphi}(\mathbb{R}^{n})\subset
H^{s-\varepsilon}(\mathbb{R}^{n})\quad\mbox{for every}\quad\varepsilon>0.
$$
Therefore the number parameter $s$ defines the main (power) smoothness, whereas the
function parameter $\varphi$ determines an additional (subpower) smoothness on the
class of spaces \eqref{eq1.6}. Specifically, if $\varphi(t)\rightarrow\infty$ (or
$\varphi(t)\rightarrow0$) as $t\rightarrow\infty$, then $\varphi$ determines an
additional positive (or negative) smoothness. Thus, the parameter $\varphi$ refines
the main smoothness $s$. Therefore the class of spaces \eqref{eq1.6} is naturally
called the refined Sobolev scale.

This scale possesses the following important property: every space
$H^{s,\varphi}(\mathbb{R}^{n})$ is a result of the interpolation, with an
appropriate function parameter, of the couple of Sobolev spaces
$H^{s-\varepsilon}(\mathbb{R}^{n})$ and $H^{s+\delta}(\mathbb{R}^{n})$, with
$\varepsilon,\delta>0$. The parameter of the interpolation is a function that varies
regularly (in the Karamata sense) of index $\theta\in(0,\,1)$ at $+\infty$; namely
$\theta:=\varepsilon/(\varepsilon+\delta)$. Moreover, the refined Sobolev scale
proves to be closed with respect to this interpolation.

Thus, every H\"ormander space $H^{s,\varphi}(\mathbb{R}^{n})$ possesses the
interpolation property with respect to the Sobolev Hilbert scale. This means that
each linear operator bounded on both the spaces $H^{s-\varepsilon}(\mathbb{R}^{n})$
and $H^{s+\delta}(\mathbb{R}^{n})$ is also bounded on
$H^{s,\varphi}(\mathbb{R}^{n})$. The interpolation property plays a decisive role
here; namely, it permits us to establish some important properties of the refined
Sobolev scale. They enable this scale to be applied in the theory of elliptic
equations. Thus, we can prove with the help of the interpolation that each space
$H^{s,\varphi}(\mathbb{R}^{n})$, as the Sobolev spaces, is invariant with respect to
diffeomorphic transformations of $\mathbb{R}^{n}$. This permits the space
$H^{s,\varphi}(\Gamma)$ to be well defined over a smooth closed manifold $\Gamma$
because the set of distributions and the topology in this space does not depend on a
choice of local charts covering~$\Gamma$. The spaces $H^{s,\varphi}(\mathbb{R}^{n})$
and $H^{s,\varphi}(\Gamma)$ are useful in the theory of elliptic operators on
manifolds and in the theory of elliptic boundary-value problems; these spaces are
present implicitly in a number of problems appearing in calculus.

Let us dwell on some results that demonstrate advantages of the introduced scale as
compared with the Sobolev scale. These results deal with the examples considered
above. As before, let $A$ be an elliptic differential operator given on $\Gamma$,
with $m:=\mathrm{ord}\,A$. Then $A$ sets the bounded and Fredholm operators
$$
A:\,H^{s+m,\varphi}(\Gamma)\rightarrow H^{s,\varphi}(\Gamma)\quad\mbox{for all}\quad
s\in\mathbb{R},\;\;\varphi\in\mathcal{M}.
$$
Here $\mathcal{M}$ is the class of slowly varying function parameters $\varphi$ used
in \eqref{eq1.6}. Note that the differential operator $A$ leaves invariant the
function parameter $\varphi$, which refines the main smoothness $s$. Besides, we
have the following regularity property of a solution to the elliptic equation
$Au=f$:
$$
(f\in H^{s,\varphi}(\Gamma)\;\;\mbox{for
some}\;\;s\in\mathbb{R},\;\varphi\in\mathcal{M})\;\Rightarrow\; u\in
H^{s+m,\varphi}(\Gamma).
$$

For the refined Sobolev scale, we have the following sharpening of Sobolev's
Imbedding Theorem: given an integer $r\geq0$ and function $\varphi\in\mathcal{M}$,
the embedding $H^{r+n/2,\varphi}(\Gamma)\subset C^{r}(\Gamma)$ is equivalent to that
\begin{equation}\label{eq1.7}
\int\limits_{1}^{\infty}\frac{dt}{t\,\varphi^{2}(t)}<\infty.
\end{equation}
Therefore, if $f\in H^{r-m+n/2,\varphi}(\Gamma)$ for some parameter
$\varphi\in\nobreak\mathcal{M}$ satisfying \eqref{eq1.7}, then the solution $u\in
C^{r}(\Gamma)$.

Similar results are also valid for the elliptic systems and elliptic boundary-value
problems.

Now let us pass to the analysis of the spectral expansion \eqref{eq1.3} convergence.
We additionally suppose that the operator $A$ is of order $m>0$ and is unbounded and
self-adjoint on the space $L_{2}(\Gamma)$. Condition \eqref{eq1.4} for the
convergence of \eqref{eq1.3} almost everywhere on $\Gamma$ is equivalent to the
inclusion
$$
f\in H^{0,\varphi}(\Gamma),\quad\mbox{with}\quad\varphi(t):=\max\{1,\log t\}.
$$
The latter is much wider than the condition "$f\in H^{s}(\Gamma)$ for some
$s>\nobreak0$". We can also similarly represent conditions for unconditional
convergence almost everywhere or convergence in the H\"older space $C^{r}(\Gamma)$,
with integral $r\geq0$.

The above and some other results show that the refined Sobolev scale is helpful and
convenient. This scale can be used in different topics of the modern analysis as
well; see, e.g., the articles by M.~Hegland \cite{Hegland95, Hegland10}, P.~Math\'e
and U.~Tautenhahn \cite{MatheTautenhahn06}.

This paper is a detailed survey of our recent articles [78--94, 101--108], which are
summed up in the monograph \cite{MikhailetsMurach10} published in Russian in 2010.
In them, we have built a theory of general elliptic (both scalar and matrix)
operators and elliptic boundary-value problems on the refined Sobolev scales of
function spaces.

Let us describe the survey contents in greater detail. The paper consists of 13
sections.

Section~\ref{sec1} is Introduction, which we are presenting now.

Section~\ref{sec2} is preliminary and contains a necessary information about
regularly varying functions and about the interpolation with a function parameter.
Here we distinguish important Theorem \ref{th2.5}, which gives a description of all
interpolation parameters for the category of separable Hilbert spaces.

In Section~\ref{sec3}, we consider the H\"ormander spaces, give a definition of the
refined Sobolev scale, and study its properties. Among them, we especially note the
interpolation properties of this scale, formulated as Theorems \ref{th3.4} and
\ref{th3.5}. They are very important for applications.

Section~\ref{sec4} deals with uniformly elliptic pseudodifferential operators that
are studied on the refined Sobolev scale over $\mathbb{R}^{n}$. We get an a priory
estimate for a solution of the elliptic equation and investigate an interior
smoothness of the solution. As an application, we obtain a sufficient condition for
the existence of continuous bounded derivatives of the solution.

Next in Section~\ref{sec5}, we define a class of H\"ormander spaces, the refined
Sobolev scale, over a smooth closed manifold. We give three equivalent definitions
of these spaces: local (in terms of local properties of distributions),
interpolational (by means of the interpolation of Sobolev spaces with an appropriate
function parameter), and operational (via the completion of the set of infinitely
smooth functions with respect to the norm generated by a certain function of the
Beltrami--Laplace operator). These definitions are similar to those used for the
Sobolev spaces. We study properties of the refined Sobolev scale over the closed
manifold. Important applications of these results are given in Sections \ref{sec6}
and \ref{sec7}.

Section~\ref{sec6} deals with elliptic pseudodifferential operators on a closed
manifold. We show that they are Fredholm (i.e. have a finite index) on appropriate
couples of H\"ormander spaces. As in Section~\ref{sec4}, a priory estimates for
solutions of the elliptic equations are obtained, and the solutions regularity is
investigated. Using elliptic operators, we give equivalent norms on H\"ormander
spaces over the manifold.

In Section~\ref{sec7}, we investigate a convergence of spectral expansions
corresponding to elliptic normal operators given on the closed manifold. We find
sufficient conditions for the following types of the convergence: almost everywhere,
unconditionally almost everywhere, and in the space $C^{k}$, with integral $k\geq0$.
These conditions are formulated in constructive terms of the convergence on some
function classes, which are H\"ormander spaces.

Section~\ref{sec8} deals with the classes of H\"ormander spaces that relate to the
refined Sobolev scale and are given over Euclidean domains being open or closed. For
these classes, we study interpolation properties, embeddings, traces, and riggings
of the space of square integrable functions with H\"ormander spaces. The results of
this section are applied in next Sections \ref{sec9}--\ref{sec12}, where a regular
elliptic boundary-value problem is investigated in appropriate H\"ormander spaces.

In Section~\ref{sec9}, this problem is studied on the one-sided refined Sobolev
scale. We show that the problem generates a Fredholm operator on this scale. We
investigate some properties of the problem; namely, a priory estimates for solutions
and local regularity are given. Moreover, a sufficient condition for the weak
solution to be classical is found in terms of H\"ormander spaces.

Section~\ref{sec10} deals with semihomogeneous elliptic boundary-value problems.
They are considered on H\"ormander spaces which form an appropriate two-sided
refined Sobolev scale. We show that the operator corresponding to the problem is
bounded and Fredholm on this scale.

In Sections \ref{sec11}--\ref{sec12}, we give various theorems about a solvability
of nonhomogeneous regular elliptic boundary-value problems in H\"ormander spaces of
an arbitrary real main smoothness. Developing the methods suggested by
Ya.A.~Roitberg \cite{Roitberg96} and J.-L.~Lions, E.~Magenes \cite{LionsMagenes72},
we establish a certain generic theorem and a wide class of individual theorems on
the solvability. The generic theorem is featured by that the domain of the elliptic
operator does not depend on the coefficients of the elliptic equation and is common
for all boundary-value problems of the same order. Conversely, the individual
theorems are characterized by that the domain depends essentially on the
coefficients, even of the lower order derivatives. In Section~\ref{sec11}, we
elaborate on Roitberg's approach in connection with H\"ormander spaces and then
deduce the generic theorem about the solvability of elliptic boundary-value problems
on the two-sided refined Sobolev scale modified in the Roitberg sense.

Section~\ref{sec12} is devoted to J.-L.~Lions and E.~Magenes' approach, which we
develop for various Hilbert scales consisting of Sobolev or H\"ormander spaces. For
the space of right-hand sides of an elliptic equation, we find a sufficiently
general condition under which the operator of the problem is bounded and Fredholm
(see key Theorems \ref{th12.1} and \ref{th12.4}). As a consequence, we obtain new
various individual theorems on the solvability of elliptic boundary-value problems
considered in Sobolev or H\"ormander spaces, both nonweighted and weighted.

In final Section~\ref{sec13}, we indicate application of H\"ormander spaces to other
important classes of elliptic problems. They are nonregular boundary-value problems,
parameter-elliptic problems, certain mixed elliptic problems, elliptic systems and
corresponding boundary-value problems.

It is necessary to note that some results given in the survey are new even for the
Sobolev spaces. These results are Theorem \ref{th10.1} in the case of half-integer
$s$ and individual Theorems \ref{th12.1}, \ref{th12.2}, and \ref{th12.3}.

In addition, note that we have also investigated a certain class of H\"ormander
spaces, which is wider than the refined Sobolev scale. Interpolation properties of
this class are studied and then applied to elliptic operators \cite{08Collection1,
09Dop3, MikhailetsMurach10, arXiv:1106.2049, 09UMJ3, arXiv:1202.6156}. It is
remarkable that this class consists of all the Hilbert spaces which possess the
interpolation property with respect to the Sobolev Hilbert scale. These results fall
beyond the limits of our survey.

\section{Preliminaries}\label{sec2}

In this section we recall some important results concerning the regularly varying
functions and the interpolation with a function parameter of couples of Hilbert
spaces. These results will be necessary for us in the sequel.

\subsection{Regularly varying functions}\label{sec2.1}

We recall the following notion.

\begin{definition}\label{def2.1}
A positive function $\psi$ defined on a semiaxis $[b,+\infty)$ is said to be
regularly varying of index $\theta\in\mathbb{R}$ at $+\infty$ if $\psi$ is Borel
measurable on $[b_{0},+\infty)$ for some number $b_{0}\geq b$ and
$$
\lim_{t\rightarrow+\infty}\;\frac{\psi(\lambda\,t)}{\psi(t)}=
\lambda^{\theta}\quad\mbox{for each}\quad \lambda>0.
$$
A function regularly varying of the index $\theta=0$ at $+\infty$ is called slowly
varying at $+\infty$.
\end{definition}

The theory of regularly varying functions was founded by Jovan Karamata
\cite{Karamata30a, Karamata30b, Karamata33} in the 1930s. These functions are
closely related to the power functions and have numerous applications, mainly due to
their special role in Tauberian-type theorems (see the monographs
\cite{BinghamGoldieTeugels89, GelukHaan87, Haan70, Maric00, Resnick87, Seneta76} and
references therein).

\begin{example}\label{ex2.1}
The well-known standard case of functions regularly varying of the index $\theta$ at
$+\infty$ is
\begin{equation}\label{eq2.1}
\psi(t):=t^{\theta}\,(\log t)^{r_{1}}\,(\log\log t)^{r_{2}} \ldots (\log\ldots\log
t)^{r_{k}}\quad\mbox{for}\quad t\gg1
\end{equation}
with arbitrary parameters $k\in\mathbb{Z}_{+}$ and $r_{1},
r_{2},\ldots,r_{k}\in\mathbb{R}$. In the case where $\theta=0$ these functions form
the logarithmic multiscale, which has a number of applications in the theory of
function spaces.
\end{example}

We denote by $\mathrm{SV}$ the set of all functions slowly varying at $+\infty$. It is
evident that $\psi$ is a function regularly varying at $+\infty$ of index $\theta$ if and
only if $\psi(t)=t^{\theta}\varphi(t)$, $t\gg1$, for some function
$\varphi\in\mathrm{SV}$. Thus, the investigation of regularly varying functions is
reduced to the case of slowly varying functions.

The study and application of regularly varying functions are based on two
fundamental theorems: the Uniform Convergence Theorem and Representation Theorem.
They were proved by Karamata \cite{Karamata30a} in the case of continuous functions
and, in general, by a number of authors later (see the monographs cited above).

\begin{theorem}[Uniform Convergence Theorem]\label{th2.1}
Let $\varphi\in\mathrm{SV}$; then $\varphi(\lambda
t)/\varphi(t)\rightarrow\nobreak1$ as $t\rightarrow+\infty$ uniformly on each
compact $\lambda$-set in $(0,\infty)$.
\end{theorem}

\begin{theorem}[Representation Theorem]\label{th2.2}
A function $\varphi$ belongs to $\mathrm{SV}$ if and only if it can be written in
the form
\begin{equation}\label{eq2.2}
\varphi(t)=\exp\Biggl(\beta(t)+\int\limits_{b}^{\:t}\frac{\alpha(\tau)}{\tau}\,d\tau\Biggr),
\quad t\geq b,
\end{equation}
for some number $b>0$, continuous function $\alpha:[b,\infty)\rightarrow\mathbb{R}$
approaching zero at $\infty$, and Borel measurable bounded function
$\beta:\nobreak[b,\infty)\rightarrow\mathbb{R}$ that has the finite limit at
$\infty$.
\end{theorem}

The Representation Theorem implies the following sufficient condition for a function
to be slowly varying at infinity \cite[Sec. 1.2]{Seneta76}.

\begin{theorem}\label{th2.3}
Suppose that a function $\varphi:(b,\infty)\rightarrow(0,\infty)$ has a continuous
derivative and satisfies the condition $t\varphi\,'(t)/\varphi(t)\rightarrow0$ as
$t\rightarrow \infty$. Then $\varphi\in\mathrm{SV}$.
\end{theorem}

Using Theorem~\ref{th2.3} one can give many interesting examples of slowly varying
functions. Among them we mention the following.

\begin{example}\label{ex2.2}
Let $\varphi(t):=\exp\psi(t)$, with $\psi$ being defined according to \eqref{eq2.1},
where $\theta=0$ and $r_{1}<1$. Then $\varphi\in\mathrm{SV}$.
\end{example}

\begin{example}\label{ex2.3}
Let $\alpha,\beta,\gamma\in\mathbb{R}$, $\beta\neq0$, and $0<\gamma<1$. We set
$\omega(t):=\alpha+\beta\sin\,\log^{\gamma}t$  and $\varphi(t):=(\log
t)^{\omega(t)}$ for $t>1$. Then $\varphi\in\mathrm{SV}$.
\end{example}

\begin{example}\label{ex2.4}
Let $\alpha,\beta,\gamma\in\mathbb{R}$, $\alpha\neq0$, $0<\gamma<\beta<1$, and
$$
\varphi(t):=\exp(\alpha(\log t)^{1-\beta}\,\sin\log^{\gamma}t)\quad\mbox{for} \quad
t>1.
$$
Then $\varphi\in\mathrm{SV}$.
\end{example}

The last two examples show that a function $\varphi$ varying slowly at $+\infty$ may
exhibit infinite oscillation, that is
$$
\liminf_{t\rightarrow+\infty}\,\varphi(t)=0\quad\mbox{and}\quad
\limsup_{t\rightarrow+\infty}\,\varphi(t)=+\infty.
$$

We will use regularly varying functions as parameters when we define certain Hilbert
spaces. If the function parameters are equivalent in a neighbourhood of $+\infty$,
we get the same space up to equivalence of norms. Therefore it is useful to
introduced the following notion \cite[p.~90]{08MFAT1}.

\begin{definition}\label{def2.2}
We say that a positive function $\psi$ defined on a semiaxis $[b,+\infty)$ is
quasiregularly varying of index $\theta\in\mathbb{R}$ at $+\infty$ if there exist a
number $b_{1}\geq b$ and a function $\psi_{1}:[b_{1},+\infty)\rightarrow
(0,+\infty)$ regularly varying of the same index $\theta\in\mathbb{R}$ at $+\infty$
such that $\psi\asymp\psi_{1}$ on $[b_{1},+\infty)$. A function quasiregularly
varying of the index $\theta=0$ at $+\infty$ is called quasislowly varying at
$+\infty$.
\end{definition}

As usual, the notation $\psi\asymp\psi_{1}$ on $[b_{1},+\infty)$ means that the
functions $\psi$ and $\psi_{1}$ are equivalent there, that is both the functions
$\psi/\psi_{1}$ and  $\psi_{1}/\psi$ are bounded on $[b_{1},+\infty)$.

We denote by $\mathrm{QSV}$ the set of all functions varying quasislowly at
$+\infty$. It is evident that $\psi$ is quasiregularly varying of the index $\theta$
at $+\infty$ if and only if $\psi(t)=t^{\theta}\varphi(t)$, $t\gg1$, for some
function $\varphi\in\mathrm{QSV}$.

We note the following properties of the class $\mathrm{QSV}$.

\begin{theorem}\label{th2.4}
Let $\varphi,\chi\in\mathrm{QSV}$. The next assertions are true:
\begin{enumerate}
\item[i)] There is a function $\varphi_{1}\in
C^{\infty}((0;+\infty))\cap\mathrm{SV}$ such that $\varphi\asymp\varphi_{1}$ in a
neighbourhood of $+\infty$.
\item[ii)] If $\theta>0$, then both $t^{-\theta}\varphi(t)\rightarrow0$
and $t^{\theta}\varphi(t)\rightarrow+\infty$ as $t\rightarrow+\infty$.
\item[iii)] All the functions $\varphi+\chi$, $\varphi\,\chi$, $\varphi/\chi$ and
$\varphi^{\sigma}$, with $\sigma\in\mathbb{R}$, belong to $\mathrm{QSV}$.
\item[iv)] Let $\theta\geq0$, and in the case where $\theta=0$ suppose that
$\varphi(t)\rightarrow+\infty$ as $t\rightarrow+\infty$. Then the composite function
$\chi(t^{\theta}\varphi(t))$ of $t$ belongs to $\mathrm{QSV}$.
\end{enumerate}
\end{theorem}

Theorem~\ref{th2.4} are known for slowly varying functions, even with the strong
equivalence $\varphi(t)\sim\varphi_{1}(t)$ as $t\rightarrow+\infty$ being in
assertion i); see, e.g., \cite[Sec. 1.3]{BinghamGoldieTeugels89} and \cite[Sec.
1.5]{Seneta76}. This implies the case when $\varphi,\chi\in\mathrm{QSV}$
\cite[p.~91]{08MFAT1}.

\subsection{The interpolation with a function parameter of Hilbert spaces}\label{sec2.2}

It is a natural generalization of the classical interpolation method by J.-L.~Lions
and S.G.~Krein (see, e.g., \cite[Ch.~IV, \S~9]{FunctionalAnalysis72} and
\cite[Ch.~1, Sec. 2 and 5]{LionsMagenes72}) to the case when a general enough
function is used as an interpolation parameter instead of a number parameter. The
generalization appeared in the paper by C.~Foia\c{s} and J.-L.~Lions
\cite[p.~278]{FoiasLions61} and then was studied by W.F.~Donoghue \cite{Donoghue67},
E.I.~Pustyl`nik \cite{Pustylnik82}, V.I.~Ovchinnikov \cite[Sec.
11.4]{Ovchinnikov84}, and the authors \cite{08MFAT1}.

We recall the definition of this interpolation. For our purposes, it is sufficient to
restrict ourselves to the case of separable Hilbert spaces.

Let an ordered couple $X:=[X_{0},X_{1}]$ of complex Hilbert spaces $X_{0}$ and
$X_{1}$ be such that these spaces are separable and that the continuous dense
embedding $X_{1}\hookrightarrow X_{0}$ holds true. We call this couple admissible.
For the couple $X$ there exists an isometric isomorphism $J:X_{1}\leftrightarrow
X_{0}$ such that $J$ is a self-adjoint positive operator on the space $X_{0}$ with
the domain $X_{1}$ (see \cite[Ch.~1, Sec. 2.1]{LionsMagenes72} and \cite[Ch.~IV,
Sec. 9.1]{FunctionalAnalysis72}). The operator $J$ is said to be generating for the
couple $X$ and is uniquely determined by $X$.

We denote by $\mathcal{B}$ the set of all functions
$\psi:(0,\infty)\rightarrow(0,\infty)$ such that:
\begin{enumerate}
\item[a)] $\psi$ is Borel measurable on the semiaxis $(0,+\infty)$;
\item[b)] $\psi$ is bounded on each compact interval $[a,b]$ with $0<a<b<+\infty$;
\item[c)] $1/\psi$ is bounded on each set $[r,+\infty)$ with $r>0$.
\end{enumerate}

Let $\psi\in\mathcal{B}$. Generally, the unbounded operator $\psi(J)$ is defined in
the space $X_{0}$ as a function of $J$. We denote by $[X_{0},X_{1}]_{\psi}$ or
simply by $X_{\psi}$ the domain of the operator $\psi(J)$ endowed with the inner
product $(u,v)_{X_{\psi}}:=(\psi(J)u,\psi(J)v)_{X_{0}}$ and the corresponding norm
$\|u\|_{X_{\psi}}:=(u,u)_{X_{\psi}}^{1/2}$. The space $X_{\psi}$ is Hilbert and
separable.

\begin{definition}\label{def2.3}
We say that a function $\psi\in\mathcal{B}$ is an interpolation parameter if the
following property is fulfilled for all admissible couples $X=[X_{0},X_{1}]$,
$Y=[Y_{0},Y_{1}]$ of Hilbert spaces and an arbitrary linear mapping $T$ given on
$X_{0}$. If the restriction of the mapping $T$ to the space $X_{j}$ is a bounded
operator $T:X_{j}\rightarrow Y_{j}$ for each $j=0,\,1$, then the restriction of the
mapping $T$ to the space $X_{\psi}$ is also a bounded operator
$T:X_{\psi}\rightarrow Y_{\psi}$.
\end{definition}

Otherwise speaking, $\psi$ is an interpolation parameter if and only if the mapping
$X\mapsto X_{\psi}$ is an interpolation functor given on the category of all
admissible couples $X$ of Hilbert spaces. (For the notion of interpolation functor,
see, e.g., \cite[Sec. 2.4]{BerghLefstrem76} and \cite[Sec. 1.2.2]{Triebel95}) In the
case where $\psi$ is an interpolation parameter, we say that the space $X_{\psi}$ is
obtained by the interpolation with the function parameter $\psi$ of the admissible
couple $X$. Then the continuous dense embeddings $X_{1}\hookrightarrow
X_{\psi}\hookrightarrow X_{0}$ are fulfilled.

The classical result by J.-L. Lions and S.G.~Krein  consists in the fact that the
power function $\psi(t):=t^{\theta}$ is an interpolation parameter whenever
$0<\theta<1$; see \cite[Ch.~IV, \S~9, Sec.~3]{FunctionalAnalysis72} and \cite[Ch.~1,
Sec. 5.1]{LionsMagenes72}.

We have the following criterion for a function to be an interpolation parameter.

\begin{theorem}\label{th2.5}
A function $\psi\in\mathcal{B}$ is an interpolation parameter if and only if $\psi$
is pseudoconcave in a neighbourhood of $+\infty$, i.e. $\psi\asymp\psi_{1}$ for some
concave positive function $\psi_{1}$
\end{theorem}

This theorem follows from Peetre's results \cite{Peetre68} on interpolations
functions (see also the monograph \cite[Sec. 5.4]{BerghLefstrem76}). The
corresponding proof is given in \cite[Sec. 2.7]{08MFAT1}.

For us, it is important the next consequence of Theorem \ref{th2.5}.

\begin{corollary}\label{cor2.1}
Suppose a function $\psi\in\mathcal{B}$ to be quasiregularly varying of index
$\theta$ at $+\infty$, with $0<\theta<1$. Then $\psi$ is an interpolation parameter.
\end{corollary}

The direct proof of this assertion is given in \cite[Sec. 2]{06UMJ2}.

\section{H\"ormander spaces}\label{sec3}

In 1963 Lars H\"ormander \cite[Sec. 2.2]{Hermander63} introduced the spaces
$B_{p,\mu}(\mathbb{R}^{n})$, which consist of distributions in $\mathbb{R}^{n}$ and
are parametrized by a number $p\in[1,\infty]$ and a general enough weight function
$\mu$ of argument $\xi\in\mathbb{R}^{n}$; see also \cite[Sec. 10.1]{Hermander83}.
The number parameter $p$ characterizes integrability properties of the
distributions, whereas the function parameter $\mu$ describes their smoothness
properties. In this section, we recall the definition of the spaces
$B_{p,\mu}(\mathbb{R}^{n})$, some their properties, and an application to
constant-coefficient partial differential equations. Further we consider the
important case where the H\"ormander space $B_{p,\mu}(\mathbb{R}^{n})$ is Hilbert,
i.e. $p=2$, and $\mu$ is a quasiregularly varying function of $(1+|\xi|^{2})^{1/2}$
at infinity.

\subsection{The spaces $B_{p,\mu}(\mathbb{R}^{n})$}\label{sec3.1}

Let an integer $n\geq1$ and a parameter $p\in[1,\infty]$.  We use the following
conventional designations, where $\Omega$ is an nonempty open set in
$\mathbb{R}^{n}$, in particular $\Omega=\mathbb{R}^{n}$:

\begin{enumerate}
\item [a)] $L_{p}(\Omega):=L_{p}(\Omega,d\xi)$ is the Banach space of complex-valued
functions $f(\xi)$ of $\xi\in\Omega$ such that $|f|^{p}$ is integrable over $\Omega$
(if $p=\infty$, then f is essentially bounded in $\Omega$);
\item [b)] $C^{k}_{\mathrm{b}}(\Omega)$ is the Banach space of functions
$u:\Omega\rightarrow\nobreak\mathbb{C}$ having continuous and bounded derivatives of
order $\leq k$ on $\Omega$;
\item [c)] $C^{\infty}_{0}(\Omega)$ is the linear topological space of infinitely
differentiable functions $u:\mathbb{R}^{n}\rightarrow\mathbb{C}$ such that their
supports are compact and belong to $\Omega$; we will identify functions from
$C^{\infty}_{0}(\Omega)$ with their restrictions to $\Omega$;
\item [d)] $\mathcal{D}'(\Omega)$ is the linear topological space of all distributions
given in $\Omega$; we always suppose that distributions are antilinear
complex-valued functionals;
\item [e)] $\mathcal{S}'(\mathbb{R}^{n})$ is the linear topological Schwartz space of
tempered distributions given in $\mathbb{R}^{n}$;
\item [f)] $\widehat{u}:=\mathcal{F}u$ is the Fourier transform of a distribution
$u\in\mathcal{S}'(\mathbb{R}^{n})$; $\mathcal{F}^{-1}f$ is the inverse Fourier
transform of $f\in\mathcal{S}'(\mathbb{R}^{n})$;
\item [g)] $\langle\xi\rangle:=(1+|\xi|^{2})^{1/2}$ is a smoothed modulus of
$\xi\in\mathbb{R}^{n}$.
\end{enumerate}

Suppose a continuous function $\mu:\mathbb{R}^{n}\rightarrow(0,\infty)$ to be such
that, for some numbers $c\geq1$ and $l>0$, we have
\begin{equation}\label{eq3.1}
\frac{\mu(\xi)}{\mu(\eta)}\leq c\,(1+|\xi-\eta|)^{l}\quad\mbox{for
all}\quad\xi,\eta\in\mathbb{R}^{n}.
\end{equation}
The function $\mu$ is called a weight function.

\begin{definition}\label{def3.1}
The H\"ormander space $B_{p,\mu}(\mathbb{R}^{n})$ is a linear space of the
distributions $u\in\mathcal{S}'(\mathbb{R}^{n})$ such that the Fourier transform
$\widehat{u}$ is locally Lebesgue integrable on $\mathbb{R}^{n}$ and, moreover,
$\mu\,\widehat{u}\in L_{p}(\mathbb{R}^{n})$. The space $B_{p,\mu}(\mathbb{R}^{n})$
is endowed with the norm
$\|u\|_{B_{p,\mu}(\mathbb{R}^{n})}:=\|\mu\,\widehat{u}\|_{L_{p}(\mathbb{R}^{n})}$.
\end{definition}

The space $B_{p,\mu}(\mathbb{R}^{n})$ is complete and continuously embedded in
$\mathcal{S}'(\mathbb{R}^{n})$. If $1\leq p<\infty$, then this space is separable,
and the set $C^{\infty}_{0}(\mathbb{R}^{n})$ is complete in it \cite[Sec.
2.2]{Hermander63}. Of special interest is the $p=2$ case, when
$B_{p,\mu}(\mathbb{R}^{n})$ becomes a Hilbert space.

\begin{remark}\label{rem3.1}
H\"ormander assumes initially that $\mu$ satisfies a stronger condition than
\eqref{eq3.1}; namely, there exist some positive numbers $c$ and $l$ such that
\begin{equation}\label{eq3.2}
\frac{\mu(\xi)}{\mu(\eta)}\leq(1+c\,|\xi-\eta|)^{l}\quad\mbox{for
all}\quad\xi,\eta\in\mathbb{R}^{n}.
\end{equation}
But he notices that two sets of functions satisfying either \eqref{eq3.1} or
\eqref{eq3.2} lead to the same class of spaces $B_{p,\mu}(\mathbb{R}^{n})$ \cite[the
remark at the end of Sec. 2.1]{Hermander63}.
\end{remark}

The term `H\"ormander space' was suggested by H.~Triebel in \cite[Sec.
4.11.4]{Triebel95}.

The following H\"ormander's theorem establishes an important relation between the
spaces $B_{p,\mu}(\mathbb{R}^{n})$ and $C^{k}_{\mathrm{b}}(\mathbb{R}^{n})$
\cite[Sec. 2.2, Theorem 2.2.7]{Hermander63}.

\begin{theorem}[H\"ormander's Embedding Theorem]\label{th3.1}
Let $p,q\in[1,\infty]$, $1/p+1/q=\nobreak1$, and an integer $k\geq0$. Then the
condition
\begin{equation}\label{eq3.3}
\langle\xi\rangle^{k}\,\mu^{-1}(\xi)\in L_{q}(\mathbb{R}^{n},d\xi)
\end{equation}
entails the continuous embedding $B_{p,\mu}(\mathbb{R}^{n})\hookrightarrow
C^{k}_{\mathrm{b}}(\mathbb{R}^{n})$. Conversely, if
$$
\{u\in B_{p,\mu}(\mathbb{R}^{n}):\mathrm{supp}\,u\subset V\}\subset
C^{k}(\mathbb{R}^{n})
$$
for some nonempty open set $V\subseteq\mathbb{R}^{n}$, then \eqref{eq3.3} is valid.
\end{theorem}

The spaces $B_{p,\mu}(\mathbb{R}^{n})$ were applied by H\"ormander to investigation
of regularity properties of solutions to some partial differential equations (see
\cite[Ch.~IV, VII]{Hermander63} and \cite[Ch.~11, 13]{Hermander83}). We state one of
his results relating to elliptic equations \cite[Sec 7.4]{Hermander63}.

Let $\Omega$ be a nonempty open set in $\mathbb{R}^{n}$. In $\Omega$, consider a partial
differential equation $P(x,D)u=f$ of an order $r$ with coefficients belonging to
$C^{\infty}(\Omega)$. Introduce the local H\"ormander space over $\Omega$:
$$
B_{p,\mu}^{\mathrm{loc}}(\Omega):=\{f\in\mathcal{D}'(\Omega):\,\chi f\in
B_{p,\mu}(\Omega)\;\;\forall\;\;\chi\in C^{\infty}_{0}(\Omega)\}.
$$
Here $B_{p,\mu}(\Omega)$ is the space of restrictions of all the distributions $u\in
B_{p,\mu}(\mathbb{R}^{n})$ to~$\Omega$.

\begin{theorem}[H\"ormander's Regularity Theorem]\label{th3.2}
Let the operator $P(x,D)$ be elliptic in $\Omega$, and $u\in\mathcal{D}'(\Omega)$.
If $P(x,D)u\in B_{p,\mu}^{\mathrm{loc}}(\Omega)$ for some $p\in[1,\infty]$ and
weight function $\mu$, then $u\in B_{p,\mu_{r}}^{\mathrm{loc}}(\Omega)$ with
$\mu_{r}(\xi):=\langle\xi\rangle^{r}\mu(\xi)$.
\end{theorem}

For applications of the spaces $B_{p,\mu}(\mathbb{R}^{n})$, the Hilbert case of
$p=2$ is the most interesting. This case was investigated by B.~Malgrange
\cite{Malgrange57} and L.R.~Volevich, B.P.~Paneah \cite{VolevichPaneah65} (see also
Paneah's monograph \cite[Sec. 1.4]{Paneah00}). Specifically, if
$\mu(\xi)=\langle\xi\rangle^{s}$ for all $\xi\in\mathbb{R}^{n}$ with some
$s\in\mathbb{R}$, then $B_{2,\mu}(\mathbb{R}^{n})$ becomes the Sobolev inner product
space $H^{s}(\mathbb{R}^{n})$ of order $s$.

In what follows we will consider the isotropic H\"ormander inner product spaces
$B_{2,\mu}(\mathbb{R}^{n})$, with $\mu(\xi)$ being a radial function, i.e. depending
only on $\langle\xi\rangle$.

\subsection{The refined Sobolev scale}\label{sec3.2}

It useful to have a class of the H\"ormander inner product spaces
$B_{2,\mu}(\mathbb{R}^{n})$ that are close to the Sobolev spaces
$H^{s}(\mathbb{R}^{n})$ with $s\in\mathbb{R}$. For this purpose we choose
$\mu(\xi):=\langle\xi\rangle^{s}\varphi(\langle\xi\rangle)$ for some functions
$\varphi\in\mathrm{QSV}$; then $\mu$ is a quasiregularly varying function of
$\langle\xi\rangle$ at infinity of index $s$. In this case it is naturally to rename
the H\"ormander space $B_{2,\mu}(\mathbb{R}^{n})$ by
$H^{s,\varphi}(\mathbb{R}^{n})$. Let us formulate the corresponding definitions.
First we introduce the following set $\mathcal{M}\subset\mathrm{QSV}$ of function
parameters $\varphi$.

By $\mathcal{M}$ we denote the set of all functions
$\varphi:[1;+\infty)\rightarrow(0;+\infty)$ such that:

\begin{enumerate}
\item [a)] $\varphi$ is Borel measurable on $[1;+\infty)$;
\item [b)] $\varphi$ and $1/\varphi$ are bounded on every compact interval $[1;b]$, where
$1<b<+\infty$;
\item [c)] $\varphi\in\mathrm{QSV}$.
\end{enumerate}

It follows from Theorem \ref{th2.2} that $\varphi\in\mathcal{M}$ if and only if
$\varphi$ can be written in the form \eqref{eq2.2} with $b=1$ for some continuous
function $\alpha:[1,\infty)\rightarrow\mathbb{R}$ approaching zero at $+\infty$ and
Borel measurable bounded function $\beta:\nobreak[1,\infty)\rightarrow\mathbb{R}$.

Let $s\in\mathbb{R}$ and $\varphi\in\mathcal{M}$.

\begin{definition}\label{def3.2}
The space $H^{s,\varphi}(\mathbb{R}^{n})$ is the H\"ormander inner product space
$B_{2,\mu}(\mathbb{R}^{n})$ with $\mu(\xi):=\langle\xi\rangle^{s}\varphi(\xi)$ for
$\xi\in\mathbb{R}^{n}$.
\end{definition}

Thus $H^{s,\varphi}(\mathbb{R}^{n})$ consists of the distributions
$u\in\mathcal{S}'(\mathbb{R}^{n})$ such that the Fourier transform $\widehat{u}$ is a
function locally Lebesgue integrable on $\mathbb{R}^{n}$ and
$$
\int\limits_{\mathbb{R}^{n}}
\langle\xi\rangle^{2s}\varphi^{2}(\langle\xi\rangle)\,|\widehat{u}(\xi)|^{2}\,
d\xi<\infty.
$$
The inner product in the space $H^{s,\varphi}(\mathbb{R}^{n})$ is defined by the formula
$$
(u_{1},u_{2})_{H^{s,\varphi}(\mathbb{R}^{n})}:=
\int\limits_{\mathbb{R}^{n}}\langle\xi\rangle^{2s}\varphi^{2}(\langle\xi\rangle)
\,\widehat{u_{1}}(\xi)\,\overline{\widehat{u_{2}}(\xi)}\,d\xi
$$
and induces the norm in the usual way, $H^{s,\varphi}(\mathbb{R}^{n})$ being a
Hilbert space.

The function $\mu$ used in Definition \ref{def3.2} is a weight function that follows
from the integral representation of the set $\mathcal{M}$ given above. We consider
the Borel measurable weight functions $\mu$, rather than continuous as H\"ormander
does. By Theorem \ref{th2.4} i) we do not obtain the spaces different from those
considered by H\"ormander.

In the simplest case where $\varphi(\cdot)\equiv1$, the space
$H^{s,\varphi}(\mathbb{R}^{n})=H^{s,1}(\mathbb{R}^{n})$ coincides with the Sobolev space
$H^{s}(\mathbb{R}^{n})$.

By Theorem \ref{th2.4} (ii), for each $\varepsilon>0$ there exist a number
$c_{\varepsilon}\geq1$ such that
$$
c_{\varepsilon}^{-1}t^{-\varepsilon}\leq\varphi(t)\leq
c_{\varepsilon}t^{\varepsilon}\quad\mbox{for all}\quad t\geq1.
$$
This implies the inclusions
\begin{equation}\label{eq3.4}
\bigcup_{\varepsilon>0}H^{s+\varepsilon}(\mathbb{R}^{n})=:H^{s+}(\mathbb{R}^{n})
\subset H^{s,\varphi}(\mathbb{R}^{n})\subset
H^{s-}(\mathbb{R}^{n}):=\bigcap_{\varepsilon>0}H^{s-\varepsilon}(\mathbb{R}^{n}).
\end{equation}
They show that in the class of spaces
\begin{equation}\label{eq3.5}
\bigl\{H^{s,\varphi}(\mathbb{R}^{n}):\,s\in\mathbb{R},\,\varphi\in\mathcal{M}\,\bigr\}
\end{equation}
the functional parameter $\varphi$ defines a supplementary (subpower) smoothness to
the basic (power) $s$-smoothness. If $\varphi(t)\rightarrow\infty$
[$\varphi(t)\rightarrow0$] as $t\rightarrow\infty$, then $\varphi$ defines a
positive [negative] supplementary smoothness. Otherwise speaking, $\varphi$
\textit{refines} the power smoothness $s$. Therefore, it is naturally to give

\begin{definition}\label{def3.3}
The class of spaces \eqref{eq3.5} is called the refined Sobolev scale
over~$\mathbb{R}^{n}$.
\end{definition}

Obviously, the scale \eqref{eq3.5} is much finer than the Hilbert scale of Sobolev
spaces. The scale \eqref{eq3.5} was considered by the authors in \cite{05UMJ5,
06UMJ3, 08MFAT1}. Let us formulate some important properties of it.

\begin{theorem}\label{th3.3}
Let $s\in\mathbb{R}$ and $\varphi,\varphi_{1}\in\mathcal{M}$. The following
assertions are true:

\begin{itemize}
\item[i)] The dense continuous embedding
$H^{s+\varepsilon,\varphi_{1}}(\mathbb{R}^{n})\hookrightarrow
H^{s,\varphi}(\mathbb{R}^{n})$ is valid for each $\varepsilon>0$.
\item[ii)] The function $\varphi/\varphi_{1}$ is bounded in a neighbourhood of
$+\infty$ if and only if $H^{s,\varphi_{1}}(\mathbb{R}^{n})\hookrightarrow
H^{s,\varphi}(\mathbb{R}^{n})$. This embedding is continuous and dense.
\item[iii)] Let an integer $k\geq0$ be given. The inequality
\begin{equation}\label{eq3.6}
\int\limits_{1}^{\infty}\frac{dt}{t\,\varphi^{\,2}(t)}<\infty
\end{equation}
is equivalent to the embedding
\begin{equation}\label{eq3.7}
H^{k+n/2,\varphi}(\mathbb{R}^{n})\hookrightarrow C^{k}_{\mathrm{b}}(\mathbb{R}^{n}).
\end{equation}
The embedding is continuous.
\item[iv)] The spaces $H^{s,\varphi}(\mathbb{R}^{n})$ and
$H^{-s,1/\varphi}(\mathbb{R}^{n})$ are mutually dual with respect to the inner
product in $L_{2}(\mathbb{R}^{n})$.
\end{itemize}
\end{theorem}

Assertion i) of this theorem follows from \eqref{eq3.4}, whereas assertions ii) --
iv) are inherited from the H\"ormander spaces properties \cite[Sec.
2.2]{Hermander63}, in particular, iii) from Theorem \ref{th3.1}. Note that
$\varphi\in\mathcal{M}\Leftrightarrow1/\varphi\in\mathcal{M}$, so the space
$H^{-s,1/\varphi}(\mathbb{R}^{n})$ in assertion iv) is defined as an element of the
refined Sobolev scale.

The refined Sobolev scale possesses the interpolation property with respect to the
Sobolev scale because every space $H^{s,\varphi}(\mathbb{R}^{n})$ is obtained by the
interpolation, with an appropriate function parameter, of a couple of inner product
Sobolev spaces.

\begin{theorem}\label{th3.4}
Let a function $\varphi\in\mathcal{M}$ and positive numbers $\varepsilon,\delta$ be
given. We set
\begin{equation}\label{eq3.8}
\psi(t):=
\begin{cases}
\;t^{\,\varepsilon/(\varepsilon+\delta)}\,
\varphi(t^{1/(\varepsilon+\delta)}) & \text{for\;\;\;$t\geq1$}, \\
\;\varphi(1) & \text{for\;\;\;$0<t<1$}.
\end{cases}
\end{equation}
Then the following assertions are true:
\begin{itemize}
\item[i)] The function $\psi$ belongs to the set $\mathcal{B}$ and is an
interpolation parameter.
\item[ii)] For an arbitrary $s\in\mathbb{R}$, we have
\begin{equation}\label{eq3.9}
[H^{s-\varepsilon}(\mathbb{R}^{n}),H^{s+\delta}(\mathbb{R}^{n})]_{\psi}
=H^{s,\varphi}(\mathbb{R}^{n})
\end{equation}
\end{itemize}
with equality of norms in the spaces.
\end{theorem}

Assertion i) holds true by Corollary \ref{cor2.1} because the function \eqref{eq3.8}
is quasiregularly varying of index
$\theta:=\varepsilon/(\varepsilon+\delta)\in(0,\,1)$ at $+\infty$. Assertion ii) is
directly verified if we note that the operator
$J:u\mapsto\mathcal{F}^{-1}(\langle\xi\rangle^{\varepsilon+\delta}\,\widehat{u}(\xi))$
is generating for the couple on the left of \eqref{eq3.9}. Then the operator
$\psi(J):u\mapsto
\mathcal{F}^{-1}(\langle\xi\rangle^{\varepsilon}\varphi(\langle\xi\rangle)\,\widehat{u}(\xi))$
maps $H^{s,\varphi}(\mathbb{R}^{n})$ onto $H^{s-\varepsilon}(\mathbb{R}^{n})$ that
means \eqref{eq3.9}; for details, see \cite[Sec.~3]{06UMJ3} or \cite[Sec.
3.2]{08MFAT1}.

The refined Sobolev scale is closed with respect to the interpolation with the functions
parameters that are quasiregularly varying at $+\infty$.

\begin{theorem}\label{th3.5}
Let $s_{0},s_{1}\in\mathbb{R}$, $s_{0}\leq s_{1}$, and
$\varphi_{0},\varphi_{1}\in\mathcal{M}$. In the case where $s_{0}=s_{1}$ we suppose
that the function $\varphi_{0}/\varphi_{1}$ is bounded in a neighbourhood of
$\infty$. Let $\psi\in\mathcal{B}$ be a quasiregularly varying function of an index
$\theta\in(0,\,1)$ at $\infty$. We represent $\psi(t)=t^{\theta}\chi(t)$ with
$\chi\in\mathrm{QSV}$ and set $s:=(1-\theta)s_{0}+\theta s_{1}$,
$$
\varphi(t):=\varphi_{0}^{1-\theta}(t)\,\varphi_{1}^{\theta}(t)\,
\chi\Bigl(t^{s_{1}-s_{0}}\,\frac{\varphi_{1}(t)}{\varphi_{0}(t)}\Bigr)\quad\mbox{for}\quad
t\geq1.
$$
Then $\varphi\in\mathcal{M}$, and
\begin{equation}\label{eq3.10}
[H^{s_{0},\varphi_{0}}(\mathbb{R}^{n}),
H^{s_{1},\varphi_{1}}(\mathbb{R}^{n})]_{\psi}= H^{s,\varphi}(\mathbb{R}^{n})
\end{equation}
with equality of norms in the spaces.
\end{theorem}

This theorem can be proved by means of the repeated application of Theorem
\ref{th3.4} if we employ the reiteration formula $[X_{f},X_{g}]_{\psi}=X_{\omega}$,
where $X$ is an admissible couple of Hilbert spaces, $f,g,\psi\in\mathcal{B}$, $f/g$
is bounded in a neighbourhood of $\infty$, and $\omega(t):=f(t)\,\psi(g(t)/f(t))$
for $t>0$; see \cite[Sec. 2.3]{08MFAT1}. Besides, it is possible to give the direct
proof, which is similar to that used for Theorem~\ref{th3.4}.

\begin{remark}\label{rem3.2}
The interpolation of the H\"ormander spaces $B_{p,\mu}(\mathbb{R}^{n})$, with $1\leq
p\leq\infty$, was studied by M.~Schechter \cite{Schechter67} with the help of the
complex method of interpolation. C.~Merucci \cite{Merucci84} and F.~Cobos,
D.L.~Fernandez \cite{CobosFernandez88} considered the interpolation of various
Banach spaces of generalized smoothness by means of the real method involving a
function parameter.
\end{remark}

\section{Elliptic operators in $\mathbb{R}^{n}$}\label{sec4}

In this section we consider an arbitrary uniformly elliptic classical
pseudodifferential operator (PsDO) $A$ on the scale \eqref{eq3.5}. We establish an a
priory estimate for a solution to the equation $Au=f$ and investigate the solution
smoothness in this scale. Our results refine the classical theorems on elliptic
operators on the Sobolev scale; see, e.g., \cite[Sec. 1.8]{Agranovich94} or
\cite[Sec. 18.1]{Hermander85}.

Following \cite[Sec. 1.1]{Agranovich94}, we denote by $\Psi^{r}(\mathbb{R}^{n})$
with $r\in\mathbb{R}$ the class of all the PsDOs $A$ in $\mathbb{R}^{n}$ (generally,
not classical) such that their symbols $a(x,\xi)$ are complex-valued infinitely
smooth functions satisfying the following condition. For arbitrary multi-indexes
$\alpha$ and $\beta$, there exist a number $c_{\alpha,\beta}>0$ such that
$$
|\,\partial_{x}^{\alpha}\,\partial_{\xi}^{\beta}\,a(x,\xi)\,|
\leq\,c_{\alpha,\beta}\,\langle\xi\rangle^{r-|\beta|}\quad\mbox{for every}\quad
x,\xi\in\mathbb{R}^{n}.
$$

\begin{lemma}\label{lem4.1}
Let $A\in\Psi^{r}(\mathbb{R}^{n})$ with $r\in\mathbb{R}$. Then the restriction of
the mapping $u\mapsto Au$, $u\in\mathcal{S}'(\mathbb{R}^{n})$, to the space
$H^{s,\varphi}(\mathbb{R}^{n})$ is a bounded linear operator
$$
A:H^{s,\varphi}(\mathbb{R}^{n})\rightarrow H^{s-r,\,\varphi}(\mathbb{R}^{n})
$$
for each $s\in\mathbb{R}$ and $\varphi\in\mathcal{M}$.
\end{lemma}

This lemma follows from the Sobolev $\varphi\equiv1$ case \cite[Sec. 1.1, Theorem
1.1.2]{Agranovich94} by the interpolation formula \eqref{eq3.9}.

By $\Psi^{r}_{\mathrm{ph}}(\mathbb{R}^{n})$ we denote the subset in
$\Psi^{r}(\mathbb{R}^{n})$ that consists of all the classical (polyhomogeneous)
PsDOs of the order $r$; see \cite[Sec. 1.5]{Agranovich94}. An important example of
PsDO from $\Psi^{r}_{\mathrm{ph}}(\mathbb{R}^{n})$ is given by a partial
differential operator of order $r$ with coefficients belonging to
$C^{\infty}_{\mathrm{b}}(\mathbb{R}^{n})$.

\begin{definition}\label{def4.1}
A PsDO $A\in\Psi^{r}_{\mathrm{ph}}(\mathbb{R}^{n})$ is called uniformly elliptic in
$\mathbb{R}^{n}$ if there exists a number $c>0$ such that $|a_{0}(x,\xi)|\geq c$ for
each $x,\xi\in\mathbb{R}^{n}$ with $|\xi|=1$. Here $a_{0}(x,\xi)$ is the principal
symbol of $A$.
\end{definition}

Let $r\in\mathbb{R}$. Suppose a PsDO $A\in\Psi^{r}_{\mathrm{ph}}(\mathbb{R}^{n})$ to be
uniformly elliptic in $\mathbb{R}^{n}$.

\begin{theorem}\label{th4.1}
Let $s\in\mathbb{R}$, $\varphi\in\mathcal{M}$, and $\sigma<s$. The following
a~priori estimate holds true:
\begin{equation}\label{eq4.1}
\|u\|_{H^{s,\varphi}(\mathbb{R}^{n})}\leq
c\,\bigr(\,\|Au\|_{H^{s-r,\varphi}(\mathbb{R}^{n})}+
\|u\|_{H^{\sigma,\varphi}(\mathbb{R}^{n})}\,\bigl)\quad\mbox{for all}\quad u\in
H^{s,\varphi}(\mathbb{R}^{n}).
\end{equation}
Here $c=c(s,\varphi,\sigma)$ is a positive number not depending on $u$.
\end{theorem}

We prove this theorem with the help of the left parametrix of $A$ if we apply Lemma
\ref{lem4.1}. As knows \cite[Sec. 1.8, Theorem 1.8.3]{Agranovich94} there exists a
PsDO $B\in\Psi^{-r}_{\mathrm{ph}}(\mathbb{R}^{n})$ such that $BA=I+T$, where $I$ is
identical operator and $T\in\Psi^{-\infty}:=\bigcap_{m\in\mathbb{R}}\,
\Psi^{m}(\mathbb{R}^{n})$. The operator $B$ is called the left parametrix of $A$.
Writing $u=BAu-Tu$, we easily get \eqref{eq4.1} by Lemma \ref{lem4.1}.

Let $\Omega$ be an arbitrary nonempty open subset in $\mathbb{R}^{n}$. We study an
interior smoothness of a solution to the equation $Au=f$ in $\Omega$.

Let us introduce some relevant spaces. By $H^{-\infty}(\mathbb{R}^{n})$ we denote the
union of all the spaces $H^{s,\varphi}(\mathbb{R}^{n})$ with $s\in\mathbb{R}$ and
$\varphi\in\mathcal{M}$. The linear space $H^{-\infty}(\mathbb{R}^{n})$ is endowed with
the inductive limit topology. We set
\begin{gather}\notag
H^{s,\varphi}_{\mathrm{int}}(\Omega):=\bigl\{f\in H^{-\infty}(\mathbb{R}^{n}):
\,\chi\,f\in H^{s,\varphi}(\mathbb{R}^{n})\\
\mbox{for all}\;\;\chi\in
C^{\infty}_{\mathrm{b}}(\mathbb{R}^{n}),\;\mathrm{supp}\,\chi\subset \Omega,\;
\mathrm{dist}(\mathrm{supp}\,\chi,\partial \Omega)>0\bigr\}.\label{eq4.2}
\end{gather}
A topology in $H^{s,\varphi}_{\mathrm{int}}(\Omega)$ is defined by the seminorms
$f\mapsto\|\chi\,f\|_{H^{s,\varphi}(\mathbb{R}^{n})}$ with $\chi$ being the same as
in \eqref{eq4.2}.

\begin{theorem}\label{th4.2}
Let $u\in H^{-\infty}(\mathbb{R}^{n})$ be a solution to the equation $Au=f$ in
$\Omega$ with $f\in H^{s,\varphi}_{\mathrm{int}}(\Omega)$ for some $s\in\mathbb{R}$
and $\varphi\in\mathcal{M}$. Then $u\in H^{s+r,\varphi}_{\mathrm{int}}(\Omega)$.
\end{theorem}

The special case when $\Omega=\mathbb{R}^{n}$ (global smoothness) follows at once
from the equality $u=Bf-Tu$, with $B$ being the left parametrix, and Lemma
\ref{lem4.1}. In general, we deduce Theorem \ref{th4.2} from this case if we
rearrange $A$ and the operator of multiplication by a function $\chi$ satisfying
\eqref{eq4.2}. Then we write
\begin{equation}\label{eq4.3}
A\chi u=A\chi\eta u=\chi\,A\eta u+A'\eta u=\chi f+\chi\,A(\eta-1)u+A'\eta u,
\end{equation}
where $A'\in\Psi^{r-1}(\mathbb{R}^{n})$, and the function $\eta$ has the same
properties as $\chi$ and is equal to 1 in a neighbourhood of $\mathrm{supp}\,\chi$.
Now, if $u\in H^{s+r-k,\varphi}_{\mathrm{int}}(\Omega)$ for some integer $k\geq1$,
then the right-hand side of \eqref{eq4.3} belongs to
$H^{s-k+1,\varphi}(\mathbb{R}^{n})$ that implies $\chi u\in
H^{s+r-k+1,\varphi}(\mathbb{R}^{n})$, i.e. $u\in
H^{s+r-k+1,\varphi}_{\mathrm{int}}(\Omega)$. By induction in $k$ we have $u\in
H^{s+r,\varphi}_{\mathrm{int}}(\Omega)$.

It is useful to compare Theorem \ref{th4.2} with H\"ormander's Regularity Theorem.
If $A$ is a partial \emph{differential} operator, and $\Omega$ is bounded, then
Theorem \ref{th4.2} is a consequence of the H\"ormander theorem.

Applying Theorems \ref{th4.2} and \ref{th3.3} iii) we get the following sufficient
condition for the solution $u$ to have continuous and bounded derivatives of the
prescribed order.

\begin{theorem}\label{th4.3}
Let $u\in H^{-\infty}(\mathbb{R}^{n})$ be a solution to the equation $Au=f$ in
$\Omega$, with $f\in H^{k-r+n/2,\varphi}_{\mathrm{int}}(\Omega)$ for some integer
$k\geq0$ and function parameter $\varphi\in\mathcal{M}$. Suppose that $\varphi$
satisfies \eqref{eq3.6}. Then $u$ has the continuous partial derivatives on $\Omega$
up to the order $k$, and they are bounded on every set $\Omega_{0}\subset \Omega$
with $\mathrm{dist}(\Omega_{0},\partial \Omega)>0$. In particular, if
$\Omega=\mathbb{R}^{n}$, then $u\in C^{k}_{\mathrm{b}}(\mathbb{R}^{n})$.
\end{theorem}

This theorem shows an advantage of the refined Sobolev scale over the Sobolev scale
when a classical smoothness of a solution is under investigation. Indeed, if we
restrict ourselves to the Sobolev case of $\varphi\equiv1$, then we have to replace
the condition $f\in H^{k-r+n/2,\varphi}_{\mathrm{int}}(\Omega)$ with the condition
$f\in H^{k-r+\varepsilon+n/2,1}_{\mathrm{int}}(\Omega)$ for some $\varepsilon>0$.
The last condition is far stronger than previous one.

Note that the condition \eqref{eq3.6} not only is sufficient in Theorem 3.3 but also
is necessary on the class of all the considered solutions $u$. Namely, \eqref{eq3.6}
is equivalent to the implication
\begin{equation}\label{eq4.4}
\bigl(\,u\in H^{-\infty}(\mathbb{R}^{n}),\,\;\;f:=Au\in
H^{k-r+n/2,\varphi}_{\mathrm{int}}(\Omega)\,\bigr)\;\;\Rightarrow\;\;u\in
C^{k}(\Omega).
\end{equation}
Indeed, if $u\in H^{k+n/2,\varphi}_{\mathrm{int}}(\Omega)$, then $f=Au\in
H^{k-r+n/2,\varphi}_{\mathrm{int}}(\Omega)$, whence $u\in C^{k}(\Omega)$ if
\eqref{eq4.4} holds. Thus \eqref{eq4.4} entails \eqref{eq3.6} in view of
H\"ormander's Theorem \ref{th3.1}.

The analogs of Theorems \ref{th4.1}--\ref{th4.3} were proved in \cite{08UMB3}
for uniformly elliptic matrix PsDOs.

\section{H\"ormander spaces over a closed manifold}\label{sec5}

In this section we introduce a certain class of H\"ormander spaces over a closed
(compact) smooth manifold. Namely, using the spaces $H^{s,\varphi}(\mathbb{R}^{n})$
with $s\in\mathbb{R}$ and $\varphi\in\mathcal{M}$ we construct their analogs for the
manifold. We give three equivalent definitions of the analogs; these definitions are
similar to those used for the Sobolev spaces (see, e.g., \cite[Ch.~1,
Sec.~5]{Taylor81}).

\subsection{The equivalent definitions}\label{sec5.1}

In what follows except Subsection \ref{sec7.1}, $\Gamma$ is a closed (i.e. compact
and without a boundary) infinitely smooth oriented manifold of an arbitrary
dimension $n\geq1$. We suppose that a certain $C^{\infty}$-density $dx$ is defined
on $\Gamma$. As usual, $\mathcal{D}'(\Gamma)$ denotes the linear topological space
of all distributions on $\Gamma$. The space $\mathcal{D}'(\Gamma)$ is antidual to
the space $C^{\infty}(\Gamma)$ with respect to the natural extension of the scalar
product in $L_{2}(\Gamma):=L_{2}(\Gamma,dx)$ by continuity. This extension is
denoted by $(f,w)_{\Gamma}$ for $f\in\mathcal{D}'(\Gamma)$ and $w\in
C^{\infty}(\Gamma)$.

Let $s\in\mathbb{R}$ and $\varphi\in\mathcal{M}$. We give the following three equivalent
definitions of the H\"ormander space $H^{s,\varphi}(\Gamma)$.

The first definition exhibits the local properties of those distributions
$f\in\mathcal{D}'(\Gamma)$ that form $H^{s,\varphi}(\Gamma)$. From the
$C^{\infty}$-structure on $\Gamma$, we arbitrarily choose a finite collection of the
local charts $\alpha_{j}:\mathbb{R}^{n}\leftrightarrow\Gamma_{j}$,
$j=1,\ldots,\varkappa$, such that the open sets $\Gamma_{j}$ form the finite covering of
$\Gamma$. Let functions $\chi_{j}\in C^{\infty}(\Gamma)$, $j=1,\ldots,\varkappa$, form a
partition of unity on $\Gamma$ satisfying the condition
$\mathrm{supp}\,\chi_{j}\subset\Gamma_{j}$.

\begin{definition}\label{def5.1}
The linear space $H^{s,\varphi}(\Gamma)$ is defined by the formula
$$
H^{s,\varphi}(\Gamma):=\bigl\{f\in\mathcal{D}'(\Gamma):\;
(\chi_{j}f)\circ\alpha_{j}\in
H^{s,\varphi}(\mathbb{R}^{n})\;\;\forall\;j=1,\ldots\varkappa\bigr\}.
$$
Here $(\chi_{j}f)\circ\alpha_{j}$ is the representation of the distribution
$\chi_{j}f$ in the local chart $\alpha_{j}$. The inner product in the space
$H^{s,\varphi}(\Gamma)$ is introduced by the formula
$$
(f_{1},f_{2})_{H^{s,\varphi}(\Gamma)}:=\sum_{j=1}^{\varkappa}\,((\chi_{j}f_{1})\circ\alpha_{j},
(\chi_{j}\,f_{2})\circ\alpha_{j})_{H^{s,\varphi}(\mathbb{R}^{n})}
$$
and induces the norm in the usual way.
\end{definition}

In the special case where $\varphi\equiv1$ the space $H^{s,\varphi}(\Gamma)$
coincides with the inner product Sobolev space $H^{s}(\Gamma)$ of order $s$. The
Sobolev spaces on $\Gamma$ are known to be complete and independent (up to
equivalence of norms) of the choice of the local charts and the partition of unity.

The second definition connects the space $H^{s,\varphi}(\Gamma)$ with the Sobolev scale
by means of the interpolation.

\begin{definition}\label{def5.2}
Let two integers $k_{0}$ and $k_{1}$ be such that $k_{0}<s<k_{1}$. We define
\begin{equation}\label{eq5.1}
H^{s,\varphi}(\Gamma):=[H^{k_{0}}(\Gamma),H^{k_{1}}(\Gamma)]_{\psi},
\end{equation}
where the interpolation parameter $\psi$ is given by the formula \eqref{eq3.8} with
$\varepsilon:=s-k_{0}$ and $\delta:=k_{1}-s$.
\end{definition}

It is useful in the spectral theory to have the third definition of
$H^{s,\varphi}(\Gamma)$ that connects the norm in $H^{s,\varphi}(\Gamma)$ with a
certain function of $1-\Delta_{\Gamma}$. As usual, $\Delta_{\Gamma}$ is the
Beltrami-Laplace operator on the manifold $\Gamma$ endowed with the Riemannian
metric that induces the density $dx$; see, e.g., \cite[Sec. 22.1]{Shubin01}.

\begin{definition}\label{def5.3}
The space $H^{s,\varphi}(\Gamma)$ is defined to be the completion of
$C^{\infty}(\Gamma)$ with respect to the Hilbert norm
\begin{equation}\label{eq5.2}
f\,\mapsto\,
\|(1-\Delta_{\Gamma})^{s/2}\varphi((1-\Delta_{\Gamma})^{1/2})\,f\|_{L_{2}(\Gamma)},\quad
f\in C^{\infty}(\Gamma).
\end{equation}
\end{definition}

\begin{theorem}\label{th5.1}
Definitions $\ref{def5.1}$, $\ref{def5.2}$, and $\ref{def5.3}$ are mutually
equivalent, that is they define the same Hilbert space $H^{s,\varphi}(\Gamma)$ up to
equivalence of norms.
\end{theorem}

Let us explain how to prove this fundamental theorem.

The equivalence of Definitions \ref{def5.1} and \ref{def5.2}. We use Definition
\ref{def5.1} as a starting point and show that the equality \eqref{eq5.1} holds true
up to equivalence of norms. We apply the $\mathbb{R}^{n}$-analog of \eqref{eq5.1},
due to Theorem \ref{th3.4}, and pass to local coordinates on~$\Gamma$. Namely, let
the mapping $T$ take each $f\in\mathcal{D}'(\Gamma)$ to the vector with components
$(\chi_{j}f)\circ\alpha_{j}$, $j=1,\ldots,\varkappa$. We get the bounded linear
operator
\begin{equation}\label{eq5.3}
T:\,H^{s,\varphi}(\Gamma)\rightarrow(H^{s,\varphi}(\mathbb{R}^{n}))^{\varkappa}.
\end{equation}
It has the right inverse bounded linear operator
\begin{equation}\label{eq5.4}
K:\,(H^{s,\varphi}(\mathbb{R}^{n}))^{\varkappa}\rightarrow H^{s,\varphi}(\Gamma),
\end{equation}
where the mapping $K$ can be constructed with the help of the local charts and is
independent of parameters $s$ and $\varphi$. If we consider these operators in the
Sobolev case of $\varphi\equiv1$ and use the $\mathbb{R}^{n}$-analog of
\eqref{eq5.1}, then we get the bounded operators
\begin{gather}\label{eq5.5}
T:\,[H^{k_{0}}(\Gamma),H^{k_{1}}(\Gamma)]_{\psi}
\rightarrow(H^{s,\varphi}(\mathbb{R}^{n}))^{\varkappa}, \\
K:\,(H^{s,\varphi}(\mathbb{R}^{n}))^{\varkappa}
\rightarrow[H^{k_{0}}(\Gamma),H^{k_{1}}(\Gamma)]_{\psi}. \label{eq5.6}
\end{gather}
Now it follows from \eqref{eq5.3} and \eqref{eq5.6} that the identity mapping $KT$
establishes the continuous embedding of $H^{s,\varphi}(\Gamma)$ into the
interpolation space $[H^{k_{0}}(\Gamma),H^{k_{1}}(\Gamma)]_{\psi}$, whereas the
inverse continuous embedding is valid by \eqref{eq5.4} and \eqref{eq5.5}. Thus the
equality \eqref{eq5.1} holds true up to equivalence of norms; for details, see
\cite[Sec.~3]{06UMJ3} or \cite[Sec. 3.3]{08MFAT1}.

By equivalence of Definitions \ref{def5.1} and \ref{def5.2}, the space
$H^{s,\varphi}(\Gamma)$ is complete and independent (up to equivalence of norms) of
the choice of the local charts and the partition of unity on $\Gamma$. Moreover, the
set $C^{\infty}(\Gamma)$ is dense in $H^{s,\varphi}(\Gamma)$.

The equivalence of Definitions \ref{def5.2} and \ref{def5.3}. Let us use Definition
\ref{def5.2} as a starting point. If $s>0$, then we choose $k_{0}:=0$ and
$k_{1}:=2k>s$ for some integer $k$ in \eqref{eq5.1}. By $\Lambda_{k}(\Gamma)$ we
denote the Sobolev space $H^{2k}(\Gamma)$ endowed with the equivalent norm
$\|(1-\Delta_{\Gamma})^{k}f\|_{L_{2}(\Gamma)}$ of $f$. We have
\begin{eqnarray*}
\|f\|_{H^{s,\varphi}(\Gamma)}&=&\|f\|_{[H^{0}(\Gamma),H^{2k}(\Gamma)]_{\psi}}\asymp
\|f\|_{[L_{2}(\Gamma),\Lambda_{k}(\Gamma)]_{\psi}}=
\|\psi((1-\Delta_{\Gamma})^{k})f\|_{L_{2}(\Gamma)}\\
&=&\|(1-\Delta_{\Gamma})^{s/2}\varphi((1-\Delta_{\Gamma})^{1/2})\,f\|_{L_{2}(\Gamma)},
\end{eqnarray*}
with $f\in C^{\infty}(\Gamma)$, because $(1-\Delta_{\Gamma})^{k}$ is the generating
operator for the couple $[L_{2}(\Gamma),\Lambda_{k}(\Gamma)]$. Thus the norm
\eqref{eq5.2} is equivalent to the norm in $H^{s,\varphi}(\Gamma)$ on the dense set
$C^{\infty}(\Gamma)$ if $s>0$. The case of $s\leq0$ can be reduced to the previous
one with the help of the homeomorphism
$$
(1-\Delta_{\Gamma})^{k}:\,H^{s+2k,\varphi}(\Gamma)\leftrightarrow H^{s,\varphi}(\Gamma),
$$
with $s+2k>0$ for some integer $k\geq1$. This homeomorphism follows from the Sobolev
case of $\varphi\equiv1$ by the interpolation formula \eqref{eq5.1} (for details,
see \cite[Sec. 3.4]{08MFAT1}).

Now we can give

\begin{definition}\label{def5.4}
The class of Hilbert spaces
\begin{equation}\label{eq5.7}
\bigl\{H^{s,\varphi}(\Gamma):\,s\in\mathbb{R},\,\varphi\in\mathcal{M}\,\bigr\}
\end{equation}
is called the refined Sobolev scale over the manifold~$\Gamma$.
\end{definition}

\subsection{The properties}\label{sec5.2}

We consider some important properties of the scale \eqref{eq5.7}. They are inherited
from the refined Sobolev scale over $\mathbb{R}^{n}$.

\begin{theorem}\label{th5.2}
Let $s\in\mathbb{R}$ and $\varphi,\varphi_{1}\in\mathcal{M}$. The following
assertions are true:

\begin{enumerate}
\item[i)] The dense compact embedding
$H^{s+\varepsilon,\varphi_{1}}(\Gamma)\hookrightarrow H^{s,\varphi}(\Gamma)$ is
valid for each $\varepsilon>0$.
\item[ii)] The function $\varphi/\varphi_{1}$ is bounded in a neighbourhood of
$+\infty$ if and only if $H^{s,\varphi_{1}}(\Gamma)\hookrightarrow
H^{s,\varphi}(\Gamma)$. This embedding is continuous and dense. It is compact if and
only if $\varphi(t)/\varphi_{1}(t)\rightarrow0$ as $t\rightarrow+\infty$.
\item[iii)] Let integer $k\geq0$ be given. Inequality \eqref{eq3.6} is equivalent to the
embedding $H^{k+n/2,\varphi}(\Gamma)\hookrightarrow C^{k}(\Gamma)$. The embedding is
compact.
\item[iv)] The spaces $H^{s,\varphi}(\Gamma)$ and $H^{-s,1/\varphi}(\Gamma)$ are
mutually dual (up to equivalence of norms) with respect to the inner product in
$L_{2}(\mathbb{R}^{n})$.
\end{enumerate}
\end{theorem}

This theorem except the statements on compactness follows from Theorem \ref{th3.3}
in view of Definition \ref{def5.1}. The compactness of the embeddings is a
consequence of the compactness of~$\Gamma$. Namely, the statements in assertions i)
and ii) follow from the next proposition. For each number $r>0$, the embedding
$$
\{u\in H^{s,\varphi_{1}}(\mathbb{R}^{n}):\mathrm{dist}(0,\mathrm{supp}\,u)\leq
r\}\hookrightarrow H^{s,\varphi}(\mathbb{R}^{n})
$$
is compact if and only if $\varphi(t)/\varphi_{1}(t)\rightarrow0$ as
$t\rightarrow\infty$. This proposition is a special case of H\"ormander's result
\cite[Sec.~2, Theorem 2.2.3]{Hermander63}.  Now we get the compactness of the
embedding in assertion iii) if we write
$$
H^{k+n/2,\varphi}(\Gamma)\hookrightarrow H^{k+n/2,\varphi_{1}}(\Gamma)\hookrightarrow
C^{k}(\Gamma)
$$
for a function $\varphi_{1}$ such that the first embedding is compact.

\begin{theorem}\label{th5.3}
Theorems $\ref{th3.4}$ and $\ref{th3.5}$ remain true if we replace the designation
$\mathbb{R}^{n}$ with $\Gamma$, and the phrase `equality of norms' with `equivalence
of norms'.
\end{theorem}

Theorem \ref{th5.3} can be proved by means of a repeated application of the
interpolation formula \eqref{eq5.1}. We can also deduce this theorem from its
$\mathbb{R}^{n}$-analogs with the help of operators $T$ and $K$ used above.

\section{Elliptic operators on a closed manifold}\label{sec6}

Recall that $\Gamma$ is a closed infinitely smooth oriented manifold. In this
section we study an arbitrary elliptic classical PsDO $A$ on the refined Sobolev
scale over $\Gamma$. We prove that $A$ is a bounded and Fredholm operator on the
respective pairs of H\"ormander spaces, and investigate the smoothness of a solution
to the equation $Au=f$. Our results \cite[Sec. 4 and 5]{07UMJ6} refine the classical
theorems on elliptic operators on the Sobolev scale over a closed smooth manifold
(see, e.g., \cite[Sec.~2]{Agranovich94}, \cite[Sec.~19]{Hermander85}, and
\cite[\S~8]{Shubin01}). We also use some elliptic PsDOs to get an important class of
equivalent norms in $H^{s,\varphi}(\Gamma)$ \cite[Sec. 3.4]{08MFAT1}.

\subsection{The main properties}\label{sec6.1}

By $\Psi^{r}(\Gamma)$ with $r\in\mathbb{R}$ we denote the class of all the PsDOs $A$
on $\Gamma$ (generally, not classical) such that the image of $A$ in each local
chart on $\Gamma$ belongs to $\Psi^{r}(\mathbb{R}^{n})$; see \cite[Sec.
2.1]{Agranovich94}.

\begin{lemma}\label{lem6.1}
Let $A\in\Psi^{r}(\Gamma)$, with $r\in\mathbb{R}$. Then the restriction of the
mapping $u\mapsto Au$, $u\in\mathcal{D}'(\Gamma)$, to the space
$H^{s,\varphi}(\Gamma)$ is the bounded linear operator
\begin{equation}\label{eq6.1}
A:H^{s,\varphi}(\Gamma)\rightarrow H^{s-r,\varphi}(\Gamma)
\end{equation}
for each $s\in\mathbb{R}$ and $\varphi\in\mathcal{M}$.
\end{lemma}

This lemma follows from the Sobolev $\varphi\equiv1$ case \cite[Sec. 2.1, Theorem
2.1.2]{Agranovich94} by the interpolation in view of Theorem \ref{th5.3}.

By $\Psi^{r}_{\mathrm{ph}}(\Gamma)$ we denote the subset in $\Psi^{r}(\Gamma)$ that
consists of all classical PsDOs of the order $r$; see \cite[Sec. 2.1]{Agranovich94}.
The image of PsDO $A\in\Psi^{r}_{\mathrm{ph}}(\Gamma)$ in every local chart on
$\Gamma$ belongs to $\Psi^{r}_{\mathrm{ph}}(\mathbb{R}^{n})$.

\begin{definition}\label{def6.1}
A PsDO $A\in\Psi^{r}_{\mathrm{ph}}(\mathbb{R}^{n})$ is called elliptic on $\Gamma$
if $a_{0}(x,\xi)\neq0$ for each point $x\in\Gamma$ and covector $\xi\in
T^{\ast}_{x}\Gamma\setminus\{0\}$. Here $a_{0}(x,\xi)$ is the principal symbol of
$A$, and $T^{\ast}_{x}\Gamma$ is the cotangent space to $\Gamma$ at $x$.
\end{definition}

Let $r\in\mathbb{R}$. Suppose a PsDO $A\in\Psi^{r}_{\mathrm{ph}}(\Gamma)$ to be
elliptic on $\Gamma$.

By $A^{+}$ we denote the PsDO formally adjoint to $A$ with respect to the
sesquilinear form $(\cdot,\cdot)_{\Gamma}$. Since both $A$ and $A^{+}$ are elliptic
on $\Gamma$, both the spaces
\begin{gather*}
N:=\{\,u\in C^{\infty}(\Gamma):\,Au=0\;\;\mbox{on}\;\;\Gamma\,\},\\
N^{+}:=\{\,v\in C^{\infty}(\Gamma):\,A^{+}v=0\;\;\mbox{on}\;\;\Gamma\,\}
\end{gather*}
are finite-dimensional \cite[Sec. 8.2, Theorem 8.1]{Shubin01}.

Recall the following definition.

\begin{definition}\label{def6.2}
 Let $X$ and $Y$ be Banach spaces. The bounded linear
operator $T:X\rightarrow Y$ is said to be Fredholm if its kernel $\ker T$ and
co-kernel $\mathrm{coker}\,T:=Y/\,T(X)$ are finite-dimensional. The number
$\mathrm{ind}\,T:=\dim\ker T-\dim\mathrm{coker}\,T$ is called the index of the
Fredholm operator~$T$.
\end{definition}

If the operator $T:X\rightarrow Y$ is Fredholm, then its range $T(X)$ is closed in
$Y$; see, e.g., \cite[Sec. 19.1, Lemma 19.1.1]{Hermander85}.

\begin{theorem}\label{th6.1}
Let $s\in\mathbb{R}$ and $\varphi\in\mathcal{M}$. For the elliptic PsDO $A$, the
operator \eqref{eq6.1} is Fredholm, has the kernel $N$ and the range
\begin{equation}\label{eq6.2}
A(H^{s,\varphi}(\Gamma))=\bigl\{f\in
H^{s-r,\varphi}(\Gamma):\,(f,v)_{\Gamma}=0\;\;\mbox{for all}\;\;v\in N^{+}\bigr\}.
\end{equation}
The index of the operator \eqref{eq6.1} is equal to $\dim N-\dim N^{+}$ and
independent of $s$ and $\varphi$.
\end{theorem}

Theorem \ref{th6.1} is well known in the Sobolev case of $\varphi\equiv1$; see,
e.g., \cite[Sec. 19.2, Theorem 19.2.1]{Hermander85}. For an arbitrary
$\varphi\in\mathcal{M}$, we deduce this theorem  from the Sobolev case with the help
of the interpolation. Indeed, consider the bounded Fredholm operators
$A:H^{s\mp1}(\Gamma)\rightarrow H^{s\mp1-r}(\Gamma)$. By applying Theorem
\ref{th5.3}, we have the bounded operator
$$
A:\,H^{s,\varphi}(\Gamma)=[H^{s-1}(\Gamma),H^{s+1}(\Gamma)]_{\psi}\rightarrow
[H^{s-1-r}(\Gamma),H^{s+1-r}(\Gamma)]_{\psi}=H^{s-r,\varphi}(\Gamma).
$$
Here $\psi$ is the interpolation parameter defined by the formula \eqref{eq3.8} with
$\varepsilon=\delta=\nobreak1$. This operator has the properties stated in Theorem
\ref{th6.1} because of the next proposition.

\begin{proposition}\label{prop6.1}
Let $X=[X_{0},X_{1}]$ and $Y=[Y_{0},Y_{1}]$ be admissible couples of Hilbert spaces,
and a linear mapping $T$ be given on $X_{0}$. Suppose we have the Fredholm bounded
operators $T:X_{j}\rightarrow Y_{j}$, with $j=0,\,1$, that possess the common kernel
and the common index. Then, for an arbitrary interpolation parameter
$\psi\in\mathcal{B}$, the bounded operator $T:X_{\psi}\rightarrow Y_{\psi}$ is
Fredholm, has the same kernel and the same index, moreover its range
$T(X_{\psi})=Y_{\psi}\cap T(X_{\,0})$.
\end{proposition}

This proposition was proved by G.~Geymonat \cite[p.~280, Proposition
5.2]{Geymonat65} for arbitrary interpolation functors given on the category of
couples of Banach spaces. The proof for Hilbert spaces is analogous.

If both the spaces $N$ and $N^{+}$ are trivial, then the operator \eqref{eq6.1} is a
homeomorphism. Generally, the index of \eqref{eq6.1} is equal to zero provided that
$\dim\Gamma\geq2$ (see \cite{AtiyahSinger63} and \cite[Sec. 2.3~f]{Agranovich94}).
In the case where $\dim\Gamma=1$, the index can be nonzero. If $A$ is a differential
operator, then the index is always zero.

The Fredholm property of $A$ implies the following a priory estimate.

\begin{theorem}\label{th6.2}
Let $s\in\mathbb{R}$, $\varphi\in\mathcal{M}$, and $\sigma<s$. Then
$$
\|u\|_{H^{s,\varphi}(\Gamma)}\leq c\,\bigr(\,\|Au\|_{H^{s-r,\varphi}(\Gamma)}+
\|u\|_{H^{\sigma,\varphi}(\Gamma)}\,\bigl)\quad\mbox{for all}\quad u\in
H^{s,\varphi}(\Gamma);
$$
here the number $c>0$ is independent of $u$.
\end{theorem}

We obtain the above implication if we use the compactness of the embedding
$H^{s,\varphi}(\Gamma)\hookrightarrow H^{\sigma,\varphi}(\Gamma)$ for $\sigma<s$ and
apply the following proposition \cite[Sec. 2.3, Theorem 2.3.4]{Agranovich94}.

\begin{proposition}\label{prop6.2}
Let $X$, $Y$, and $Z$ be Banach spaces. Suppose that the compact embedding
$X\hookrightarrow Z$ is valid, and a bounded linear operator $T:X\rightarrow Y$ is
given. Then $\ker T$ is finite-dimensional and $T(X)$ is closed in $Y$ if and only
if there exists a number $c>0$ such that
$$
\|u\|_{X}\leq c\,(\,\|Tu\|_{Y}+\|u\|_{Z}\,)\quad\mbox{for all}\quad u\in X.
$$
\end{proposition}

Now we study a local smoothness of a solution to the elliptic equation $Au=f$. Let
$\Gamma_{0}$ be an nonempty open set on the manifold $\Gamma$, and define
\begin{equation}\label{eq6.3}
H^{s,\varphi}_{\mathrm{loc}}(\Gamma_{0})
:=\bigl\{f\in\mathcal{D}'(\Gamma):\,\chi\,f\in
H^{s,\varphi}(\Gamma),\;\;\forall\;\chi\in
C^{\infty}(\Gamma),\;\mathrm{supp}\,\chi\subseteq \Gamma_{0}\bigr\}.
\end{equation}

\begin{theorem}\label{th6.3}
Let $u\in\mathcal{D}'(\Gamma)$ be a solution to the equation $Au=f$ on $\Gamma_{0}$
with $f\in H^{s,\varphi}_{\mathrm{loc}}(\Gamma_{0})$ for some $s\in\mathbb{R}$ and
$\varphi\in\mathcal{M}$. Then $u\in H^{s+r,\varphi}_{\mathrm{loc}}(\Gamma_{0})$.
\end{theorem}

The special case when $\Gamma_{0}=\Gamma$ (global smoothness) follows from Theorem
\ref{th6.1}. Indeed, using \eqref{eq6.2} we can write $f=Av$ for some $v\in
H^{s+r,\varphi}(\Gamma)$, whence $u-v\in N$ and $u\in H^{s+r,\varphi}(\Gamma)$. In
general, we deduce Theorem \ref{th6.3} from this case reasoning similar to the proof
of Theorem~\ref{th4.2}.

If we bring Theorems \ref{th6.3} and \ref{th5.2} iii) together, then we get the
following sufficient condition for the solution $u$ to have continuous derivatives
of the prescribed order on $\Gamma_{0}$. Recall that $n:=\dim\Gamma$.

\begin{theorem}\label{th6.4}
Let $u\in\mathcal{D}'(\Gamma)$ be a solution to the equation $Au=f$ on $\Gamma_{0}$,
with $f\in H^{k-r+n/2,\varphi}_{\mathrm{loc}}(\Gamma_{0})$ for some integer $k\geq0$
and function parameter $\varphi\in\mathcal{M}$. If $\varphi$ satisfies
\eqref{eq3.6}, then $u\in C^{k}(\Gamma_{0})$.
\end{theorem}

Here it is important that the condition \eqref{eq3.6} not only is sufficient for $u$
to belong to $C^{k}(\Gamma_{0})$ but also is necessary on the class of all the
considered solutions $u$.

\subsection{The equivalent norms induced by elliptic operators}\label{sec6.2}

Let $r>0$, and a PsDO $A\in\Psi^{r}_{\mathrm{ph}}(\Gamma)$ be elliptic on $\Gamma$.
We may consider $A$ as a closed unbounded operator on $L_{2}(\Gamma)$ with the
domain $H^{r}(\Gamma)$; see, e.g., \cite[Sec. 2.3, Theorem 2.3.5]{Agranovich94}.
Suppose the operator $A$ to be positive in $L_{2}(\Gamma)$. Then $A$ is self-adjoint
in $L_{2}(\Gamma)$ \cite[Sec. 2.3, Theorem 2.3.7]{Agranovich94}.

For $s\in\mathbb{R}$ and $\varphi\in\mathcal{M}$, we set
$$
\varphi_{s,r}(t):=t^{s/r}\varphi(t^{1/r})\;\;\mbox{for}\;\;t\geq1,\quad\mbox{and}
\quad\varphi_{s,r}(t):=\varphi(1)\;\;\mbox{for}\;\;0<t<1.
$$
The operator $\varphi_{s,r}(A)$ is regarded as the Borel function $\varphi_{s,r}$ of the
positive self-adjoint operator~$A$ in $L_{2}(\Gamma)$. Consider the norm
\begin{equation}\label{eq6.4}
f\,\mapsto\,\|\varphi_{s,r}(A)f\|_{L_{2}(\Gamma)},\quad f\in C^{\infty}(\Gamma).
\end{equation}

\begin{theorem}\label{th6.5}
Let $s\in\mathbb{R}$ and $\varphi\in\mathcal{M}$. The norm in the space
$H^{s,\varphi}(\Gamma)$ is equivalent to the norm \eqref{eq6.4} on the set
$C^{\infty}(\Gamma)$. Thus $H^{s,\varphi}(\Gamma)$ is the completion of
$C^{\infty}(\Gamma)$ with respect to the norm \eqref{eq6.4}.
\end{theorem}

The proof of this theorem is quite similar to the reasoning we did to demonstrate
the equivalence of Definitions \ref{def5.2} and \ref{def5.3}.

If $H^{s,\varphi}(\Gamma)\hookrightarrow L_{2}(\Gamma)$, then Theorem \ref{th6.5}
entails the following.

\begin{corollary}\label{cor6.1}
Let $s\geq0$ and $\varphi\in\mathcal{M}$. In the case where $s=0$ we suppose that
the function $1/\varphi$ is bounded in a neighbourhood of $\infty$. Then the space
$H^{s,\varphi}(\Gamma)$ coincides with the domain of the operator
$\varphi_{s,r}(A)$, and the norm in the space $H^{s,\varphi}(\Gamma)$ is equivalent
to the graphics norm of $\varphi_{s,r}(A)$.
\end{corollary}

It is useful to have an analog of Theorem \ref{th6.5} formulated in terms of
sequences. For this purpose, we recall some spectral properties of the operator $A$;
see, e.g., \cite[Sec. 6.1]{Agranovich94} or \cite[Sec. 8.3 and 15.2]{Shubin01}.

There is an orthonormal basis $(h_{j})_{j=1}^{\infty}$ of $L_{2}(\Gamma)$ formed by
eigenfunctions $h_{j}\in C^{\infty}(\Gamma)$ of the operator~$A$. Let
$\lambda_{j}>0$ is the eigenvalue corresponding to $h_{j}$; the enumeration is such
that $\lambda_{j}\leq\lambda_{j+1}$. Then the spectrum of $A$ coincides with the set
$\{\lambda_{1},\lambda_{2},\lambda_{3},\ldots\}$ of all eigenvalues of $A$, and the
asymptotics formula holds: $\lambda_{j}\sim c\,j^{\,r/n}$ as $j\rightarrow\infty$.
Each distribution $f\in\mathcal{D}'(\Gamma)$ is expanded into the Fourier series
\begin{equation}\label{eq6.5}
f=\sum_{j=1}^{\infty}\;c_{j}(f)\,h_{j}
\end{equation}
convergent in $\mathcal{D}'(\Gamma)$; here $c_{j}(f):=(f,h_{j})_{\Gamma}$ are the
Fourier coefficients of $f$.

\begin{theorem}\label{th6.6}
Let $s\in\mathbb{R}$ and $\varphi\in\mathcal{M}$. Then the next equality of spaces
with equivalence of norms in them holds:
\begin{gather}\label{eq6.6}
H^{s,\varphi}(\Gamma)=\Bigl\{f\in\mathcal{D}'(\Gamma):\,
\sum_{j=1}^{\infty}\;j^{\,2s/n}\,\varphi^{2}(j^{\,1/n})\,|c_{j}(f)|^{2}<\infty\Bigr\},\\
\|f\|_{H^{s,\varphi}(\Gamma)}\asymp\Bigl(\;\sum_{j=1}^{\infty}\;
j^{\,2s/n}\,\varphi^{2}(j^{\,1/n})\,|c_{j}(f)|^{2}\;\Bigr)^{1/2}.\label{eq6.7}
\end{gather}
\end{theorem}

This theorem follows from Theorem \ref{th6.5} since
$$
\varphi_{s,r}(A)\,f=\sum_{j=1}^{\infty}\,\varphi_{s,r}(\lambda_{j})\,c_{j}(f)\,h_{j}
\quad\mbox{(convergence in $L_{2}(\Gamma)$)}
$$
for each function $f$ from the domain of the operator $\varphi_{s,r}(A)$. Here
$\varphi_{s,r}(\lambda_{j})\asymp j^{\,s/n}\varphi(j^{\,1/n})$ with integers
$j\geq1$ in view of the asymptotics formula mentioned above. By applying Parseval's
equality we get \eqref{eq6.7} and then can deduce \eqref{eq6.6}.

\begin{corollary}\label{cor6.2}
Suppose $f\in H^{s,\varphi}(\Gamma)$, then the series \eqref{eq6.5} converges to $f$
in the space $H^{s,\varphi}(\Gamma)$.
\end{corollary}

This is a simple consequence of \eqref{eq6.7}.

\begin{example}\label{ex6.1}
Let $\Gamma$ be a unit circle and $A:=1-d^{2}/dt^{2}$, where $t$ sets the natural
parametrization on $\Gamma$. The eigenfunctions
$h_{j}(t):=(2\pi)^{-1}\mathrm{e}^{ijt}$, $j\in\mathbb{Z}$, of $A$ form an
orthonormal basis in $L_{2}(\Gamma)$, with $\lambda_{j}:=1+j^{\,2}$ being the
corresponding eigenvalues. Therefore the equivalence \eqref{eq6.7} becomes
$$
\|f\|_{H^{s,\varphi}(\Gamma)}\asymp\|\varphi_{s,2}(A)f\|_{L_{2}(\Gamma)}=
\Bigl(\;\sum_{j=-\infty}^{\infty}\;(1+j^{\,2})^{s}
\,\varphi^{2}((1+j^{\,2})^{1/2})\,|c_{j}(f)|^{2}\;\Bigr)^{1/2}.
$$
Note that we can chose the basis formed by the real-valued eigenfunctions
$h_{0}(t):=(2\pi)^{-1}$, $h_{j}(t):=\pi^{-1}\cos jt$, and $h_{-j}(t):=\pi^{-1}\sin
jt$, with integral $j\geq\nobreak1$. Then
$$
\|f\|_{H^{s,\varphi}(\Gamma)}^{2}\asymp|a_{0}(f)|^{2}+
\sum_{j=1}^{\infty}\;j^{\,2s}\,\varphi^{2}(j)\,\bigl(|a_{j}(f)|^{2}+|b_{j}(f)|^{2}\bigr),
$$
where $a_{0}(f)$, $a_{j}(f)$, and $b_{j}(f)$ are the Fourier coefficients of $f$
with respect to these eigenfunctions. In the considered case,
$H^{s,\varphi}(\Gamma)$ is closely related to the spaces of periodic functions
considered by A.I.~Stepanets \cite[Ch.~I, \S~7]{Stepanets87}, \cite[Part~I, Ch.~3,
Sec. 7.1]{Stepanets05}.
\end{example}

\section{Applications to spectral expansions}\label{sec7}

In this section, we investigate the convergence of expansions in eigenfunctions of
normal (in particular, self-adjoint) elliptic PsDOs given on a compact smooth
manifold. Using the refined Sobolev scale, we find new sufficient conditions for the
convergence almost everywhere; they are expressed in constructive terms of
regularity of functions. We also give a criterion for convergence in the metrics of
$C^{k}$ on the classes being H\"ormander spaces. Beforehand let us recall some
classical results concerning the convergence almost everywhere of arbitrary
orthogonal series. The results will be used below.

\subsection{The classical results}\label{sec7.1}

In this subsection, $\Gamma$ is an arbitrary set with a finite measure~$\mu$.
Suppose that $(h_{j})_{j=1}^{\infty}$ is an orthonormal system of functions in
$L_{2}(\Gamma):=L_{2}(\Gamma,d\mu)$, generally, complex-valued. The following
proposition is a general version of the well-known Menshov-Rademacher convergence
theorem.

\begin{theorem}\label{th7.1}
Let a sequence of complex numbers $(a_{j})_{j=1}^{\infty}$ be such that
\begin{equation}\label{eq7.1}
\sum_{j=2}^{\infty}\;|a_{j}|^{2}\,\log^{2}j<\infty.
\end{equation}
Then the series
\begin{equation}\label{eq7.2}
\sum_{j=1}^{\infty}\;a_{j}\,h_{j}(x)
\end{equation}
converges almost everywhere on $\Gamma$.
\end{theorem}

This theorem was proved independently by D.E.~Menshov \cite{Menschoff23} and
H.~Rad\-e\-ma\-cher \cite{Rademacher22} in the classical case where $\Gamma$ is a
bounded interval on $\mathbb{R}$, $\mu$ is the Lebesgue measure, and the functions
$h_{j}$ are real-valued. The proof of the Menshov--Rademacher theorem given in
\cite[Ch.~8, \S~1]{KashinSaakyan89} remains valid in the general situation that we
consider (apparently, the most general case is treated in \cite{11MFAT4}).

It is important that the Menshov--Rademacher theorem is precise. Menshov
\cite{Menschoff23} gave an example of an orthonormal system $(h_{j})_{j=1}^{\infty}$
in $L_{2}((0,1))$ for which Theorem \ref{th7.1} will not be true if one replaces, in
\eqref{eq7.1}, the sequence $(\log^{2}j)_{j=1}^{\infty}$ by any increasing sequence
of positive numbers $\omega_{j}=o(\log^{2}j)$ with $j\rightarrow\infty$. This result
is set forth, e.g., in the monograph \cite[Ch.~8, \S~1, Theorem~2]{KashinSaakyan89}.

Note that, for series \eqref{eq7.2} with coefficients subject to \eqref{eq7.1}, the
convergence almost everywhere need not be unconditional; see, e.g., \cite[Ch.~8,
\S~2]{KashinSaakyan89}. Recall that a series of functions is said to be
unconditionally convergent almost everywhere on a set if it remains convergent
almost everywhere on the set after arbitrary permutation of its terms (the null
measure set of divergence may vary.) The following proposition is a general version
of the Orlicz theorem on unconditional convergence of orthogonal series of
functions.

\begin{theorem}\label{th7.2}
Let a sequence of complex numbers $(a_{j})_{j=1}^{\infty}$ and increasing sequence
of positive numbers $(\omega_{j})_{j=1}^{\infty}$ satisfy the conditions
$$
\sum_{j=2}^{\infty}\;|a_{j}|^{2}\,(\log^{2}j)\,\omega_{j}<\infty,\quad
\sum_{j=2}^{\infty}\;\frac{1}{j\,(\log j)\,\omega_{j}}<\infty.
$$
Then the series \eqref{eq7.2} converges unconditionally almost everywhere on
$\Gamma$.
\end{theorem}

This equivalent statement of W.~Orlicz' theorem \cite{Orlicz27} was given by
P.L.~Uly\-a\-nov \cite[\S~4]{Uljanov63}; they considered the classical case
mentioned above. In our (more general) case, Theorem 6.2 follows from K.~Tandori's
theorem \cite{Tandori62}, which remains valid for arbitrary measure space
\cite{11UMJ11, 11MFAT4}. As Tandori proved \cite{Tandori62}, the Orlicz theorem is
the best possible in the sense that its condition on the sequence
$(\omega_{j})_{j=1}^{\infty}$ cannot be weaken.

\subsection{The convergence of spectral expansions}\label{sec7.2}

Further, $\Gamma$ is a closed infinitely smooth oriented manifold, and
$n=\dim\Gamma$. A $C^{\infty}$-density $dx$ is given on $\Gamma$ and defines the
finite measure there. Let a PsDO $A\in\Psi^{r}_{\mathrm{ph}}(\Gamma)$, with $r>0$,
be elliptic on~$\Gamma$. Suppose that $A$ is a normal (unbounded) operator on
$L_{2}(\Gamma)=L_{2}(\Gamma,dx)$. Let $(h_{j})_{j=1}^{\infty}$ be a complete
orthonormal system of eigenfunctions of this operator. They are enumerated so that
$|\lambda_{j}|\leq|\lambda_{j+1}|$ for $j=1,2,3,\ldots$, where
$Ah_{j}=\lambda_{j}h_{j}$. For an arbitrary function $f\in L_{2}(\Gamma)$, we
consider its expansion into the Fourier series \eqref{eq6.5} with respect to the
system $(h_{j})_{j=1}^{\infty}$.

We say that the series \eqref{eq6.5} converges on a function class $X(\Gamma)$ in
the indicated sense if, for every function $f\in X(\Gamma)$, the series converges to
$f$ in the indicated manner.

We investigate the convergence almost everywhere of the spectral expansion
\eqref{eq6.5} with the help of Theorems \ref{th6.6}, \ref{th7.1}, and \ref{th7.2}.
Let $\log^{\ast}t:=\max\{1,\log t\}$ for $t\geq1$.

\begin{theorem}\label{th7.3}
The series \eqref{eq6.5} converges almost everywhere on $\Gamma$ on the class
$H^{0,\log^{\ast}}(\Gamma)$.
\end{theorem}

Indeed, if $A$ is a positive operator on $L_{2}(\Gamma)$, then by Theorem
\ref{th6.6} we have
$$
|c_{1}(f)|^{2}+ \sum_{j=2}^{\infty}\;|c_{j}(f)|^{2}\,\log^{2}j\;\asymp\;
\|f\|_{H^{0,\log^{\ast}}(\Gamma)}^{2}<\infty\quad\mbox{for}\quad f\in
H^{0,\log^{\ast}}(\Gamma).
$$
This and Theorem \ref{th7.1} yields Theorem \ref{th7.3}. In general, if $A$ is a
normal operator, we should exchange $A$ for the positive elliptic PsDO
$B:=1+A^{\ast}A$ in our consideration and use the fact that $(h_{j})_{j=1}^{\infty}$
is a complete system of eigenfunctions of $B$.

Similarly, if we bring together Theorems \ref{th6.6} and \ref{th7.2}, we will get
the following result.

\begin{theorem}\label{th7.4}
Let an increasing function $\varphi\in\mathcal{M}$ be such that
$$
\int\limits_{2}^{\infty}\frac{dt}{t\,(\log t)\,\varphi^{2}(t)}<\infty.
$$
Then the series \eqref{eq6.5} converges unconditionally almost everywhere on
$\Gamma$ on the class $H^{0,\varphi\log^{\ast}}(\Gamma)$.
\end{theorem}

Note that the applying of H\"ormander spaces permits us to use the conditions of
Theorems \ref{th7.1} and \ref{th7.2} in an exhaustive manner. If we remain in the
framework of the Sobolev scale, we deduce that the series \eqref{eq6.5} converges
(unconditionally) almost everywhere on $\Gamma$ on the class $H^{0+}(\Gamma):=
\bigcup_{\varepsilon>0}H^{\varepsilon}(\Gamma)$. The result is far rougher than
those formulated in Theorems \ref{th7.3} and \ref{th7.4}. In the special case of
$A=\Delta_{\Gamma}$, this result was proved by C.~Meaney \cite{Meaney82}. (As Meaney
noted, it has "all the qualities of a folk theorem".)

To compete this section we give a criterion for the convergence of the spectral
expansions in the metrics of $C^{k}(\Gamma)$ on the classes $H^{s,\varphi}(\Gamma)$.

\begin{theorem}\label{th7.5}
Let an integer $k\geq0$ and function $\varphi\in\mathcal{M}$ be given. The series
\eqref{eq6.5} converges in $C^{k}(\Gamma)$ on the class $H^{k+n/2,\varphi}(\Gamma)$
if and only if $\varphi$ satisfies \eqref{eq3.6}.
\end{theorem}

This theorem results from Corollary \ref{cor6.2} and Theorem \ref{th5.2} iii).

Note that the convergence conditions in Theorems \ref{th7.3}, \ref{th7.4}, and
\ref{th7.5} are given in constructive terms of regularity of functions. The
regularity properties can be determined locally on $\Gamma$ according to Definition
\ref{def5.1}.

\section{H\"ormander spaces over Euclidean domains}\label{sec8}

In the next sections, we will consider some applications of H\"ormander spaces to
elliptic boundary problems in a bounded domain $\Omega\subset\mathbb{R}^{n}$. For
this purpose we need the H\"ormander spaces that consists of distributions given
in~$\Omega$. The spaces of distributions supported on the closure
$\overline{\Omega}$ of the domain $\Omega$ is also of use. These spaces are
constructed from the H\"ormander spaces over $\mathbb{R}^{n}$ in the standard way
\cite[Ch.~1, \S~3]{VolevichPaneah65}. We are interested in the H\"ormander spaces
that form the refined Sobolev scales over $\Omega$ and $\overline{\Omega}$. In this
section, we give the definitions of these spaces and consider their properties,
among them the interpolation properties being of great importance for us. We also
study a connection between the refined Sobolev scales over $\Omega$ and its boundary
(the trace theorems) and introduce riggings of $L_{2}(\Omega)$ with some H\"ormander
spaces.

\subsection{The definitions}\label{sec8.1}

Let $s\in\mathbb{R}$ and $\varphi\in\mathcal{M}$.

\begin{definition}\label{def8.1}
Suppose that $Q$ is a nonempty closed set in $\mathbb{R}^{n}$. The linear space
$H^{s,\varphi}_{Q}(\mathbb{R}^{n})$ is defined to consist of the distributions $u\in
H^{s,\varphi}(\mathbb{R}^{n})$ such that $\mathrm{supp}\,u\subseteq Q$. The space
$H^{s,\varphi}_{Q}(\mathbb{R}^{n})$ is endowed with the inner product and norm from
$H^{s,\varphi}(\mathbb{R}^{n})$.
\end{definition}

The space $H^{s,\varphi}_{Q}(\mathbb{R}^{n})$ is complete (i.e., Hilbert) because of
the continuous embedding of the Hilbert space $H^{s,\varphi}(\mathbb{R}^{n})$ into
$\mathcal{S}'(\mathbb{R}^{n})$.

\begin{definition}\label{def8.2}
Suppose that $\Omega$ is a nonempty open set in $\mathbb{R}^{n}$. The linear space
$H^{s,\varphi}(\Omega)$ is defined to consist of the restrictions
$v=u\!\upharpoonright\!\Omega$ of all the distributions $u\in
H^{s,\varphi}(\mathbb{R}^{n})$ to $\Omega$. The space $H^{s,\varphi}(\Omega)$ is
endowed with the norm
\begin{equation}\label{eq8.1}
\|v\|_{H^{s,\varphi}(\Omega)}:=
\inf\,\bigl\{\,\|u\|_{H^{s,\varphi}(\mathbb{R}^{n})}:\,u\in
H^{s,\varphi}(\mathbb{R}^{n}),\;\;v=u\;\;\mbox{in}\;\;\Omega\,\bigr\}.
\end{equation}
\end{definition}

By Definition \ref{def8.2}, $H^{s,\varphi}(\Omega)$ is a factor space
$H^{s,\varphi}(\mathbb{R}^{n})/H^{s,\varphi}_{\widehat{\Omega}}(\mathbb{R}^{n})$,
where $\widehat{\Omega}:=\mathbb{R}^{n}\setminus\Omega$. Hence, the space
$H^{s,\varphi}(\Omega)$ is Hilbert; the norm \ref{def8.1} is induced by the inner
product
$$
\bigl(v_{1},v_{2}\bigr)_{H^{s,\varphi}(\Omega)}:= \bigl(u_{1}-\Pi u_{1},u_{2}-\Pi
u_{2}\bigr)_{H^{s,\varphi}(\mathbb{R}^{n})}.
$$
Here $u_{j}\in H^{s,\varphi}(\mathbb{R}^{n})$, $u_{j}=v_{j}$ in $\Omega$ for $j=1,\,2$,
and $\Pi$ is the orthogonal projector of the space $H^{s,\varphi}(\mathbb{R}^{n})$ onto
the subspace $H^{s,\varphi}_{\widehat{\Omega}}(\mathbb{R}^{n})$.

Both the spaces $H^{s,\varphi}_{Q}(\mathbb{R}^{n})$ and $H^{s,\varphi}(\Omega)$ are
separable. In the Sobolev case of $\varphi\equiv1$ we will omit the index $\varphi$ in
the designations of these and other $H^{s,\varphi}$-type spaces.

In what follows, we suppose that $\Omega$ is a bounded domain in $\mathbb{R}^{n}$ with
$n\geq2$, and its boundary $\partial\Omega$ is an infinitely smooth closed manifold of
the dimension $n-1$. (Note that domains are defined to be an open and connected sets.)
Consider the classes of H\"ormander spaces
\begin{equation}\label{eq8.2}
\bigl\{H^{s,\varphi}(\Omega):\,s\in\mathbb{R},\,\varphi\in\mathcal{M}\,\bigr\}\quad
\mbox{and}\quad \bigl\{H^{s,\varphi}_{\overline{\Omega}}(\mathbb{R}^{n}):
\,s\in\mathbb{R},\,\varphi\in\mathcal{M}\,\bigr\}.
\end{equation}
The space $H^{s,\varphi}(\Omega)$ consists of distributions given in $\Omega$,
whereas the space $H^{s,\varphi}_{\overline{\Omega}}(\mathbb{R}^{n})$ consists of
distributions supported on $\overline{\Omega}$.

\begin{definition}\label{def8.3}
The classes appearing in \eqref{eq8.2} are called the refined Sobolev scales over
$\Omega$ and $\overline{\Omega}$ respectively.
\end{definition}

\subsection{The interpolation properties}\label{sec8.2}

The scales \eqref{eq8.2} have the interpolation properties analogous to those the
refined Sobolev scale over $\mathbb{R}^{n}$ possesses.

\begin{theorem}\label{th8.1}
Let a function $\varphi\in\mathcal{M}$ and positive numbers $\varepsilon,\delta$ be
given, and let the interpolation parameter $\psi\in\mathcal{B}$ be defined by
\eqref{eq3.8}. Then, for each $s\in\mathbb{R}$, the following equalities of spaces
with equivalence of norms in them hold:
\begin{gather}\label{eq8.3}
\bigl[H^{s-\varepsilon}(\Omega),H^{s+\delta}(\Omega)\bigl]_{\psi}=H^{s,\varphi}(\Omega)\\
\bigl[H^{s-\varepsilon}_{\overline{\Omega}}(\mathbb{R}^{n}),
H^{s+\delta}_{\overline{\Omega}}(\mathbb{R}^{n})\bigr]_{\psi}
=H^{s,\varphi}_{\overline{\Omega}}(\mathbb{R}^{n}).\label{eq8.4}
\end{gather}
\end{theorem}

We will deduce this theorem from Theorem \ref{th3.4} with the help of the following
result concerning the interpolation of subspaces and factor spaces.

\begin{proposition}\label{prop8.1}
Let $X=[X_{0},X_{1}]$ be an admissible couple of Hilbert spaces, and $Y_{0}$ be a
subspace in $X_{0}$. Then $Y_{1}:=X_{1}\cap Y_{0}$ is a subspace in $X_{1}$. Suppose
that there exists a linear mapping $P$ which is a projector of $X_{j}$ onto $Y_{j}$
for $j=0,\,1$. Then the couples $[Y_{0},Y_{1}]$ and $[X_{0}/Y_{0},X_{1}/Y_{1}]$ are
admissible, and, for each interpolation parameter $\psi\in\mathcal{B}$, the
following equalities of spaces up to equivalence of norms in them hold:
$$
[Y_{0},Y_{1}]_{\psi}=X_{\psi}\cap Y_{0}, \quad
[X_{0}/Y_{0},X_{1}/Y_{1}]_{\psi}=X_{\psi}/(X_{\psi}\cap Y_{0}).
$$
Here $X_{\psi}\cap Y_{0}$ is a subspace in $X_{\psi}$.
\end{proposition}

Recall that, by definition, subspaces of a Hilbert space are closed, and projectors
on subspaces are, generally, nonorthogonal. Proposition \ref{prop8.1} was proved in
H.~Triebel's monograph \cite[Sec. 1.17]{Triebel95} for arbitrary interpolation
functors given on the category of couples of Banach spaces. The proof for Hilbert
spaces is quite similar.

Let us explain how to prove Theorem \ref{th8.1}. It is known \cite[Sec. 2.10.4,
Theorem~2]{Triebel95} that, for each integer $k>0$, there exists a linear mapping
$P_{k,Q}$ which is a projector of every space $H^{\sigma}(\mathbb{R}^{n})$, with
$|\sigma|<k$, onto the subspace $H^{\sigma}_{Q}(\mathbb{R}^{n})$, where $Q$ is a
closed half-space in $\mathbb{R}^{n}$. Using the local coordinates methods and
$P_{k,Q}$, we can construct a linear mapping that projects
$H^{\sigma}(\mathbb{R}^{n})$ onto $H^{\sigma}_{\widehat{\Omega}}(\mathbb{R}^{n})$
(or onto $H^{\sigma}_{\overline{\Omega}}(\mathbb{R}^{n})$) with $|\sigma|<k$. Now
Theorem \ref{th8.1} follows from  Theorem \ref{th3.4} and Proposition~\ref{prop8.1}.

Reasoning as above we can also deduce analogs of Theorem \ref{th3.5} (on
interpolation) for the refined Sobolev scale given over $\Omega$ or
$\overline{\Omega}$. We will not formulate them.

\subsection{Embeddings and other properties}\label{sec8.3}

Let us consider some other important properties of scales \eqref{eq8.2}. Among these
properties, there are the following embeddings.

\begin{theorem}\label{th8.2}
Let $s\in\mathbb{R}$ and $\varphi,\varphi_{1}\in\mathcal{M}$. The next assertions
are true:
\begin{enumerate}
\item[i)] The set $C^{\infty}(\,\overline{\Omega}\,)$ is dense in
$H^{s,\varphi}(\Omega)$, whereas the set $C^{\infty}_{0}(\Omega)$ is dense in
$H^{s,\varphi}_{\overline{\Omega}}(\mathbb{R}^{n})$.
\item[ii)] For each $\varepsilon>0$, the dense compact embeddings hold:
\begin{equation}\label{eq8.5}
H^{s+\varepsilon,\,\varphi_{1}}(\Omega)\hookrightarrow H^{s,\varphi}(\Omega), \quad
H^{s+\varepsilon,\,\varphi_{1}}_{\overline{\Omega}}(\mathbb{R}^{n})\hookrightarrow
H^{s,\varphi}_{\overline{\Omega}}(\mathbb{R}^{n}).
\end{equation}
\item[iii)] The function $\varphi/\varphi_{1}$ is bounded in a neighbourhood of $+\infty$
if and only if the embeddings \eqref{eq8.5} are valid for $\varepsilon=0$. The
embeddings are continuous and dense. They are compact if and only if
$\varphi(t)/\varphi_{1}(t)\rightarrow0$ as $t\rightarrow+\infty$.
\item[iv)] For every fixed integer $k\geq0$, the inequality \eqref{eq3.6} is equivalent
to the embedding $H^{k+n/2,\varphi}(\Omega)\hookrightarrow
C^{k}(\,\overline{\Omega}\,)$. This embedding is compact.
\end{enumerate}
\end{theorem}

Assertion i) can be deduced directly from the Sobolev case of $\varphi\equiv1$ with
the help of the interpolation Theorem \ref{th8.1}.  Assertions ii)--iv) follow from
Theorem \ref{th3.3} with the exception of the statements on compactness. The
compactness of the embeddings is a consequences of the boundedness of $\Omega$ and
can be proved similarly to the argument of Theorem \ref{th5.2}. The general analogs
of assertions i)--iv) for the H\"ormander inner product spaces parametrized by
arbitrary weight functions were proved by L.R.~Volevich and B.P.~Paneach in
\cite[\S~3, 7, and 8]{VolevichPaneah65}.

Further we examine the properties that exhibit a relation between the refined Sobolev
scales over $\Omega$ and $\overline{\Omega}$. Denote by $H^{s,\varphi}_{0}(\Omega)$ the
closure of $C^{\infty}_{0}(\Omega)$ in $H^{s,\varphi}(\Omega)$. We consider
$H^{s,\varphi}_{0}(\Omega)$ as a Hilbert space with respect to the inner product in
$H^{s,\varphi}(\Omega)$.

\begin{theorem}\label{th8.3}
Let $s\in\mathbb{R}$ and $\varphi\in\mathcal{M}$. The following assertions are true:
\begin{enumerate}
\item[i)] If $s<1/2$, then $C^{\infty}_{0}(\Omega)$ is dense in
$H^{s,\varphi}(\Omega)$, and therefore
$H^{s,\varphi}(\Omega)=H^{s,\varphi}_{0}(\Omega)$.
\item[ii)] If $s>-1/2$ and $s+1/2\notin\mathbb{Z}$, then the restriction mapping
$u\rightarrow u\!\upharpoonright\!\Omega$, $u\in\mathcal{D}'(\mathbb{R}^{n})$,
establishes a homeomorphism of $H^{s,\varphi}_{\overline{\Omega}}(\mathbb{R}^{n})$
onto $H^{s,\varphi}_{0}(\Omega)$.
\item[iii)] The spaces $H^{s,\varphi}(\Omega)$
and $H^{-s,1/\varphi}_{\overline{\Omega}}(\mathbb{R}^{n})$ are mutually dual with
respect to the inner product in $L_{2}(\Omega)$.
\item[iv)] Suppose that $s<1/2$ and $s-1/2\notin\mathbb{Z}$. Then the spaces
$H^{s,\varphi}(\Omega)$ and $H^{-s,1/\varphi}_{0}(\Omega)$ are mutually dual, up to
equivalence of norms, with respect to the inner product in $L_{2}(\Omega)$.
Therefore the space $H^{s,\varphi}(\Omega)$ coincides, up to equivalence of norms,
with the factor space $H^{s,\varphi}_{\overline{\Omega}}(\mathbb{R}^{n})/
H^{s,\varphi}_{\partial\Omega}(\mathbb{R}^{n})$ dual to
$H^{-s,1/\varphi}_{0}(\Omega)$; i.e.,
$H^{s,\varphi}(\Omega)=\{u\!\upharpoonright\!\Omega:\,u\in
H^{s,\varphi}_{\overline{\Omega}}(\mathbb{R}^{n})\}$.
\end{enumerate}
\end{theorem}

This theorem is known in the Sobolev case of $\varphi\equiv1$; see, e.g., \cite[Sec.
4.3.2 and 4.8]{Triebel95}. In general, assertion i) follows from the Sobolev case in
view of Theorem \ref{th8.2}~ii), where $\varphi_{1}\equiv1$; assertion ii) is
deduced with the help of the interpolation Theorem \ref{th8.1}; assertion iii)
results from Theorem \ref{th3.3} iv); finally, assertion iv) follows from ii) and
iii). To deduce assertion ii) we need, besides \eqref{eq8.4}, the interpolation
formula
\begin{equation}\label{eq8.6}
\bigl[H^{s-\varepsilon}_{0}(\Omega),H^{s+\delta}_{0}(\Omega)\bigl]_{\psi}=
H^{s,\varphi}(\Omega)\cap H^{s-\varepsilon}_{0}(\Omega)=H^{s,\varphi}_{0}(\Omega).
\end{equation}
Here $\varepsilon$, $\delta$ are positive numbers such that
$[s-\varepsilon,s+\delta]$ is disjoint from $\mathbb{Z}-1/2$, and $\psi$ is the
interpolation parameter defined by \eqref{eq3.8}. Formula \eqref{eq8.6} follows from
\eqref{eq8.3} and Proposition \ref{prop8.1} (interpolation of subspaces) because
there exist a linear mapping that projects $H^{\sigma}(\Omega)$ onto
$H^{\sigma}_{0}(\Omega)$ if $\sigma$ runs over $[s-\varepsilon,s+\delta]$. The
mapping is constructed in \cite[Lemma 5.4.4 with regard for Theorem
4.7.1]{Triebel95}.

Note that if $s$ is half-integer, then assertions ii) and iv) fail to hold at least
for $\varphi\equiv1$ \cite[Sec. 4.3.2, Remark~2]{Triebel95}.

\begin{remark}\label{rem8.1}
In the literature, there are three different definitions of the Sobolev space of
negative order $s$ over $\Omega$. The first of them coincides with Definition
\ref{def8.2} for $\varphi\equiv1$ (\cite[Sec. A.4]{Grubb96} and \cite[Sec.
4.2.1]{Triebel95}). The second defines this space as the dual of
$H^{-s}_{0}(\Omega)$ (\cite[Ch.~II, \S~1, Sec.~5]{FunctionalAnalysis72} and
\cite[Ch.~1, Sec. 12.1]{LionsMagenes72}), whereas the third defines it as the dual
of $H^{-s}(\Omega)$ (\cite[Ch. XIV, \S~3]{BerezanskySheftelUs96b} and \cite[Sec.
1.10]{Roitberg96}), the both duality being with respect to the inner product in
$L_{2}(\Omega)$. By Theorem \ref{th8.3} iii) and iv), the first and second
definitions are tantamount if (and only if) $s-1/2\notin\mathbb{Z}$, but the third
gives $H^{s}_{\overline{\Omega}}(\mathbb{R}^{n})$ and, therefore, are not equivalent
to them for $s<-1/2$. If $-1/2<s<0$, then all the three definitions are tantamount
in view of Theorem \ref{th8.3} i) and ii). They are suitable in various situations
appearing in the theory of elliptic boundary problems. The situations will occur
below when we will investigate these problems in the refined Sobolev scale. We chose
Definition \ref{def8.2} to introduce the H\"ormander spaces over $\Omega$ because it
is universal; i.e., it allows us to define the space
$X(\Omega)\hookrightarrow\mathcal{D}'(\Omega)$ if we have an arbitrary function
Banach space $X(\mathbb{R}^{n})\hookrightarrow\mathcal{D}'(\mathbb{R}^{n})$ instead
of $H^{s,\varphi}(\mathbb{R}^{n})$ (embeddings being continuous).
\end{remark}

\subsection{Traces}\label{sec8.4}

We now study the traces of functions $f\in H^{s,\varphi}(\Omega)$ and their normal
derivatives on the boundary $\partial\Omega$. The traces belong to certain spaces
from the refined Sobolev scale on $\partial\Omega$. This scale was defined in
Section \ref{sec5.1} because $\partial\Omega$ is a closed infinitely smooth oriented
manifold (of dimension $n-1$). Further we use the notation
$D_{\nu}:=i\,\partial/\partial\nu$, where $\nu$ is the field of unit vectors of
inner normals to the boundary $\partial\Omega$; this field is given in a
neighbourhood of $\partial\Omega$. (For us it will be more suitable to use $D_{\nu}$
instead of $\partial/\partial\nu$; see Sec. \ref{sec11} below.)

\begin{theorem}\label{th8.4}
Let an integer $r\geq1$, real number $s>r-1/2$, and function $\varphi\in\mathcal{M}$
be given. Then the mapping
\begin{equation}\label{eq8.7}
R_{r}:u\mapsto\bigl((D_{\nu}^{k-1}u)\!\upharpoonright\!\partial\Omega:\,
k=1,\ldots,r\bigr),\quad u\in C^{\infty}(\,\overline{\Omega}),
\end{equation}
extends uniquely to a continuous linear operator
$$
R_{r}:\,H^{s,\varphi}(\Omega)\rightarrow
\bigoplus_{k=1}^{r}\,H^{s-k+1/2,\,\varphi}(\partial\Omega)=:
\mathcal{H}_{s,\varphi}^{r}(\partial\Omega).
$$
The operator \eqref{eq8.7} has a right inverse continuous linear operator
$\Upsilon_{r}:\mathcal{H}_{s,\varphi}^{r}(\partial\Omega)\rightarrow
H^{s,\varphi}(\Omega)$ such that the mapping $\Upsilon_{r}$ does not depend on $s$
and $\varphi$.
\end{theorem}

Theorem \ref{th8.4} is known in the Sobolev case of $\varphi\equiv1$; see, e.g.,
\cite[Ch.~1, Sec. 9.2]{LionsMagenes72} or \cite[Sec. 4.7.1]{Triebel95}. For
arbitrary $\varphi\in\mathcal{M}$, the theorem follows from this case by the
interpolation Theorems \ref{th5.3} and \ref{th8.1}.

It useful to note that
\begin{equation}\label{eq8.8}
H^{s,\varphi}_{0}(\Omega)=\bigl\{u\in
H^{s,\varphi}(\Omega):\,R_{r}u=0\bigr\}\quad\mbox{if}\quad r-1/2<s<r+1/2;
\end{equation}
here the integer $r\geq1$. This formula is known in the $\varphi\equiv1$ case, the
equality $r+1/2=s$ being possible; see, e.g., \cite[Ch.~1, Sec.
11.4]{LionsMagenes72} or \cite[Sec. 4.7.1]{Triebel95}. In general, \eqref{eq8.8}
follows from the Sobolev case by \eqref{eq8.6}.

If $s>1/2$ and $\varphi\in\mathcal{M}$, then, by Theorem \ref{th8.4}, a trace
$u\!\upharpoonright\!\partial\Omega:=R_{1}u$ of each function $u\in
H^{s,\varphi}(\Omega)$ on the boundary $\partial\Omega$ exists and belongs to the
space $H^{s-1/2,\varphi}(\partial\Omega)$. Moreover, we get the following
description of this space in terms of traces. Put $\sigma:=s-1/2$.

\begin{corollary}\label{cor8.1}
Let $\sigma>0$ and $\varphi\in\mathcal{M}$. Then
\begin{gather*}
H^{\sigma,\varphi}(\partial\Omega)=\{g:=u\!\upharpoonright\!\partial\Omega:\,
u\in H^{\sigma+1/2,\varphi}(\Omega)\},\\
\|g\|_{H^{\sigma,\varphi}(\partial\Omega)}\asymp
\inf\,\bigl\{\|u\|_{H^{\sigma+1/2,\varphi}(\Omega)}:
\,u\!\upharpoonright\!\partial\Omega=g\bigr\}.
\end{gather*}
\end{corollary}

If $s<1/2$ and $\varphi\in\mathcal{M}$, then the trace mapping
\begin{equation}\label{eq8.9}
R_{1}:u\mapsto u\!\upharpoonright\!\partial\Omega,\quad u\in
C^{\infty}(\,\overline{\Omega}),
\end{equation}
has not a continuous extension
$R_{1}:H^{s,\varphi}(\Omega)\rightarrow\mathcal{D}'(\partial\Omega)$. Indeed, if
this extension existed, we would get, by Theorem \ref{th8.3} i), the equality
$R_{1}u=0$ on $\partial\Omega$ for each $u\in H^{s,\varphi}(\Omega)$. But this
equality fails to hold, e.g., for the function $u\equiv1$.

In the limiting case of $s=1/2$, we have the next criterion for the trace operator
$R_{1}$ to be well defined on $H^{1/2,\varphi}(\Omega)$.

\begin{theorem}\label{th8.5}
Let $\varphi\in\mathcal{M}$. The following assertions are true:
\begin{enumerate}
\item[i)] The function $\varphi$ satisfies \eqref{eq3.6} if an only if the mapping
\eqref{eq8.9} is a continuous operator from the space
$C^{\infty}(\,\overline{\Omega})$ endowed with the topology of
$H^{1/2,\varphi}(\Omega)$ to the space $\mathcal{D}'(\partial\Omega)$.
\item[ii)] Moreover, if $\varphi$ satisfies \eqref{eq3.6}, then the mapping \eqref{eq8.9}
extends uniquely to a continuous linear operator
$R_{1}:H^{1/2,\varphi}(\Omega)\rightarrow H^{0,\varphi_{0}}(\partial\Omega)$, where
$\varphi_{0}\in\mathcal{M}$ is given by the formula
\begin{equation}\label{eq8.10}
\varphi_{0}(\tau):=
\biggl(\:\int\limits_{\tau}^{\infty}\frac{d\,t}{t\,\varphi^{2}(t)}\:\biggr)^{-1/2}
\quad\mbox{for}\quad\tau\geq1.
\end{equation}
This operator has a right inverse continuous linear operator
$$
\Upsilon_{1,\varphi}:H^{0,\varphi_{0}}(\partial\Omega)\rightarrow
H^{1/2,\varphi}(\Omega),
$$
the map $\Upsilon_{1,\varphi}$ depending on $\varphi$.
\end{enumerate}
\end{theorem}

Theorem \ref{th8.5} follows from the trace theorems proved by L.~H\"ormander
\cite[Sec. 2.2, Theorem 2.2.8]{Hermander63} and L.R.~Volevich, B.P.~Paneah
\cite[\S~6, Theorems 6.1 and 6.2]{VolevichPaneah65}. Indeed, consider a H\"ormander
space $B_{p,\mu}(\mathbb{R}^{n})$ for some weight function $\mu$, and let $U$ be a
neighbourhood of the origin in $\mathbb{R}^{n}$. We write points
$x\in\mathbb{R}^{n}$ as $x=(x',x_{n})$ with $x'\in\mathbb{R}^{n-1}$ and
$x_{n}\in\mathbb{R}$. According to the trace theorems, the condition
\begin{equation}\label{eq8.11}
\nu^{-2}(\xi\,'):= \int\limits_{-\infty}^{\infty}\mu^{-2}(\xi',\xi_{n})\,d\xi_{n}
<\infty\quad\mbox{for all}\quad\xi'\in\mathbb{R}^{n-1}
\end{equation}
holds true if and only if the mapping $u(x',x_{n})\rightarrow u(x',0)$ is a
continuous operator from the space $C^{\infty}_{0}(U)$ endowed with the topology of
$B_{2,\mu}(\mathbb{R}^{n})$ to the space $\mathcal{D}'(\mathbb{R}^{n-1})$. (In
\eqref{eq8.11} the phrase `for all' can be replaced with `for a certain'.) Moreover,
if \eqref{eq8.11} holds and $U=\mathbb{R}^{n}$, then this mapping extends by
continuity to a bounded operator from $B_{2,\mu}(\mathbb{R}^{n})$ to
$B_{2,\nu}(\mathbb{R}^{n-1})$; the operator has a linear bounded right-inverse.
Whence, by setting $\mu(\xi):=\langle\xi\rangle^{1/2}\varphi(\langle\xi\rangle)$ and
observing that \eqref{eq8.11} $\Leftrightarrow$ \eqref{eq3.6} with
$\nu(\xi')\asymp\varphi_{0}(\langle\xi'\rangle)$, we can deduce Theorem \ref{th8.5}
with the help of the local coordinates method.

Let us note that the domain $\Omega$ is a special case of an infinitely smooth
compact manifold with boundary. The refined Sobolev scale over such a manifold was
introduced and investigated by authors in \cite[Sec.~3]{06UMJ3}. Specifically, the
interpolation Theorem \ref{th8.1} and the traces Theorems \ref{th8.4} (for $r=1$)
and \ref{th8.5} were proved.

\subsection{Riggings}

We recall the important notion of a Hilbert rigging, which has various applications,
specifically, in the theory of elliptic operators; see, e.g., \cite[Ch.~I,
\S~1]{Berezansky68} and \cite[Ch.~XIV, \S~1]{BerezanskySheftelUs96b}. Let $H$ and
$H_{+}$ be Hilbert spaces such that the dense continuous embedding
$H_{+}\hookrightarrow H$ holds. Denote by $H_{-}$ the completion of $H$ with respect
to the norm
$$
\|f\|_{H_{-}}:=\sup\,\biggl\{\,\frac{|(f,u)_{H}|}{\|u\|_{H_{+}}}:\,u\in
H_{+}\,\biggr\},\quad f\in H.
$$
It is known the following: this norm and the space $H_{-}$ are Hilbert; the spaces
$H_{+}$ and $H_{-}$ are mutually dual with respect to the inner product in $H$; the dense
continuous embeddings $H_{+}\hookrightarrow H\hookrightarrow H_{-}$ hold.

\begin{definition}\label{def8.4}
The chain $H_{-}\hookleftarrow H\hookleftarrow H_{+}$ is said to be a (Hilbert)
rigging of $H$ with $H_{+}$ and $H_{-}$. In this rigging, $H_{-}$, $H$ and $H_{+}$
are called negative, zero, and positive spaces respectively.
\end{definition}

According to Theorem \ref{th8.3} iii) we have the following rigging of
$L_{2}(\Omega)$ with some H\"ormander spaces:
\begin{equation}\label{eq8.12}
H^{-s,1/\varphi}_{\overline{\Omega}}(\mathbb{R}^{n})\hookleftarrow
L_{2}(\Omega)\hookleftarrow H^{s,\varphi}(\Omega),\quad s>0,\;\varphi\in\mathcal{M}.
\end{equation}
Here we naturally identify $L_{2}(\Omega)$ with
$H^{0}_{\overline{\Omega}}(\mathbb{R}^{n})$ (extending the functions $f\in
L_{2}(\Omega)$ by zero).

In some applications to elliptic problems, it is useful to have a scale consisting
of both negative and positive spaces in \eqref{eq8.12}. For this purpose we
introduce the uniform notation
\begin{equation}\label{eq8.13}
H^{s,\varphi,(0)}(\Omega):=
\begin{cases}
\;H^{s,\varphi}(\Omega)\;\; & \text{for}\;\;s\geq0, \\
\;H^{s,\varphi}_{\overline{\Omega}}(\mathbb{R}^{n}) & \text{for}\;\;s<0,
\end{cases}
\end{equation}
with $\varphi\in\mathcal{M}$, and form the scale of Hilbert spaces
\begin{equation}\label{eq8.14}
\bigl\{H^{s,\varphi,(0)}(\Omega):\,s\in\mathbb{R},\,\varphi\in\mathcal{M}\,\bigr\}.
\end{equation}

In the Sobolev case of $\varphi\equiv1$ the rigging \eqref{eq8.12} and the scale of
spaces $H^{s,(0)}(\Omega):=H^{s,1,(0)}(\Omega)$, $s\in\mathbb{R}$, were used by
Yu.M.~Berezansky, S.G.~Krein, Ya.A.~Roitberg \cite{Berezansky68,
BerezanskyKreinRoitberg63, BerezanskySheftelUs96b, Roitberg64} and M.~Schechter
\cite{Schechter63} in the elliptic theory. (They also considered the Banach
$L_{p}$-analogs of these spaces with $1<p<\infty$ and denoted negative spaces in the
same manner as positive ones but with negative index $s$, e.g. $H^{s}(\Omega)$; see
Remark \ref{rem8.1} above.)

Properties of the scale \eqref{eq8.14} are inherited from the refined Sobolev scales
over $\Omega$ and $\overline{\Omega}$. Now we dwell on the properties that link
negative and positive spaces to each other.

When dealing with the scale \eqref{eq8.14}, it is suitable to identify each function
$f\in C^{\infty}(\,\overline{\Omega}\,)$ with its extension by zero
\begin{equation}\label{eq8.15}
\mathcal{O}f(x):=
\begin{cases}
\;f(x) &\;\; \text{for}\;\; x\in\overline{\Omega}, \\
\;0 &\;\; \text{for}\;\; x\in\mathbb{R}^{n}\setminus\overline{\Omega}.
\end{cases}
\end{equation}
The extension is a regular distribution belonging to
$H^{s,\varphi}_{\overline{\Omega}}(\mathbb{R}^{n})$ for $s<0$. Now the set
$C^{\infty}(\,\overline{\Omega}\,)$ is dense in every space
$H^{s,\varphi,(0)}(\Omega)$ from \eqref{eq8.14} in view of Theorem \ref{th8.2} i).
This allow us to consider the continuous embeddings of spaces pertaining to
\eqref{eq8.14} and viewed as the completions of the same linear manifold,
$C^{\infty}(\,\overline{\Omega}\,)$, with respect to different norms. (The general
theory of such embeddings is in \cite[Ch. XIV, \S~7]{BerezanskySheftelUs96b}). So,
by Theorem \ref{th8.2} ii) and formula \eqref{eq8.17} given below, we have the dense
compact embeddings in the scale \eqref{eq8.14}:
\begin{equation}\label{eq8.16}
H^{s_{1},\varphi_{1},(0)}(\Omega)\hookrightarrow
H^{s,\varphi,(0)}(\Omega),\quad-\infty<s<s_{1}<\infty\;\;\mbox{and}
\;\;\varphi,\varphi_{1}\in\mathcal{M}.
\end{equation}

Note that in the $|s|<1/2$ case the spaces $H^{s,\varphi}(\Omega)$ and
$H^{s,\varphi}_{\overline{\Omega}}(\mathbb{R}^{n})$ can be considered as completions
of $C^{\infty}_{0}(\Omega)$ with respect to equivalent norms in view of Theorem
\ref{th8.3} i) and ii). Hence, up to equivalence of norms,
\begin{equation}\label{eq8.17}
H^{s,\varphi}_{\overline{\Omega}}(\mathbb{R}^{n})=H^{s,\varphi,(0)}(\Omega)=
H^{s,\varphi}(\Omega)\quad\mbox{for}\;\;|s|<1/2,\;\;\varphi\in\mathcal{M}.
\end{equation}

It follows from this result and Theorem \ref{th8.3} iii) that the spaces
$H^{s,\varphi,(0)}(\Omega)$ and $H^{-s,1/\varphi,(0)}(\Omega)$ are mutually dual
with respect to the inner product in $L_{2}(\Omega)$ for every $s\in\mathbb{R}$ and
$\varphi\in\mathcal{M}$, the duality being up to equivalence of norms if $s=0$.

The scale \eqref{eq8.14} has an interpolation property analogous to that stated in
Theorem \ref{th8.1}.

\begin{theorem}\label{th8.6}
Under the conditions of Theorem $\ref{th8.1}$ we have
\begin{equation}\label{eq8.18}
\bigl[H^{s-\varepsilon,(0)}(\Omega),H^{s+\delta,(0)}(\Omega)\bigl]_{\psi}=
H^{s,\varphi,(0)}(\Omega)\quad\mbox{for all}\;\;s\in\mathbb{R}
\end{equation}
up to equivalence of norms in the spaces.
\end{theorem}

If $s-\varepsilon>-1/2$ or $s+\delta<1/2$, then \eqref{eq8.18} holds by Theorem
\ref{th8.1} and \eqref{eq8.17}. If $\varphi\equiv1$, then \eqref{eq8.18} is proved
in \cite[Ch.~1, Sec. 12.5, Theorem 12.5]{LionsMagenes72}. The general case can be
reduced to the previous ones by the reiterated interpolation.

\section{Elliptic boundary-value problems on the one-sided scale}\label{sec9}

In Sections \ref{sec9}--\ref{sec12}, we will investigate regular elliptic
boundary-value problems on various scales of H\"ormander spaces. We begin with the
one-sided refined Sobolev scale consisting of the spaces $H^{s,\varphi}(\Omega)$ for
sufficiently large~$s$.

\subsection{The statement of the boundary-value problem}\label{sec9.1}

Recall that $\Omega$ is a boun\-ded domain in $\mathbb{R}^{n}$ with $n\geq2$, and
its boundary $\partial\Omega$ is an infinitely smooth closed manifold of the
dimension $n-1$. Let $\nu(x)$ denote the unit vector of the inner normal to
$\partial\Omega$ at a point $x\in\partial\Omega$.

We consider the following boundary-value problem in $\Omega$:
\begin{gather}\label{eq9.1}
L\,u\equiv\sum_{|\mu|\leq2q}\,l_{\mu}(x)\,D^{\mu}u=f\quad\mbox{in}\;\;\Omega,\\
B_{j}\,u\equiv\sum_{|\mu|\leq m_{j}}\,b_{j,\mu}(x)\,D^{\mu}u=
g_{j}\quad\mbox{on}\;\;\partial\Omega,\quad j=1,\ldots,q. \label{eq9.2}
\end{gather}
Here $L=L(x,D)$, $x\in\overline{\Omega}$, and $B_{j}=B_{j}(x,D)$, $x\in\partial\Omega$,
are linear partial differential expressions with complex-valued coefficients $l_{\mu}\in
C^{\infty}(\,\overline{\Omega}\,)$ and $b_{j,\mu}\in C^{\infty}(\partial\Omega)$. We
suppose that $\mathrm{ord}\,L=2q$ is an even positive number and
$\mathrm{ord}\,B_{j}=m_{j}\leq2q-1$ for all $j=1,\ldots,q$. Set $B:=(B_{1},\ldots,B_{q})$
and $m:=\max\,\{m_{1},\ldots,m_{q}\}$.

Note that we use the standard notation in \ref{eq9.1} and \ref{eq9.2}; namely, for a
multi-index $\mu=(\mu_{1},\ldots,\mu_{n})$ we let $|\mu|:=\mu_{1}+\ldots+\mu_{n}$
and $D^{\mu}:=D_{1}^{\mu_{1}}\ldots D_{n}^{\mu_{n}}$, with
$D_{k}:=i\,\partial/\partial x_{k}$ for $k=1,\ldots,n$ and
$x=(x_{1},\ldots,x_{n})\in\mathbb{R}^{n}$.

\begin{lemma}\label{lem9.1}
Let $s>m+1/2$ and $\varphi\in\mathcal{M}$.  Then the mapping
\begin{equation}\label{eq9.3}
(L,B):\,u\rightarrow(Lu,B_{1}u,\ldots,B_{q}u),\quad u\in
C^{\infty}(\,\overline{\Omega}\,),
\end{equation}
extends uniquely to a continuous linear operator
\begin{equation}\label{eq9.4}
(L,B):\,H^{s,\varphi}(\Omega)\rightarrow
H^{s-2q,\,\varphi}(\Omega)\oplus\bigoplus_{j=1}^{q}\,
H^{s-m_{j}-1/2,\,\varphi}(\partial\Omega)=:
\mathcal{H}_{s,\varphi}(\Omega,\partial\Omega).
\end{equation}
\end{lemma}

Note the differential operator $L$ maps $H^{s,\varphi}(\Omega)$ continuously into
$H^{s-2g,\varphi}(\Omega)$ for each real~$s$, whereas the boundary differential
operator $B_{j}$ maps $H^{s,\varphi}(\Omega)$ continuously into
$H^{s-m_{j}-1/2,\varphi}(\partial\Omega)$ provided that $s>m_{j}+1/2$. This is well
known in the Sobolev case of $\varphi\equiv1$ (see, e.g., \cite[Sec.
B.2]{Hermander85}), whence the case of an arbitrary $\varphi\in\mathcal{M}$ is got
by the interpolation in view of Theorems \ref{th5.3} and \ref{th8.1}. The mentioned
continuity of $B_{j}$ also results from the trace Theorem \ref{th8.4} (the $r=1$
case). If $s\leq m+1/2$, then Lemma \ref{lem9.1} fails to hold (see Section
\ref{sec8.4}). In the limiting case of $s=m+1/2$ an analog of the lemma is valid
provided the function $\varphi$ satisfies \eqref{eq3.6} and we exchange the space
$\mathcal{H}_{s,\varphi}(\Omega,\partial\Omega)$ for another (see Section
\ref{sec9.3} below).

We will investigate properties of the operator \eqref{eq9.4} under the assumption
that the boundary-value problem \eqref{eq9.1}, \eqref{eq9.2} is regular elliptic in
$\Omega$. Recall some relevant notions; see, e.g., \cite[Ch.~III,
\S~6]{FunctionalAnalysis72} or \cite[Ch.~2, Sec. 1 and 2]{LionsMagenes72}.

The principal symbols of the partial differential expressions $L(x,D)$, with
$x\in\overline{\Omega}$, and $B_{j}(x,D)$, with $x\in\partial\Omega$, are defined as
follows:
$$
L^{(0)}(x,\xi):=\sum_{|\mu|=2q}l_{\mu}(x)\,\xi^{\mu},\quad\quad
B_{j}^{(0)}(x,\xi):=\sum_{|\mu|=\,m_{j}}b_{j,\mu}(x)\,\xi^{\mu}.
$$
They are homogeneous polynomials in $\xi=(\xi_{1},\ldots,\xi_{n})\in\mathbb{C}^n$; here
as usual $\xi^{\mu}:=\xi_{1}^{\mu_{1}}\ldots\xi_{n}^{\mu_{n}}$.

\begin{definition}\label{def9.1}
The boundary-value problem \eqref{eq9.1}, \eqref{eq9.2} is said to be regular (or
normal) elliptic in $\Omega$ if the following conditions are satisfied:
\begin{enumerate}
\item[i)] The expression $L$ is proper elliptic on $\overline{\Omega}$; i.e.,
for each point $x\in\overline{\Omega}$ and for all linearly independent vectors
$\xi',\xi''\in\mathbb{R}^{n}$, the polynomial $L^{(0)}(x,\xi'+\tau\xi'')$ in
$\tau\in\mathbb{C}$ has $q$ roots $\tau^{+}_{j}(x;\xi',\xi'')$, $j=1,\ldots,q$, with
positive imaginary part and $q$ roots with negative imaginary part, each root being
taken the number of times equal to its multiplicity.
\item[ii)] The system $\{B_{1},\ldots,B_{q}\}$ satisfies the Lopatinsky
condition with respect to $L$ on $\partial\Omega$; i.e., for an arbitrary point
$x\in\partial\Omega$ and for each vector $\xi\neq0$ tangent to $\partial\Omega$ at
$x$, the polynomials $B_{j}^{(0)}(x,\xi+\tau\nu(x))$, $j=1,\ldots,q$, in $\tau$ are
linearly independent modulo
$\prod_{j=1}^{q}\bigl(\tau-\tau^{+}_{j}(x;\xi,\nu(x))\bigr)$.
\item[iii)] The system $\{B_{1},\ldots,B_{q}\}$ is normal; i.e., the orders
$m_{j}$ are pairwise distinct, and $B_{j}^{(0)}(x,\nu(x))\neq\nobreak0$ for each
$x\in\partial\Omega$.
\end{enumerate}
\end{definition}

\begin{remark}\label{rem9.1}
It follows from condition i) that $L^{(0)}(x,\xi)\neq0$ for each point
$x\in\overline{\Omega}$ and nonzero vector $\xi\in\mathbb{R}^{n}$, i.e. $L$ is
elliptic on $\overline{\Omega}$. If $n\geq3$, then the ellipticity condition is
equivalent to i). The equivalence also holds if $n=2$ and all coefficients of $L$
are real-valued. Not more that there are various equivalent statements of the
Lopatinsky condition; see \cite[Sec. 1.2 and 1.3]{Agranovich97}.
\end{remark}

\begin{example}\label{ex9.1}
Let $L$ satisfy condition i), and let
$B_{j}u:=\partial^{k+j-1}u/\partial\nu^{k+j-1}$, with $j=1,\ldots,q$, for some
$k\in\mathbb{Z}$, $0\leq k\leq q$. Then the boundary-value problem \eqref{eq9.1},
\eqref{eq9.2} is regular elliptic. If $k=0$, we have the Dirichlet boundary-value
problem.
\end{example}

In what follows the boundary-value problem \eqref{eq9.1}, \eqref{eq9.2} is supposed
to be regular elliptic in $\Omega$.

To describe the range of the operator \eqref{eq9.4} we consider the boundary-value
problem
\begin{gather}\label{eq9.5}
L^{+}v\equiv\sum_{|\mu|\leq\,2q}D^{\mu}(\overline{l_{\mu}(x)}\,v)=\omega
\quad\mbox{in}\;\;\Omega,\\
B^{+}_{j}v=h_{j}\quad\mbox{on}\;\;\partial\Omega,\quad j=1,\ldots,q, \label{eq9.6}
\end{gather}
that is formally adjoint to the problem \eqref{eq9.1}, \eqref{eq9.2} with respect to
the Green formula
\begin{equation}\label{eq9.7}
(Lu,v)_{\Omega}+\sum_{j=1}^{q}\;(B_{j}u,\,C_{j}^{+}v)_{\partial\Omega}
=(u,L^{+}v)_{\Omega}+\sum_{j=1}^{q}\;(C_{j}u,\,B_{j}^{+}v)_{\partial\Omega},
\end{equation}
where $u,v\in C^{\infty}(\,\overline{\Omega}\,)$. Here $\{B^{+}_{j}\}$, $\{C_{j}\}$,
$\{C^{+}_{j}\}$ are some normal systems of linear partial differential boundary
expressions with coefficients from $C^{\infty}(\partial\Omega)$; the orders of
expressions $B_{j}^{\pm}$, $C_{j}^{\pm}$, $j=1,\ldots,q$, run over the set
$\{0,1,\ldots,2q-1\}$ and satisfy the equality
$$
\mathrm{ord}\,B_{j}+\mathrm{ord}\,C^{+}_{j}=
\mathrm{ord}\,C_{j}+\mathrm{ord}\,B^{+}_{j}=2q-1.
$$
We denote $m_{j}^{+}:=\mathrm{ord}\,B_{j}^{+}$. In \eqref{eq9.7} and bellow, the
notations $(\cdot,\cdot)_{\Omega}$ and $(\cdot,\cdot)_{\partial\Omega}$ stand for
the inner products in the spaces $L_{2}(\Omega)$ and $L_{2}(\partial\Omega)$
respectively, and for extensions by continuity of these products as well.

The expression $L^{+}$ is said to be formally adjoint to $L$, whereas the system
$\{B^{+}_{j}\}$ is said to be adjoint to $\{B_{j}\}$ with respect to $L$. Note that
$\{B^{+}_{j}\}$ is not uniquely defined.

We set
\begin{gather*}
\mathcal{N}:=\{u\in
C^{\infty}(\,\overline{\Omega}\,):\,Lu=0\;\,\mbox{in}\;\,\Omega,\;
B_{j}u=0\;\,\mbox{on}\;\,\partial\Omega\;\,\mbox{for all}\;\,j=1,\ldots,q\},\\
\mathcal{N}^{+}:=\{v\in
C^{\infty}(\,\overline{\Omega}\,):\,L^{+}v=0\;\,\mbox{in}\;\,\Omega,\;
B^{+}_{j}v=0\;\,\mbox{on}\;\,\partial\Omega\;\,\mbox{for all}\;\,j=1,\ldots,q\}.
\end{gather*}
Since both the problems \eqref{eq9.1}, \eqref{eq9.2} and \eqref{eq9.5},
\eqref{eq9.6} are regular elliptic, both the spaces $\mathcal{N}$ and
$\mathcal{N}^{+}$ are finite dimensional. Note that the space $\mathcal{N}^{+}$ does
not not depend on the choice of the system $\{B^{+}_{j}\}$ adjoint to $\{B_{j}\}$.

\begin{example}\label{ex9.2}
If the boundary-value problem \eqref{eq9.1}, \eqref{eq9.2} is Dirichlet, then the
formally adjoint problem is also Dirichlet, with
$\dim\mathcal{N}=\dim\mathcal{N}^{+}$.
\end{example}

\subsection{The operator properties}\label{sec9.2}

Now we investigate properties of the operator \eqref{eq9.4} corresponding to the
regular elliptic boundary-value problem \eqref{eq9.1}, \eqref{eq9.2}.

\begin{theorem}\label{th9.1}
Let $s>m+1/2$ and $\varphi\in\mathcal{M}$. Then the bounded linear operator
\eqref{eq9.4} is Fredholm, its kernel coincides with $\mathcal{N}$, and its range
consists of all the vectors
$(f,g_{1},\ldots,g_{q})\in\mathcal{H}_{s,\varphi}(\Omega,\partial\Omega)$ such that
\begin{equation}\label{eq9.8}
(f,v)_{\Omega}+\sum_{j=1}^{q}\,(g_{j},C^{+}_{j}v)_{\partial\Omega}=0\quad\mbox{for
all}\quad v\in \mathcal{N}^{+}.
\end{equation}
The index of \eqref{eq9.4} is $\dim\mathcal{N}-\dim\mathcal{N}^{+}$ and does not
depend on $s$ and $\varphi$.
\end{theorem}

In the Sobolev case of $\varphi\equiv1$ Theorem \ref{th9.1} is a classical result if
$s\geq2q$; see, e.g., \cite[Ch.~III, \S~6, Subsec.~4]{FunctionalAnalysis72} or
\cite[Ch.~2, Sec.~5.4]{LionsMagenes72}. If $m+1/2<s<2q$, then this theorem is also
true \cite[Ch.~III, Sec. 2.2]{Egorov86}. For an arbitrary $\varphi\in\mathcal{M}$
the theorem can be deduced from the Sobolev case by the interpolation with a
function parameter if we apply Proposition \ref{prop6.1} and Theorems \ref{th5.3}
and \ref{th8.1}.

\begin{remark}\label{rem9.2}
G.~Slenzak \cite{Slenzak74} proved an analog of Theorem \ref{th9.1} for a different
scale of H\"ormander inner product spaces. These spaces are not attached to Sobolev
spaces with the number parameter; the class of the weight functions used by Slenzak
is not described constructively.
\end{remark}

\begin{theorem}\label{th9.2}
Let $s>m+1/2$, $\varphi\in\mathcal{M}$, and $\sigma<s$. Then the following a~priori
estimate holds:
$$
\|u\|_{H^{s,\varphi}(\Omega)}\leq
c\,\bigr(\,\|(L,B)u\|_{\mathcal{H}_{s,\varphi}(\Omega,\partial\Omega)}+
\|u\|_{H^{\sigma,\varphi}(\Omega)}\,\bigl)\quad\mbox{for all}\quad u\in
H^{s,\varphi}(\Omega);
$$
here the number $c>0$ is independent of $u$.
\end{theorem}

This theorem follows from the Fredholm property of the operator \eqref{eq9.4} in
view of Proposition \ref{prop6.2} and the compactness of the embedding
$H^{s,\varphi}(\Omega)\hookrightarrow H^{\sigma,\varphi}(\Omega)$ for $\sigma<s$.

Now we study a local smoothness of a solution $u$ to the boundary-value problem
\eqref{eq9.1}, \eqref{eq9.2} in the refined Sobolev scale. The relevant property
will be formulated as a theorem on the local increase in smoothness.

Let $U$ be an open set in $\mathbb{R}^{n}$; we put
$\Omega_{0}:=U\cap\Omega\neq\varnothing$ and $\Gamma_{0}:=U\cap\partial\Omega$ (the
case were $\Gamma_{0}=\varnothing$ is possible). We introduce the following local
analog of the space $H^{s,\varphi}(\Omega)$ with $s\in\mathbb{R}$ and
$\varphi\in\mathcal{M}$:
\begin{gather*}
H^{s,\varphi}_{\mathrm{loc}}(\Omega_{0},\Gamma_{0}):=
\bigl\{u\in\mathcal{D}'(\Omega):\chi\,u\in H^{s,\varphi}(\Omega)\\
\mbox{for all}\;\;\chi\in C^{\infty}(\overline{\Omega})\;\;\mbox{with}\;\;
\mathrm{supp}\,\chi\subseteq\Omega_{0}\cup\Gamma_{0}\bigr\}.
\end{gather*}
The other local space $H^{s,\varphi}_{\mathrm{loc}}(\Gamma_{0})$, which we need, was
defined in \eqref{eq6.3}.

\begin{theorem}\label{th9.3}
Let  $s>m+1/2$ and $\eta\in\mathcal{M}$. Suppose that the distribution $u\nobreak\in
H^{s,\eta}(\Omega)$ is a solution to the boundary-value problem \eqref{eq9.1},
\eqref{eq9.2}, with
$$
f\in H^{s-2q+\varepsilon,\varphi}_{\mathrm{loc}}(\Omega_{0},\Gamma_{0})
\quad\mbox{and}\quad g_{j}\in
H^{s-m_{j}-1/2+\varepsilon,\varphi}_{\mathrm{loc}}(\Gamma_{0}), \;\;j=1,\ldots,q,
$$
for some $\varepsilon\geq0$ and $\varphi\in\mathcal{M}$. Then $u\in
H^{s+\varepsilon,\varphi}_{\mathrm{loc}}(\Omega_{0},\Gamma_{0})$.
\end{theorem}

In the special case where $\Omega_{0}=\Omega$ and $\Gamma_{0}=\partial\Omega$ we
have the global smoothness increase (i.e. the increase in the domain $\Omega$ up to
its boundary). This case follows from Theorem \ref{th9.1}. Indeed, since the vector
$(f,g)\in\mathcal{H}_{s+\varepsilon,\varphi}(\Omega,\partial\Omega)$ satisfies
\eqref{eq9.8}, we can write $(L,B)v=(f,g)$ for some $v\in
H^{s+\varepsilon,\varphi}(\Omega)$, whence $u-v\in\mathcal{N}$ and $u\in
H^{s+\varepsilon,\varphi}(\Omega)$; here $g:=(g_{1},\ldots,g_{q})$. In general, we
can deduce Theorem \ref{th9.3} from the above case reasoning similar to the proof of
Theorem~\ref{th4.2}. Note, if $\Gamma_{0}=\varnothing$, then we get an interior
smoothness increase (in neighbourhoods of interior points of $\Omega$).

Theorem \ref{th9.3} specifies, with regard to the refined Sobolev scale, the
classical results on a local smoothness of solutions to elliptic boundary-value
problems \cite[Ch. 3, Sec.~4]{Berezansky68}, \cite{Browder56, Nirenberg55,
Schechter61}.

Theorems \ref{th9.3} and \ref{th8.2} iv) imply the following sufficient condition
for the solution $u$ to be classical.

\begin{theorem}\label{th9.4}
Let $s>m+1/2$ and $\eta\in\mathcal{M}$. Suppose that the distribution $u\nobreak\in
H^{s,\eta}(\Omega)$ is a solution to the boundary-value problem \eqref{eq9.1},
\eqref{eq9.2}, where
\begin{gather*}
f\in H^{n/2,\,\varphi}_{\mathrm{loc}}(\Omega,\varnothing)\cap
H^{m-2q+n/2,\,\varphi}(\Omega),\\
g_{j}\in H^{m-m_{j}+(n-1)/2,\,\varphi}(\partial\Omega),\;\;j=1,\ldots,q,
\end{gather*}
for some $\varphi\in\mathcal{M}$. If $\varphi$ satisfies \eqref{eq3.6}, then the
solution $u$ is classical, i.e. $u\in C^{2q}(\Omega)\cap
C^{m}(\,\overline{\Omega}\,)$.
\end{theorem}

Note that the condition \eqref{eq3.6} not only is sufficient for $u$ to be a
classical solution but also is necessary on the class of all the considered
solutions $u$. This follows from Theorem \ref{th3.1}.

\subsection{The limiting case}\label{sec9.3}

This case is $s=m+1/2$. We study it, e.g., for the Dirichlet problem for the Laplace
equation:
$$
\Delta\,u=f\;\;\mbox{in}\;\;\Omega,\quad R_{1}u:=u\!\upharpoonright\!\partial\Omega=g.
$$
This problem is regular elliptic, with $m=0$. Let $\varphi\in\mathcal{M}$. By
Theorem \ref{th8.5}, the mapping $u\mapsto(\Delta\,u,R_{1}u)$, $u\in
C^{\infty}(\,\overline{\Omega}\,)$, extends uniquely to a continuous linear operator
from $H^{1/2,\varphi}(\Omega)$ to
$H^{-3/2,\varphi}(\Omega)\times\mathcal{D}'(\partial\Omega)$ if and only if
$\varphi$ satisfies \eqref{eq3.6}. Suppose the inequality \eqref{eq3.6} is
fulfilled, and $\varphi_{0}\in\mathcal{M}$ is defined by \eqref{eq8.10}. Then we get
the bounded linear operator
\begin{equation}\label{eq9.9}
(\Delta,R_{1}):\,H^{1/2,\varphi}(\Omega)\rightarrow H^{-3/2,\varphi}(\Omega)\oplus
H^{0,\varphi_{0}}(\partial\Omega)=:\mathcal{H}(\Omega,\partial\Omega),
\end{equation}
with $R_{1}(H^{1/2,\varphi}(\Omega))$ being equal to
$H^{0,\varphi_{0}}(\partial\Omega)$. It is reasonable to ask whether this operator
is Fredholm or not. The answer is no because the range of \eqref{eq9.9} is not
closed in $\mathcal{H}(\Omega,\partial\Omega)$.

To prove this let us suppose the contrary, i.e., the range of \eqref{eq9.9} to be
closed in $\mathcal{H}(\Omega,\partial\Omega)$. Then the restriction of
\eqref{eq9.9} to the subspace
$$
K_{\Delta}^{1/2,\varphi}(\Omega):=\bigl\{\,u\in
H^{1/2,\varphi}(\Omega):\,\Delta\,u=0\;\;\mbox{in}\;\;\Omega\,\bigr\}
$$
has a closed range in $H^{0,\varphi_{0}}(\partial\Omega)$. But, according to Theorem
\ref{th10.1} given below in Section \ref{th10.1}, this restriction establishes a
homeomorphism of $K_{\Delta}^{1/2,\varphi}(\Omega)$ onto
$H^{0,\varphi}(\partial\Omega)$. Hence, $H^{0,\varphi}(\partial\Omega)$ is a
(closed) subspace of $H^{0,\varphi_{0}}(\partial\Omega)$, so that
$H^{0,\varphi}(\partial\Omega)=H^{0,\varphi_{0}}(\partial\Omega)$. We arrive at a
contradiction if we note that $\varphi_{0}(t)/\varphi(t)\rightarrow0$ as
$t\rightarrow+\infty$ and use Theorem \ref{th5.2} ii). Thus our hypothesis is false.

Given a general elliptic boundary-value problem \eqref{eq9.1}, \eqref{eq9.2}, the
reasoning is similar. If $s=m+1/2$, $\varphi$ satisfies \eqref{eq3.6}, and
$\varphi_{0}$ is defined by \eqref{eq8.10}, then we get the bounded linear operator
\eqref{eq9.4} providing the space
$H^{s-m_{j}-1/2,\varphi}(\partial\Omega)=H^{0,\varphi}(\partial\Omega)$ is replaced
by $H^{0,\varphi_{0}}(\partial\Omega)$ for $j$ such that $m_{j}=m$. This operator
has a nonclosed range and therefore is not Fredholm.

\section{Semihomogeneous elliptic problems}\label{sec10}

As we have mentioned, the results of the previous section are not valid for $s\leq
m+1/2$. But if the boundary-value problem \eqref{eq9.1}, \eqref{eq9.2} is
semihomogeneous (i.e., $f\equiv0$ or all $g_{j}\equiv0$), it establishes a bounded
and Fredholm operator on the two-sided refined Sobolev scale, in which the number
parameter $s$ runs over the whole real axis. In this section we separately consider
the case of the homogeneous elliptic equation \eqref{eq9.1} and the case of the
homogeneous boundary conditions \eqref{eq9.2}. In what follows we focuss our
attention on analogs of Theorem \ref{th9.1} on the Fredholm property of $(L,B)$.
Counterparts of Theorems \ref{th9.2}--\ref{th9.4} can be derived from the analogs
similarly to the reasoning outlined in Section~\ref{sec9} (for details, see
\cite{05UMJ5, 06UMJ11, 06UMB4}).

\subsection{A boundary-value problem for a homogeneous elliptic
equation}\label{sec10.1}

Let us consider the regular elliptic boundary-value problem \eqref{eq9.1},
\eqref{eq9.2} provided that $f\equiv0$, namely
\begin{equation}\label{eq10.1}
Lu=0\;\;\mbox{in}\;\;\Omega,\quad
B_{j}u=g_{j}\;\;\mbox{on}\;\;\partial\Omega,\;\;j=1,\ldots,q.
\end{equation}

We connect the following linear spaces with this problem:
\begin{gather*}
K_{L}^{\infty}(\Omega):=\bigl\{\,u\in
C^{\infty}(\,\overline{\Omega}\,):\,L\,u=0\;\;\mbox{in}\;\;\Omega\,\bigr\}, \\
K_{L}^{s,\varphi}(\Omega):=\bigl\{\,u\in
H^{s,\varphi}(\Omega):\,L\,u=0\;\;\mbox{in}\;\;\Omega\,\bigr\}
\end{gather*}
for $s\in\mathbb{R}$ and $\varphi\in\mathcal{M}$. Here the equality $L\,u=0$ is
understood in the distribution theory sense. It  follows from a continuity of the
embedding $H^{s,\varphi}(\Omega)\hookrightarrow\mathcal{D}'(\Omega)$ that
$K_{L}^{s,\varphi}(\Omega)$ is a (closed) subspace in $H^{s,\varphi}(\Omega)$. We
consider $K_{L}^{s,\varphi}(\Omega)$ as a Hilbert space with respect to the inner product
in $H^{s,\varphi}(\Omega)$.

\begin{theorem}\label{th10.1}
Let $s\in\mathbb{R}$ and $\varphi\in\mathcal{M}$. Then the set
$K_{L}^{\infty}(\Omega)$ is dense in the space $K_{L}^{s,\varphi}(\Omega)$, and the
mapping
$$
u\mapsto Bu=(B_{1}u,\ldots,B_{q}u),\quad u\in K_{L}^{\infty}(\Omega),
$$
extends uniquely to a continuous linear operator
\begin{equation}\label{eq10.2}
B:\,K_{L}^{s,\varphi}(\Omega)\rightarrow
\bigoplus_{j=1}^{q}\,H^{s-m_{j}-1/2,\,\varphi}(\partial\Omega)=:
\mathcal{H}_{s,\varphi}(\partial\Omega).
\end{equation}
This operator is Fredholm. Its kernel coincides with $\mathcal{N}$, whereas its
range consists of all the vectors
$(g_{1},\ldots,g_{q})\in\mathcal{H}_{s,\varphi}(\partial\Omega)$ such that
$$
\sum_{j=1}^{q}\,(g_{j},C^{+}_{j}v)_{\partial\Omega}=0\quad\mbox{for all}\quad v\in
\mathcal{N}^{+}.
$$
The index of the operator \eqref{eq10.2} is equal to
$\dim\mathcal{N}-\dim\mathcal{G}^{+}$, with
$$
\mathcal{G}^{+}:=\bigl\{\,(C_{1}^{+}v,\ldots,C_{q}^{+}v):\,v\in
\mathcal{N}^{+}\,\bigr\},
$$
and does not depend on $s$ and $\varphi$.
\end{theorem}

Theorem \ref{th10.1} was proved in \cite[Sec. 6]{05UMJ5}. In the $s>m+1/2$ case the
theorem follows plainly from Lemma \ref{lem9.1} and Theorem \ref{th9.1}. If $s\leq
m+1/2$, then the ellipticity condition is essential for the continuity of the
operator \eqref{eq10.2}. Note that $\dim\mathcal{G}^{+}\leq\dim\mathcal{N}^{+}$, the
strict inequality being possible \cite[Theorem 13.6.15]{Hermander83}.

Theorem \ref{th10.1} can be regarded as a certain analog of the Harnack theorem on
convergence of sequences of harmonic functions (see, e.g., \cite[Ch. 11,
\S~9]{Mikhlin68}), however we use the metric in $H^{s,\varphi}(\Omega)$ instead of
the uniform metric. Here it is relevant to mention R.~Seeley's investigation
\cite{Seeley66} of the Cauchy data of solutions to a homogeneous elliptic equation
in the two-sided Sobolev scale; see also the survey \cite[Sec. 5.4~b]{Agranovich97}.

Let us outline the proof of Theorem \ref{th10.1}. For the sake of simplicity, we
suppose that both $\mathcal{N}$ and $\mathcal{N}^{+}$ are trivial. Let $s<2q$ and
$\varphi\in\mathcal{M}$. Chose an integer $r\geq1$ such that $2q(1-r)<s<2q$. We need
the following Hilbert space
\begin{gather*}
D^{s,\varphi}_{L}(\Omega):=\bigl\{u\in H^{s,\varphi}(\Omega):\,L\,u\in
L_{2}(\Omega)\bigr\},\\
(u_{1},u_{2})_{D^{s,\varphi}_{L}(\Omega)}:=(u_{1},u_{2})_{H^{s,\varphi}(\Omega)}+
(L\,u_{1},L\,u_{2})_{L_{2}(\Omega)}.
\end{gather*}
The mapping \eqref{eq9.3} extends uniquely to the homeomorphisms
\begin{gather*}
(L,B):\,D^{2q(1-r)}_{L}(\Omega)\leftrightarrow
L_{2}(\Omega)\oplus\mathcal{H}_{2q(1-r)}(\partial\Omega),\\
(L,B):\,H^{2q}(\Omega)\leftrightarrow
L_{2}(\Omega)\oplus\mathcal{H}_{2q}(\partial\Omega).
\end{gather*}
The first of them follows from the Lions--Magenes theorems \cite{LionsMagenes72}
stated below in Section \ref{sec11.1}, whereas the second is a special case of
Theorem \ref{th9.1}. (Recall that we omit $\varphi$ in the notations if
$\varphi\equiv1$.) Applying the interpolation with the function parameter $\psi$
defined by \eqref{eq3.8} with $\varepsilon:=s-2q(1-r)$ and $\delta:=2q-s$ we get
another homeomorphism
\begin{equation}\label{eq10.3}
(L,B):\,\bigl[D^{2q(1-r)}_{L}(\Omega),H^{2q}(\Omega)\bigr]_{\psi}\leftrightarrow
L_{2}(\Omega)\oplus\mathcal{H}_{s,\varphi}(\partial\Omega).
\end{equation}
Now if we prove that
$Z_{\psi}:=\bigl[D^{2q(1-r)}_{L}(\Omega),H^{2q}(\Omega)\bigr]_{\psi}$ coincides with
$D^{s,\varphi}_{L}(\Omega)$ up to equivalence of norms, then the restriction of
\eqref{eq10.3} to $K_{L}^{s,\varphi}(\Omega)$ will give the homeomorphism
\eqref{eq10.2}.

The continuous embedding $Z_{\psi}\hookrightarrow D^{s,\varphi}_{L}(\Omega)$ is
evident. The inverse can be proved by the following modification of the reasoning
used by Lions and Magenes \cite[Ch.~2, Sec. 7.2]{LionsMagenes72} for $r=1$ and power
parameter~$\psi$. In view of Theorem \ref{th9.1} we have the homeomorphism
$$
L^{r}L^{r+}+I:\,\bigl\{u\in
H^{\sigma}(\Omega):(D_{\nu}^{j-1}u)\upharpoonright\partial\Omega=0\;\;\forall\;
j=1,\ldots,r\bigr\}\leftrightarrow H^{\sigma-4qr}(\Omega)
$$
for each $\sigma\geq2qr$. Here $L^{r}$ is the $r$-th iteration of $L$, $L^{r+}$ is
the formally adjoint to $L^{r}$, and $I$ is the identity operator. We regard the
domain of $L^{r}L^{r+}+I$ as a subspace of $H^{\sigma}(\Omega)$. Consider the
bounded linear inverse operators
$$
(L^{r}L^{r+}+I)^{-1}:\,H^{\sigma}(\Omega)\rightarrow H^{\sigma+4qr}(\Omega),\quad
\sigma\geq-2qr.
$$
Set $R:=L^{r-1}L^{r+}(L^{r}L^{r+}+I)^{-1}$ and $P:=-RL+I$. Since
$$
LPu=(L^{r}L^{r+}+I)^{-1}Lu\in L_{2}(\Omega)\quad\mbox{for each}\quad u\in
H^{2q(1-r)}(\Omega),
$$
the operator $P$ maps continuously $H^{\sigma}(\Omega)\rightarrow
D^{\sigma}_{L}(\Omega)$ with $\sigma\geq2q(1-r)$. Therefore, by the interpolation,
we get the bounded operator
$$
P:\,H^{s,\varphi}(\Omega)=\bigl[H^{2q(1-r)}(\Omega),H^{2q}(\Omega)\bigr]_{\psi}\rightarrow
\bigl[D^{2q(1-r)}_{L}(\Omega),H^{2q}(\Omega)\bigr]_{\psi}=Z_{\psi}.
$$
Now, for each $u\in D^{s,\varphi}_{L}(\Omega)$, we have $u=Pu+RLu$, with $Pu\in Z_{\psi}$
and $RLu\in H^{2q}(\Omega)\subset Z_{\psi}$. So $D^{s,\varphi}_{L}(\Omega)\subseteq
Z_{\psi}$, and our reasoning is complete.

Note that in the Sobolev case of $\varphi\equiv1$ Theorem \ref{th10.1} is a
consequence of the above-mentioned Lions--Magenes theorems provided $s$ is negative
and not half-integer. If negative $s$ is half-integer, then Theorem \ref{th10.1} is
new even in the Sobolev case.

\subsection{An elliptic problem with homogeneous boundary
conditions}\label{sec10.2}

Now we consider the regular elliptic boundary-value problem \eqref{eq9.1},
\eqref{eq9.2} provided that all $g_{j}\equiv0$, namely
\begin{equation}\label{eq10.4}
Lu=f\;\;\mbox{in}\;\;\Omega,\quad
B_{j}u=0\;\;\mbox{on}\;\;\partial\Omega,\;\;j=1,\ldots,q.
\end{equation}

Let us introduce some function spaces related to the boundary-value problem
\eqref{eq10.4}. For the sake of brevity, we denote by $(\mathrm{b.c.})$ the
homogeneous boundary conditions in \eqref{eq10.4}. In addition, we denote by
$(\mathrm{b.c.})^{+}$ the homogeneous adjoint boundary conditions \eqref{eq9.6}:
$$
B^{+}_{j}v=0\;\;\mbox{on}\;\;\partial\Omega,\;\;j=1,\ldots,q.
$$
We set
\begin{gather*}
C^{\infty}(\mathrm{b.c.}):=\bigl\{u\in
C^{\infty}(\,\overline{\Omega}\,):\,B_{j}u=0\;\;\mbox{on}\;\;\partial\Omega\;\;
\forall\;\;j=1,\ldots,q\bigr\}, \\
C^{\infty}(\mathrm{b.c.})^{+}:=\bigl\{v\in
C^{\infty}(\,\overline{\Omega}\,):\,B^{+}_{j}v=0\;\;\mbox{on}\;\;\partial\Omega
\;\;\forall\;\;j=1,\ldots,q\bigr\}.
\end{gather*}

Let $s\in\mathbb{R}$ and $\varphi\in\mathcal{M}$. We introduce the separable Hilbert
spaces $H^{s,\varphi}(\mathrm{b.c.})$ and $H^{s,\varphi}(\mathrm{b.c.})^{+}$ formed by
distributions satisfying the homogeneous boundary conditions $(\mathrm{b.c.})$ and
$(\mathrm{b.c.})^{+}$ respectively.

\begin{definition}\label{def10.1}
If $s\notin\{m_{j}+1/2:j=1,\ldots,q\}$, then $H^{s,\varphi}(\mathrm{b.c.})$ is
defined to be the closure of $C^{\infty}(\mathrm{b.c.})$ in
$H^{s,\varphi,(0)}(\Omega)$, the space $H^{s,\varphi}(\mathrm{b.c.})$ being regarded
as a subspace of $H^{s,\varphi,(0)}(\Omega)$. If $s\in\{m_{j}+1/2:j=1,\ldots,q\}$,
then the space $H^{s,\varphi}(\mathrm{b.c.})$ is defined by means of the
interpolation with the power parameter $\psi(t)=t^{1/2}$:
\begin{equation}\label{eq10.5}
H^{s,\varphi}(\mathrm{b.c.}):=\bigl[H^{s-1/2,\,\varphi}(\mathrm{b.c.}),
H^{s+1/2,\,\varphi}(\mathrm{b.c.})\bigr]_{t^{1/2}}.
\end{equation}
Changing $(\mathrm{b.c.})$ for $(\mathrm{b.c.})^{+}$, and $m_{j}$ for $m_{j}^{+}$ in
the last two sentences, we have the definition of the space
$H^{s,\varphi}(\mathrm{b.c.})^{+}$.
\end{definition}

The space $C^{\infty}(\mathrm{b.c.})^{+}$ and therefore
$H^{s,\varphi}(\mathrm{b.c.})^{+}$ are independent of the choice of the system
$\{B^{+}_{j}\}$ adjoint to $\{B_{j}\}$; see, e.g., \cite[Ch.~2, Sec.
2.5]{LionsMagenes72}.

Note that the case of $s\in\{m_{j}+1/2:j=1,\ldots,q\}$ is special in the definition
of $H^{s,\varphi}(\mathrm{b.c.})$. We have to resort to the interpolation formula
\eqref{eq10.5} to get the spaces for which the main result of the subsection,
Theorem \ref{th10.3}, will be valid. In this case the norms in the spaces
$H^{s,\varphi}(\mathrm{b.c.})$ and $H^{s,\varphi,(0)}(\Omega)$ are not equivalent.
The analogous fact is true for $H^{s,\varphi}(\mathrm{b.c.})^{+}$. Providing
$\varphi\equiv1$, this was proved in \cite{Grisvard67, Seeley72} (see also
\cite[Sec. 4.3.3]{Triebel95}).

The spaces just introduced admit the following constructive description.

\begin{theorem}\label{th10.2}
Let $s\in\mathbb{R}$, $s\neq m_{j}+1/2$ for all $j=1,\ldots,q$, and
$\varphi\in\mathcal{M}$. If $s>0$, then the space $H^{s,\varphi}(\mathrm{b.c.})$
consists of the functions $u\in H^{s,\varphi}(\Omega)$ such that $B_{j}u=0$ on
$\partial\Omega$ for all indices $j=1,\ldots,q$ satisfying $s>m_{j}+1/2$. If
$s<1/2$, then $H^{s,\varphi}(\mathrm{b.c.})=H^{s,\varphi,(0)}(\Omega)$. This
proposition remains true if one changes $m_{j}$ for $m_{j}^{+}$, $(\mathrm{b.c.})$
for $(\mathrm{b.c.})^{+}$, and $B_{j}$ for $B_{j}^{+}$.
\end{theorem}

Theorem \ref{th10.2} is known in the Sobolev case of $\varphi\equiv1$ \cite[Sec.
5.5.2]{Roitberg96}. In general, we can deduce it by means of the interpolation with
a function parameter. Here we only need to treat the case where
$m_{k}+1/2<s<m_{k+1}+1/2$ for some $k=1,\ldots,q$, with $m_{1}<m_{2}<\ldots<m_{q}$
and $m_{q+1}:=\infty$. Chose $\varepsilon>0$ such that
$m_{k}+1/2<s\mp\varepsilon<m_{k+1}+1/2$. Then the space
$H^{s\mp\varepsilon}(\mathrm{b.c.})$ consists of the functions $u\in
H^{s\mp\varepsilon}(\Omega)$ satisfying the condition $B_{j}u=0$ on $\partial\Omega$
for all $j=1,\ldots,k$. So there exists a projector $P_{k}$ of
$H^{s\mp\varepsilon}(\Omega)$ onto $H^{s\mp\varepsilon}(\mathrm{b.c.})$; it is
constructed in \cite[the proof of Lemma 5.4.4]{Triebel95}. Hence, by Proposition
\ref{prop8.1} and Theorem \ref{th8.1} with $\varepsilon=\delta$, we get that
$Y_{\psi}:=[H^{s-\varepsilon}(\mathrm{b.c.}),H^{s+\varepsilon}(\mathrm{b.c.})]_{\psi}$
is the subspace $H^{s,\varphi}(\Omega)\cap H^{s-\varepsilon}(\mathrm{b.c.)}$ of
$H^{s,\varphi}(\Omega)$. Now since $C^{\infty}(\mathrm{b.c.})$ is dense in
$Y_{\psi}$, we have
\begin{equation}\label{eq10.6}
\bigl[H^{s-\varepsilon}(\mathrm{b.c.}),H^{s+\varepsilon}(\mathrm{b.c.})\bigl]_{\psi}=
H^{s,\varphi}(\mathrm{b.c.})
\end{equation}
with equivalence of norms, so that $H^{s,\varphi}(\mathrm{b.c.})$ admits the
description stated in Theorem \ref{th10.1}.

The following theorem is about the Fredholm property of the boundary-value problem
\eqref{eq10.4} in the two-sided refined Sobolev scale.

\begin{theorem}\label{th10.3}
Let $s\in\mathbb{R}$ and $\varphi\in\mathcal{M}$. Then the mapping $u\mapsto Lu$,
with $u\in C^{\infty}(\mathrm{b.c.})$, extends uniquely to a continuous linear
operator
\begin{equation}\label{eq10.7}
L:H^{s,\varphi}(\mathrm{b.c.})\rightarrow (H^{2q-s,1/\varphi}(\mathrm{b.c.})^{+})'.
\end{equation}
Here $Lu$ is interpreted as the functional $(Lu,\,\cdot\,)_{\Omega}$, and
$(H^{2q-s,\,1/\varphi}(\mathrm{b.c.})^{+})'$ denotes the antidual space to
$H^{2q-s,1/\varphi}(\mathrm{b.c.})^{+}$ with respect to the inner product in
$L_{2}(\Omega)$. The operator \eqref{eq10.7} is Fredholm. Its kernel coincides with
$\mathcal{N}$, whereas its range consists of all the functionals
$f\in(H^{2q-s,\,1/\varphi}(\mathrm{b.c.})^{+})'$ such that $(f,v)_{\Omega}=0$ for
all $v\in\mathcal{N}^{+}$. The index of \eqref{eq10.7} is
$\dim\mathcal{N}-\dim\mathcal{N}^{+}$ and does not depend on $s$ and~$\varphi$.
\end{theorem}

For the Sobolev scale, where $\varphi\equiv1$, this theorem was proved by
Yu.~M.~Berezansky, S.G.~Krein, and Ya.A.~Roitberg (\cite{BerezanskyKreinRoitberg63}
and \cite[Ch. III, \S~6, Sec. 10]{Berezansky68}) in the case of integral $s$ and by
Roitberg \cite[Sec. 5.5.2]{Roitberg96} for all real $s$; see also the textbook
\cite[Ch. XVI, \S~1]{BerezanskySheftelUs96b} and the survey \cite[Sec.
7.9~c]{Agranovich97}. They formulated the theorem in an equivalent form of a
homeomorphism theorem. Note that if $s\leq m+1/2$, then the ellipticity condition is
essential for the continuity of the operator \eqref{eq10.7}.

For arbitrary $\varphi\in\mathcal{M}$, Theorem \ref{th10.3} follows from the Sobolev
case by Proposition \ref{prop6.1} if we apply the interpolation formulas
\eqref{eq10.6}, \eqref{eq10.5} and their counterparts for
$H^{2q-s,1/\varphi}(\mathrm{b.c.})^{+}$. First we should use \eqref{eq10.6} for
$s\notin\{j-1/2:j=1,\ldots,2q\}$ and a sufficiently small $\varepsilon>0$, then
should apply \eqref{eq10.5} for the rest of $s$. Moreover, we have to resort to the
interpolation duality formula $[X_{1}',X_{0}']_{\psi}=[X_{0},X_{1}]_{\chi}'$, where
$X:=[X_{0},X_{1}]$ is an admissible couple of Hilbert spaces and
$\chi(t):=t/\psi(t)$ for $t>0$. The formula follows directly from the definition of
$X_{\psi}$; see, e.g., \cite[Sec. 2.4]{08MFAT1}.

\subsection{On a connection between nonhomogeneous and semihomogeneous
elliptic problems}\label{sec10.3}

Here, for the sake of simplicity, we suppose that
$\mathcal{N}=\mathcal{N}^{+}=\{0\}$. Let $s>m+1/2$ and $\varphi\in\mathcal{M}$. It
follows from Theorems \ref{th9.1} and \ref{th10.2} that the space
$H^{s,\varphi}(\Omega)$ is the direct sum of the subspaces
$K^{s,\varphi}_{L}(\Omega)$ and $H^{s,\varphi}(\mathrm{b.c.})$. Therefore Theorem
\ref{th9.1} are equivalent to Theorems \ref{th10.1} and \ref{th10.3} taken together;
note that the antidual space $(H^{2q-s,1/\varphi}(\mathrm{b.c.})^{+})'$ coincides
with $H^{s-2q,\varphi}(\Omega)$. Thus the nonhomogeneous problem \eqref{eq9.1},
\eqref{eq9.2} can be reduced immediately to the semihomogeneous problems
\eqref{eq10.1} and \eqref{eq10.4} provided $s>m+1/2$.

This reduction fails for $s<m+1/2$. Indeed, if $0\leq s<m+1/2$, then the operator
$(L,B)$ cannot be well defined on $K_{L}^{s,\varphi}(\Omega)\cup
H^{s,\varphi}(\mathrm{b.c.})$ because $K_{L}^{s,\varphi}(\Omega)\cap
H^{s,\varphi}(\mathrm{b.c.})\neq\varnothing$. This inequality follows from Theorems
\ref{th10.1} and \ref{th10.2} if we note that the boundary-value problem
\eqref{eq10.1}, with $g_{q}\equiv1$ and $g_{j}\equiv0$ for $j<q$, has a nonzero
solution $u\in K_{L}^{\infty}(\Omega)$ belonging to $H^{s,\varphi}(\mathrm{b.c.})$.
Here we may suppose that $m_{q}=m$.

So much the more, the above reduction is impossible for negative $s$. Note if
$s<-1/2$, then solutions to the semihomogeneous problems pertain to the spaces of
distributions of the different nature. Namely, the solutions to the problem
\eqref{eq10.1} belong to $K_{L}^{s,\varphi}(\Omega)\subset H^{s,\varphi}(\Omega)$
and are distributions given in the open domain $\Omega$, whereas the solutions to
the problem \eqref{eq10.4} belong to $H^{s,\varphi}(\mathrm{b.c.})\subset
H^{s,\varphi}_{\overline{\Omega}}(\mathbb{R}^{n})$ and are distributions supported
on the closed domain $\overline{\Omega}$.

The same conclusions are valid in general, for nontrivial $\mathcal{N}$ and/or
$\mathcal{N}^{+}$.

\medskip

\section{Generic theorems for elliptic problems in two-sided scales}\label{sec11}

Let us return to the nonhomogeneous regular elliptic boundary-value problem
\eqref{eq9.1}, \eqref{eq9.2}. We aim to prove analogs of Theorem \ref{th9.1} for
\emph{arbitrary} real $s$. To get the bounded operator $(L,B)$ for such $s$ we have
to chose another space instead of $H^{s,\varphi}(\Omega)$ as a domain of the
operator. There are known two essentially different ways to construct the domain.
They were suggested by Ya.A.~Roitberg \cite{Roitberg64, Roitberg65, Roitberg96} and
J.-L.~Lions, E.~Magenes \cite{LionsMagenes62, LionsMagenes63, LionsMagenes72,
Magenes65} in the Sobolev case. These ways lead to different kinds of theorems on
the Fredholm property of $(L,B)$; we name them generic and individual theorems. In
generic theorems, the domain of $(L,B)$ does not depend on the coefficients of the
elliptic expression $L$ and is generic for all boundary-value problems of the same
order. Note that Theorem \ref{th9.1} is generic. In individual theorems, the domain
depends on coefficients of $L$, even on the coefficients of lover order derivatives.

In this section we realize Roitberg's approach with regard to the refined Sobolev
scale; namely, we modify this scale by Roitberg and prove a generic theorem about
the Fredholm property of $(L,B)$ on the two-sided modified scale. The results of the
section were obtained by the authors in \cite{08UMJ4}. Lions and Magenes's approach
led to individual theorems will be considered below in Section \ref{sec12}.

\subsection{The modification of the refined Sobolev scale}\label{sec11.1}

Let $s\in\mathbb{R}$, $\varphi\in\mathcal{M}$, and integer $r>0$. We set
$E_{r}:=\{k-1/2:k=1,\ldots,r\}$. Note that $D_{\nu}:=i\,\partial/\partial\nu$, where
$\nu$ is the field of unit vectors of inner normals to $\partial\Omega$. Let us
define the separable Hilbert spaces $H^{s,\varphi,(r)}(\Omega)$, which form the
modified scale.

\begin{definition}\label{def11.1}
If $s\in\mathbb{R}\setminus E_{r}$, then the space $H^{s,\varphi,(r)}(\Omega)$ is
defined to be the completion of $C^{\infty}(\,\overline{\Omega}\,)$ with respect to
the Hilbert norm
\begin{equation}\label{eq11.1}
\|u\|_{H^{s,\varphi,(r)}(\Omega)}:= \biggl(\,\|u\|_{H^{s,\varphi,(0)}(\Omega)}^{2}+
\sum_{k=1}^{r}\;\bigl\|(D_{\nu}^{k-1}u)\upharpoonright\partial\Omega\,\bigr\|
_{H^{s-k+1/2,\varphi}(\partial\Omega)}^{2}\,\biggr)^{1/2}.
\end{equation}
If $s\in E_{r}$, then the space $H^{s,\varphi,(r)}(\Omega)$ is defined by means of
the interpolation with the power parameter $\psi(t)=t^{1/2}$, namely
\begin{equation}\label{eq11.2}
H^{s,\varphi,(r)}(\Omega):=\bigl[\,H^{s-1/2,\varphi,(r)}(\Omega),
H^{s+1/2,\varphi,(r)}(\Omega)\,\bigr]_{t^{1/2}}.
\end{equation}
\end{definition}

In the Sobolev case of $\varphi\equiv1$ the space $H^{s,\varphi,(r)}(\Omega)$ was
introduced and investigated by Ya.A.~Roitberg; see \cite{Roitberg64, Roitberg65} and
\cite[Ch.~2]{Roitberg96}. As usual, we put $H^{s,(r)}(\Omega):=H^{s,1,(r)}(\Omega)$.

Note that the case of $s\in E_{r}$ is special in Definition \ref{def11.1} because
the norm in $H^{s,\varphi,(r)}(\Omega)$ is defined by the interpolation formula
\eqref{eq11.2} instead of \eqref{eq11.1}. These formulas give nonequivalent norms.
As in Subsection \ref{sec10.2}, we have to resort to the interpolation in the
mentioned case to get the spaces for which the main result of this section, Theorem
\ref{th11.2}, will be true.

\begin{definition}\label{def11.2}
The class of Hilbert spaces
\begin{equation}\label{eq11.3}
\{H^{s,\varphi,(r)}(\Omega):\,s\in\mathbb{R},\,\varphi\in\mathcal{M}\}
\end{equation}
is called the refined Sobolev scale modified by Roitberg. The number $r$ is called
the order of the modification.
\end{definition}

The scale \eqref{eq11.3} is found fruitful in the theory of boundary-value problems
because the trace mapping \eqref{eq8.7} extends uniquely to an operator $R_{r}$
mapping continuously
$H^{s,\varphi,(r)}(\Omega)\rightarrow\mathcal{H}^{r}_{s,\varphi}(\partial\Omega)$
for all real $s$. It is useful to compare this fact with Theorem \ref{th8.4}, in
which the condition $s>r-1/2$  cannot be neglected. Note that
\begin{equation}\label{eq11.4}
H^{s,\varphi,(r)}(\Omega)=H^{s,\varphi}(\Omega)\quad\mbox{if}\quad s>r-1/2
\end{equation}
because the spaces in \eqref{eq11.4} are completions of
$C^{\infty}(\,\overline{\Omega}\,)$ with equivalence norms due to
Theorem~\ref{th8.4}.

The spaces $H^{s,\varphi,(r)}(\Omega)$ admit the following isometric representation.
We denote by $K_{s,\varphi,(r)}(\Omega,\partial\Omega)$ the linear space of all
vectors
\begin{equation}\label{eq11.5}
(u_{0},u_{1},\ldots,u_{r})\in H^{s,\varphi,(0)}(\Omega)
\oplus\bigoplus_{k=1}^{r}\,H^{s-k+1/2,\,\varphi}(\partial\Omega)=:
\Pi_{s,\varphi,(r)}(\Omega,\partial\Omega)
\end{equation}
such that $u_{k}=(D_{\nu}^{k-1}u_{0})\!\upharpoonright\!\partial\Omega$ for each
integer $k=1,\ldots r$ satisfying $s>k-1/2$. By Theorem \ref{th8.4} we may regard
$K_{s,\varphi,(r)}(\Omega,\partial\Omega)$ as a subspace of
$\Pi_{s,\varphi,(r)}(\Omega,\partial\Omega)$.

\begin{theorem}\label{th11.1}
The mapping
$$
T_{r}:u\mapsto\bigl(\,u,u\!\upharpoonright\!\partial\Omega,\ldots,
(D_{\nu}^{r-1}u)\!\upharpoonright\!\partial\Omega\,\bigr),\quad u\in
C^{\infty}(\,\overline{\Omega}\,),
$$
extends uniquely to a continuous linear operator
\begin{equation}\label{eq11.6}
T_{r}:\,H^{s,\varphi,(r)}(\Omega)\rightarrow
K_{s,\varphi,(r)}(\Omega,\partial\Omega)
\end{equation}
for all $s\in\mathbb{R}$ and $\varphi\in\mathcal{M}$. This operator is injective.
Moreover, if $s\notin E_{r}$, then \eqref{eq11.6} is an isometric isomorphism.
\end{theorem}

We need only to argue that \eqref{eq11.6} is surjective if $s\notin E_{r}$. For
$\varphi\equiv1$ this property is proved by Ya.A.~Roitberg; see, e.g.,
\cite[Sec.~2.2]{Roitberg96}. In general, the proof is quite similar provided we
apply Theorem \ref{th8.4} and \eqref{eq8.8}.

Note that we have the following dense compact embeddings in the modified scale
\eqref{eq11.3}:
\begin{equation}\label{eq11.7}
H^{s_{1},\varphi_{1},(r)}(\Omega)\hookrightarrow
H^{s,\varphi,(r)}(\Omega),\quad-\infty<s<s_{1}<\infty\;\;\mbox{and}
\;\;\varphi,\varphi_{1}\in\mathcal{M}.
\end{equation}
They results from \eqref{eq8.16} and Theorem \ref{th5.2} (i) by Theorem \ref{th11.1}
and are understood as embeddings of spaces which are completions of the same set,
$C^{\infty}(\,\overline{\Omega}\,)$, with different norms. Suppose that $s=s_{1}$,
then the continuous embedding \eqref{eq11.7} holds if and only if
$\varphi/\varphi_{1}$ is bounded in a neighbourhood of~$+\infty$; the embedding is
compact if and only if $\varphi(t)/\varphi_{1}(t)\rightarrow0$ as
$t\rightarrow\nobreak+\infty$. This follows from the relevant properties of the
refined Sobolev scales over $\Omega$ and $\partial\Omega$.

\subsection{Roitberg's type generic theorem}\label{sec11.2}

The main result of the section is the following generic theorem about properties of
the operator $(L,B)$ on the two-sided scale \eqref{eq11.3} with $r=2q$.

\begin{theorem}\label{th11.2}
Let $s\in\mathbb{R}$ and $\varphi\in\mathcal{M}$. The mapping \eqref{eq9.3} extends
uniquely to a continuous linear operator
\begin{gather}\label{eq11.8}
(L,B):\,H^{s,\varphi,(2q)}(\Omega)\rightarrow
H^{s-2q,\varphi,(0)}(\Omega)\oplus\bigoplus_{j=1}^{q}\,
H^{s-m_{j}-1/2,\varphi}(\partial\Omega)\\
=:\mathcal{H}_{s,\varphi,(0)}(\Omega,\partial\Omega).\notag
\end{gather}
This operator is Fredholm. Its kernel coincides with $\mathcal{N}$, and its range
consists of all the vectors
$(f,g_{1},\ldots,g_{q})\in\mathcal{H}_{s,\varphi,(0)}(\Omega,\partial\Omega)$ that
satisfy \eqref{eq9.8}. The index of \eqref{eq11.8} is
$\dim\mathcal{N}-\dim\mathcal{N}^{+}$ and does not depend on $s$ and $\varphi$.
\end{theorem}

Note that Theorem \ref{th11.2} is generic because the domain of the operator
\eqref{eq11.8}, the space $H^{s,\varphi,(2q)}(\Omega)$, is independent of $L$ due to
Definition \ref{def11.1}. If $s>2q-1/2$, then generic Theorems \ref{th9.1} and
\ref{th11.2} are tantamount in view of \eqref{eq11.4} and \eqref{eq8.17}.

For the modified Sobolev scale (the $\varphi\equiv1$ case) Theorem \ref{th11.2} was
proved by Ya.A.~Roi\-tberg \cite{Roitberg64, Roitberg65}, \cite[Ch.~4 and Sec.
5.3]{Roitberg96}; see also the monograph \cite[Ch.~3, Sec.~6, Theorem
6.9]{Berezansky68}, the handbook \cite[Ch.~III, \S~6, Sec. 5]{FunctionalAnalysis72},
and the survey \cite[Sec. 7.9]{Agranovich97}.

For arbitrary $\varphi\in\mathcal{M}$ we can deduce Theorem \ref{th11.2} from the
$\varphi\equiv1$ case with the help of the interpolation in the following way. First
assume that $s\notin E_{2q}$ and let $\varepsilon>0$. We have the Fredholm bounded
operators on the modified Sobolev scale
\begin{equation}\label{eq11.9}
(L,B):\,H^{s\mp\varepsilon,(2q)}(\Omega)\rightarrow
\mathcal{H}_{s\mp\varepsilon,(0)}(\Omega,\partial\Omega).
\end{equation}
They possess the common kernel $\mathcal{N}$ and the common index
$\varkappa:=\dim\mathcal{N}-\dim\mathcal{N}^{+}$. Applying the interpolation with
the function parameter $\psi$ defined by \eqref{eq3.8} for $\varepsilon=\delta$, we
get by Proposition \ref{prop6.1} and Theorems \ref{th5.3}, \ref{th8.6} that
\eqref{eq11.9} implies the boundedness and the Fredholm property of the operator
$$
(L,B):\,\bigl[H^{s-\varepsilon,(2q)}(\Omega),H^{s+\varepsilon,(2q)}(\Omega)\bigr]_{\psi}
\rightarrow\mathcal{H}_{s,\varphi,(0)}(\Omega,\partial\Omega).
$$
It remains to prove the interpolation formula
\begin{equation}\label{eq11.10}
\bigl[H^{s-\varepsilon,(2q)}(\Omega),H^{s+\varepsilon,(2q)}(\Omega)\bigr]_{\psi}=
H^{s,\varphi,(2q)}(\Omega),
\end{equation}
where the equality of spaces is up to equivalence of norms.

Let an index $p$ be such that $s\in\alpha_{p}$, where $\alpha_{0}:=(-\infty,1/2)$,
$\alpha_{k}:=(k-1/2,\,k+1/2)$ with $k=1,\ldots,2q-1$, and
$\alpha_{2q}:=(2q-1/2,\infty)$. We chose $\varepsilon>0$ satisfying
$s\mp\varepsilon\in\alpha_{p}$. By Theorem \ref{th11.1}, the mapping
$$
T_{2q,p}:\,u\mapsto
\left(\,u,\,\{(D_{\nu}^{k-1}u)\!\upharpoonright\!\partial\Omega:\,p+1\leq k\leq
2q\}\,\right)
$$
establishes the homeomorphisms
\begin{gather} \label{eq11.11}
T_{2q,p}:\,H^{s,\varphi,(2q)}(\Omega)\leftrightarrow H^{s,\varphi,(0)}(\Omega)
\oplus\bigoplus_{p+1\leq k\leq2q}
H^{s-k+1/2,\varphi}(\partial\Omega)=:K_{s,\varphi,(2q)}^{p}(\Omega,\partial\Omega),\\
T_{2q,p}:\,H^{s\mp\varepsilon,(2q)}(\Omega)\leftrightarrow
K_{s\mp\varepsilon,(2q)}^{p}(\Omega,\partial\Omega). \label{eq11.12}
\end{gather}
Applying the interpolation with $\psi$, we deduce another homeomorphism from
\eqref{eq11.12}:
\begin{equation}\label{eq11.13}
T_{2q,p}:\,\bigl[H^{s-\varepsilon,(2q)}(\Omega),H^{s+\varepsilon,(2q)}(\Omega)\bigr]_{\psi}
\leftrightarrow K_{s,\varphi,(2q)}^{p}(\Omega,\partial\Omega).
\end{equation}
Now \eqref{eq11.11} and \eqref{eq11.13} imply the required formula \eqref{eq11.10}.

In the remaining case of $s\in E_{2q}$, we deduce Theorem \ref{th11.2} from the
$s\notin E_{2q}$ case by the interpolation with the power parameter
$\psi(t)=t^{1/2}$ if we apply \eqref{eq11.2} and the counterparts of Theorem
\ref{th3.5} for the refined Sobolev scales over $\Omega$ and $\partial\Omega$.

Note that the continuity of the operator \eqref{eq11.8} holds without the assumption
about the regular ellipticity of the boundary-value problem \eqref{eq9.1},
\eqref{eq9.2}.

\subsection{Roitberg's interpretation of generalized solutions}\label{sec11.3}

Using Theorem 10.1, we can give the following interpretation of a solution $u\in
H^{s,\varphi,(2q)}(\Omega)$ to the boun\-dary-value problem \eqref{eq9.1},
\eqref{eq9.2} in the framework of the distribution theory.

Let us write down the differential expressions $L$ and $B_{j}$ in a neighbourhood of
$\partial\Omega$ in the form
\begin{equation}\label{eq11.14}
L=\sum_{k=0}^{2q}\;L_{k}\,D_{\nu}^{k},\quad
B_{j}=\sum_{k=0}^{m_{j}}\;B_{j,k}\,D_{\nu}^{k}.
\end{equation}
Here $L_{k}$ and $B_{j,k}$ are certain tangent differential expression. Integrating
by parts, we arrive at the (special) Green formula
$$
(Lu,v)_{\Omega}=(u,L^{+}v)_{\Omega}-
i\sum_{k=1}^{2q}\;(D_{\nu}^{k-1}u,L^{(k)}v)_{\partial\Omega},\quad u,v\in
C^{\infty}(\,\overline{\Omega}\,).
$$
Here $L^{(k)}:=\sum_{r=k}^{2q}D_{\nu}^{r-k}L_{r}^{+}$, where $L_{r}^{+}$ is the
tangent differential expression formally adjoint to $L_{r}$. By passing to the limit
and using the notation
\begin{equation}\label{eq11.15}
(u_{0},u_{1},\ldots,u_{2q}):=T_{2q}u\in K_{s,\varphi,(2q)}(\Omega,\partial\Omega),
\end{equation}
we get the next equality for $u\in H^{s,\varphi,(2q)}(\Omega)$:
\begin{equation}\label{eq11.16}
(Lu,v)_{\Omega}=(u_{0},L^{+}v)_{\Omega}-
i\sum_{k=1}^{2q}\;(u_{k},L^{(k)}v)_{\partial\Omega},\quad v\in
C^{\infty}(\,\overline{\Omega}\,).
\end{equation}

Now it follows from \eqref{eq11.14} and \eqref{eq11.16} that the element $u\in
H^{s,\varphi,(2q)}(\Omega)$ is a solution to the boundary-value problem
\eqref{eq9.1}, \eqref{eq9.2} with $f\in H^{s-2q,\varphi,(0)}(\Omega)$, $g_{i}\in
H^{s-m_{j}-1/2,\,\varphi}(\partial\Omega)$ if and only if
\begin{gather} \label{eq11.17}
(u_{0},L^{+}v)_{\Omega}- i\sum_{k=1}^{2q}\;(u_{k},L^{(k)}v)_{\partial\Omega}=
(f,v)_{\Omega}\quad\mbox{for all}\quad v\in
C^{\infty}(\,\overline{\Omega}\,),\\
\sum_{k=0}^{m_{j}}\;B_{j,k}\,u_{k+1}=g_{j}\;\;\mbox{on}\;\;\partial\Omega,\quad
j=1,\ldots,q. \label{eq11.18}
\end{gather}

Note that these equalities have meaning for arbitrary distributions
\begin{gather} \label{eq11.19}
u_{0}\in\mathcal{D}'(\mathbb{R}^{n}),\;\;
\mathrm{supp}\,u_{0}\subseteq\overline{\Omega},\quad
u_{1},\ldots,u_{2q}\in\mathcal{D}'(\Gamma), \\
f\in\mathcal{D}'(\mathbb{R}^{n}),\;\;
\mathrm{supp}\,f\subseteq\overline{\Omega},\quad
g_{1},\ldots,g_{q}\in\mathcal{D}'(\Gamma). \label{eq11.20}
\end{gather}
Therefore it is useful to introduce the following notion.

\begin{definition}\label{def11.3}
Suppose that \eqref{eq11.19} and \eqref{eq11.20} are fulfilled. Then the vector
$u:=(u_{0},u_{1},\ldots,u_{2q})$ is called a generalized solution in Roitberg's
sense to the boundary-value problem \eqref{eq9.1}, \eqref{eq9.2} if the conditions
\eqref{eq11.17} and \eqref{eq11.18} are valid.
\end{definition}

This interpretation of a generalized solution is suggested by Roitberg; see, e.g,
his monograph \cite[Sec. 2.4]{Roitberg96}.

Thus, Theorem \ref{th11.2} can be regarded as a statement about the solvability of
the boundary-value problem \eqref{eq9.1}, \eqref{eq9.2} in the class of generalized
solutions in Roitberg's sense provided that we identify solutions $u\in
H^{s,\varphi,(2q)}(\Omega)$ with vectors \eqref{eq11.15}.

Roitberg's interpretation of a generalized solution and the relevant Theorem
\ref{th11.2} have been found fruitful in the theory of elliptic boundary-value
problems. Analogs of this theorem were proved by Roitberg for nonregular elliptic
boundary-value problems and for general elliptic systems of differential equations,
the modified scale of the $L_{p}$-type Sobolev spaces with $1<p<\infty$ being used.
In the literature \cite{FunctionalAnalysis72, Roitberg96, Roitberg99}, Theorem
\ref{th11.2} and its analogs are known as theorems on a complete collection of
homeomorphisms. They have various applications; among them are the theorems on an
increase in smoothness of solutions up to the boundary, application to the
investigation of Green functions of elliptic boundary-value problems, applications
to elliptic problems with power singularities, to the transmission problem, the
Odhnoff problem, and others. The investigations of Ya.A.~Roitberg, Z.G.~Sheftel' and
their disciples into this subject were summed up in Roitberg's monographs
\cite{Roitberg96}.

Note that, in the most general form, the theorem on a complete collection of
homeomorphisms was proved by A.~Kozhevnikov \cite{Kozhevnikov01} for general
elliptic pseudodifferential boundary-value problems. Analogs of Theorem \ref{th11.2}
were obtained in \cite{94UMJ12, 94Dop12} for some non-Sobolev Banach spaces
parametrized by collections of numbers; the case of a scalar elliptic equation was
treated therein. We also remark applications of the concept of a generalized
solution and relevant modified two-sided scale in the theory of elliptic
boundary-value problems in nonsmooth domains \cite{KozlovMazyaRossmann97} and in the
theory of parabolic \cite{EidelmanZhitarashu98} and hyperbolic \cite{Roitberg99}
equations.

\section{Individual theorems for elliptic problems}\label{sec12}

In this section, we generalize J.-L.~Lions and E.~Magenes's method
\cite{LionsMagenes62, LionsMagenes63, LionsMagenes72, Magenes65} for constructing of
the domain of the operator $(L,B)$. We prove new theorems on the Fredholm property
of the operator on scales of Sobolev inner product spaces and some H\"ormander
spaces. These theorems has an individual character because the domain of $(L,B)$
depends on coefficients of elliptic expression $L$, as distinguished from generic
Theorems \ref{th9.1} and \ref{th11.2}. Moreover, in the individual theorems the
operator $(L,B)$ acts on the spaces consisting of distributions given in the domain
$\Omega$, so that we do not need to modify the refined Sobolev scale as this was
done for Theorem \ref{th11.2}.

The section is organized in the following manner. First, for the sake of the
reader's convenience, we recall Lions and Magenes's theorems about elliptic
boun\-da\-ry-value problems. Then we prove a certain general form of the
Lions--Magenes theorems; we call it the key theorem. Namely, we find a general
condition on the space of right-hand sides of the elliptic equation $Lu=f$ under
which the operator $(L,B)$ is bounded and Fredholm on the corresponding pairs of
Sobolev inner product spaces of negative order. Extensive classes of the spaces
satisfying this condition will be constructed; they contain the spaces used by Lions
and Magenes and many others spaces. These results motivate statements and proofs of
individual theorems on the Fredholm property of the operator $(L,B)$ on some Hilbert
H\"ormander spaces.

\subsection{The Lions--Magenes theorems}\label{sec12.1}

As we have mentioned in Remark \ref{rem8.1}, J.-L.~Lions and E.~Magenes used a
definition of the Sobolev space of negative order $s$ over $\Omega$ which is
different from our Definition \ref{def8.2} for $\varphi\equiv1$. Namely, they
defined this space as the dual of $H^{-s}_{0}(\Omega)$ with respect to the inner
product in $L_{2}(\Omega)$. We use this definition throughout Section \ref{sec12}.
To distinguish the Sobolev spaces $H^{s}(\Omega)$ introduced above by Definition
\ref{def8.2} from ones used here, we resort to the somewhat different notation
$\mathrm{H}^{s}(\Omega)$, where the letter H is not slanted.

Thus we put
$$
\mathrm{H}^{s}(\Omega):=
\begin{cases}
\;H^{s}(\Omega)&\;\text{for}\;s\geq0, \\
\;(H^{-s}_{0}(\Omega))'&\;\text{for}\;s<0.
\end{cases}
$$
Here $(H^{-s}_{0}(\Omega))'$ denotes the Hilbert space antidual to $H^{-s}_{0}(\Omega)$
with respect to the inner product in $L_{2}(\Omega)$.

The antilinear continuous functionals from $\mathrm{H}^{s}(\Omega)$ with $s<0$ are
defined uniquely by their values on the functions in $C^{\infty}_{0}(\Omega)$.
Therefore it is reasonable to identify these functionals with distributions given in
$\Omega$. In so doing, we have \cite[Ch.~1, Remark 12.5]{LionsMagenes72}
\begin{equation}\label{eq12.1}
\mathrm{H}^{s}(\Omega)=H^{s}_{\overline{\Omega}}(\mathbb{R}^{n})/
H^{s}_{\partial\Omega}(\mathbb{R}^{n})=\bigl\{w\!\upharpoonright\!\Omega:\,w\in
H^{s}_{\overline{\Omega}}(\mathbb{R}^{n})\bigr\}\quad\mbox{for}\quad s<0.
\end{equation}

It is remarkable that the spaces $\mathrm{H}^{s}(\Omega)$ and $H^{s}(\Omega)$, with
$s<0$, coincide up to equivalence of norms provided $s+1/2\notin\mathbb{Z}$; see,
e.g., \cite[Sec. 4.8.2]{Triebel95}. If $s$ is half-integer, then
$\mathrm{H}^{s}(\Omega)$ is narrower than $H^{s}(\Omega)$. Note also that
\begin{equation}\label{eq12.2}
-1/2\leq s<0\;\Rightarrow\;\mathrm{H}^{s}(\Omega)=H^{s,(0)}(\Omega)\;\;\mbox{with
equality of norms}.
\end{equation}
This fact follows, by the duality, from the equality
$H^{-s}_{0}(\Omega)=H^{-s}(\Omega)$; see, e.g., \cite[Sec. 4.7.1]{Triebel95}.

Lions and Magenes consider the operator
\begin{equation}\label{eq12.3}
(L,B):\,D^{\sigma+2q}_{L,X}(\Omega)\rightarrow
X^{\sigma}(\Omega)\oplus\bigoplus_{j=1}^{q}\,H^{\sigma+2q-m_{j}-1/2}(\partial\Omega)
=:\mathbf{X}_{\sigma}(\Omega,\partial\Omega),
\end{equation}
with $\sigma\in\mathbb{R}$. Here $X^{\sigma}(\Omega)$ is a certain Hilbert space
consisting of distributions in $\Omega$ and embedded continuously in
$\mathcal{D}'(\Omega)$. The domain of the operator \eqref{eq12.3} is the Hilbert
space
$$
D^{\sigma+2q}_{L,X}(\Omega):=\bigl\{u\in \mathrm{H}^{\sigma+2q}(\Omega):\, Lu\in
X^{\sigma}(\Omega)\bigr\}
$$
endowed with the graph inner product
$$
(u_{1},u_{2})_{D^{\sigma+2q}_{L,X}(\Omega)}:=
(u_{1},u_{2})_{\mathrm{H}^{\sigma+2q}(\Omega)}+(Lu_{1},Lu_{2})_{X^{\sigma}(\Omega)}.
$$

In the case where $s:=\sigma+2q>m+1/2$ we may set
$X^{\sigma}(\Omega):=H^{\sigma}(\Omega)$ that leads us to Theorem \ref{th9.1} for
$\varphi\equiv1$. But in the case where $s\leq m+1/2$ we cannot do so if we want to
have the well-defined operator \eqref{eq12.3}. The space $X^{\sigma}(\Omega)$ must
be narrower than $H^{\sigma}(\Omega)$.

Lions and Magenes found some important spaces $X^{\sigma}(\Omega)$ with $\sigma<0$
such that the operator \eqref{eq12.3} is bounded and Fredholm; see
\cite{LionsMagenes62, LionsMagenes63} and \cite[Ch.~2, Sec. 6.3]{LionsMagenes72}. We
state their results in the form of two individual theorems on elliptic
boundary-value problems.

\begin{theorem}[the first Lions--Magenes theorem \cite{LionsMagenes62, LionsMagenes63}]
\label{thLM1} Let $\sigma<0$ and $X^{\sigma}(\Omega):=L_{2}(\Omega)$. Then the
mapping \eqref{eq9.3} extends uniquely to the continuous linear operator
\eqref{eq12.3}. This operator is Fredholm. Its kernel coincides with $\mathcal{N}$,
and its range consists of all the vectors
$(f,g_{1},\ldots,g_{q})\in\mathbf{X}_{\sigma}(\Omega,\partial\Omega)$ satisfying
\eqref{eq9.8}. The index of \eqref{eq12.3} is $\dim\mathcal{N}-\dim\mathcal{N}^{+}$
and does not depend on $\sigma$.
\end{theorem}

\begin{remark}
Here, the $\sigma=-2q$ case is important in the spectral theory of elliptic
operators with general boundary conditions \cite{Grubb68, Grubb96, Mikhailets82,
Mikhailets89}; see also the survey \cite[Sec. 7.7 and 9.6]{Agranovich97}. Then the
space $D^{0}_{A,L_{2}}(\Omega)=\{u\in L_{2}(\Omega):Au\in L_{2}(\Omega)\}$ is the
domain of the maximal operator $A_{\mathrm{max}}$ corresponding to the differential
expression~$A$. Even when all coefficients of $A$ are constant, this space depends
essentially on each of them \cite[Sec. 3.1, Theorem 3.1]{Hermander55}.
\end{remark}

To formulate the second Lions--Magenes theorem, we need the next weighted space
$$
\varrho\mathrm{H}^{\sigma}(\Omega):=\{f=\varrho
v:\,v\in\mathrm{H}^{\sigma}(\Omega)\,\},\quad (f_{1},f_{2})_{\varrho
\mathrm{H}^{\sigma}(\Omega)}:=
(\varrho^{-1}f_{1},\varrho^{-1}f_{2})_{\mathrm{H}^{\sigma}(\Omega)},
$$
with $\sigma<0$ and a positive function $\varrho\in C^{\infty}(\Omega)$. The space
$\varrho\mathrm{H}^{\sigma}(\Omega)$ is Hilbert and imbedded continuously in
$\mathcal{D}'(\Omega)$. Consider a weight function $\varrho:=\varrho_{1}^{-\sigma}$ such
that
\begin{equation}\label{eq12.4}
\varrho_{1}\in C^{\infty}(\,\overline{\Omega}\,),\;\;\varrho_{1}>0\;\;\mbox{in
$\Omega$},\;\;\varrho_{1}(x)=\mathrm{dist}(x,\partial\Omega)\;\;\mbox{near
$\partial\Omega$}.
\end{equation}

\begin{theorem}[the second Lions--Magenes theorem \cite{LionsMagenes72}]\label{thLM2}
Let $\sigma<0$ and
\begin{equation}\label{eq12.5}
X^{\sigma}(\Omega):=
\begin{cases}
\;\varrho_{1}^{-\sigma}\mathrm{H}^{\sigma}(\Omega)&\mbox{if}\;\;\sigma+1/2\notin\mathbb{Z},\\
\;\bigl[\,\varrho_{1}^{-\sigma+1/2}\,\mathrm{H}^{\sigma-1/2}(\Omega),\,
\varrho_{1}^{-\sigma-1/2}\,\mathrm{H}^{\sigma+1/2}(\Omega)\bigr]_{t^{1/2}}&
\mbox{if}\;\;\sigma+1/2\in\mathbb{Z}.
\end{cases}
\end{equation}
Then the conclusion of Theorem $\ref{thLM1}$ remains true.
\end{theorem}

\begin{remark}\label{rem12.2}
In the cited monograph \cite[Ch.~2, Sec. 6.3]{LionsMagenes72}, Lions and Magenes
introduced the space $X^{\sigma}(\Omega)$ in a way different from \eqref{eq12.5} and
designated $X^{\sigma}(\Omega)$ as $\Xi^{\sigma}(\Omega)$. Namely, for an integer
$\sigma\geq0$, the space $\Xi^{\sigma}(\Omega)$ is defined to consists of all
$f\in\mathcal{D}'(\Omega)$ such that $\varrho_{1}^{|\mu|}D^{\mu}f\in L_{2}(\Omega)$
for each multi-index $\mu$ with $|\mu|\leq\sigma$, and $\Xi^{\sigma}(\Omega)$ is
endowed with the Hilbert norm
$\sum_{|\mu|\leq\sigma}\|\varrho_{1}^{|\mu|}D^{\mu}f\|_{L_{2}(\Omega)}$. Then,
$\Xi^{\sigma}(\Omega):=
[\Xi^{[\sigma]}(\Omega),\Xi^{[\sigma]+1}(\Omega)]_{t^{\{\sigma\}}}$ for fractional
$\sigma>0$, with $\sigma=[\sigma]+\{\sigma\}$ and $[\sigma]$ being the integral part
of $\sigma$. Finally, $\Xi^{\sigma}(\Omega):=(\Xi^{-\sigma}(\Omega))'$ for
$\sigma<0$, the duality being with respect to the inner product in $L_{2}(\Omega)$.
It follows from the result of Lions and Magenes \cite[Ch.~2, Sec. 7.1, Corollary
7.4]{LionsMagenes72} that, for every $\sigma<0$, the space $\Xi^{\sigma}(\Omega)$
coincides with \eqref{eq12.5} up to equivalence of norms.
\end{remark}

\subsection{An extension of the Lions--Magenes theorems}\label{sec12.2}

The results presented here are got by the second author in \cite{09MFAT2}. First, we
establish the key theorem, which is a certain generalization of the Lions--Magenes
theorems stated above. The key theorem asserts that the operator \eqref{eq12.3} is
well defined, bounded, and Fredholm for $\sigma<0$ provided that a Hilbert space
$X^{\sigma}(\Omega)\hookrightarrow\mathcal{D}'(\Omega)$ satisfies the following
condition.

\begin{condition}[we name it as I$_{\sigma}$]\label{cond12.1}
The set $X^{\infty}(\Omega):=X^{\sigma}(\Omega)\cap
C^{\infty}(\,\overline{\Omega}\,)$ is dense in $X^{\sigma}(\Omega)$, and there
exists a number $c>0$ such that $\|\mathcal{O}f\|_{H^{\sigma}(\mathbb{R}^{n})}\leq
c\,\|f\|_{X^{\sigma}(\Omega)}$ for all $f\in X^{\infty}(\Omega)$, where
$\mathcal{O}f$ is defined by \eqref{eq8.15}.
\end{condition}

Note that the smaller $\sigma$ is, the weaker Condition \ref{cond12.1}
(I$_{\sigma}$) will be for the same space $X^{\sigma}(\Omega)$.

The spaces $X^{\sigma}(\Omega)$ appearing in Theorems \ref{thLM1} and \ref{thLM2}
satisfy Condition \ref{cond12.1}. This is evident for the first theorem, whereas,
for the second one, this follows from the dense continuous imbedding
$\mathrm{H}^{-\sigma}(\Omega)\hookrightarrow\Xi^{-\sigma}(\Omega)$ by the duality in
view of Theorem \ref{th8.3} iii) and Remark \ref{rem12.2}.

Our key theorem is the following.

\begin{theorem}\label{th12.1}
Let $\sigma<0$ and $X^{\sigma}(\Omega)$ be an arbitrary Hilbert space imbedded
continuously in $\mathcal{D}'(\Omega)$ and satisfying Condition $\ref{cond12.1}$
($\mathrm{I}_{\sigma}$). Then:
\begin{enumerate}
\item[i)] The set $D^{\infty}_{L,X}(\Omega):=\{u\in
C^{\infty}(\,\overline{\Omega}\,):Lu\in X^{\sigma}(\Omega)\}$ is dense in
$D^{\sigma+2q}_{L,X}(\Omega)$.
\item[ii)] The mapping $u\rightarrow(Lu,Bu)$, with $u\in
D^{\infty}_{L,X}(\Omega)$, extends uniquely to the continuous linear operator
\eqref{eq12.3}.
\item[iii)] The operator \eqref{eq12.3} is Fredholm. Its kernel is $\mathcal{N}$, and
its range consists of all the vectors
$(f,g_{1},\ldots,g_{q})\in\mathbf{X}_{\sigma}(\Omega,\partial\Omega)$ that satisfy
\eqref{eq9.8}.
\item[iv)] If $\mathcal{O}(X^{\infty}(\Omega))$ is dense in
$H^{\sigma}_{\overline{\Omega}}(\mathbb{R}^{n})$, then the index of \eqref{eq12.3}
is $\dim\mathcal{N}-\dim\mathcal{N}^{+}$.
\end{enumerate}
\end{theorem}

Let us outline the proof of Theorem \ref{th12.1}. The main idea is to derive this
theorem from Roitberg's generic theorem, i.e. from Theorem \ref{th11.2} considered
in the $\varphi\equiv1$ case. For the sake of simplicity, suppose that
$\mathcal{N}=\mathcal{N}^{+}=\{0\}$.

We get from Condition \ref{cond12.1} ($\mathrm{I}_{\sigma}$) that the mapping
$f\mapsto\mathcal{O}f$, $f\in X^{\infty}(\Omega)$, extends by a continuity to a
bounded linear injective operator $\mathcal{O}:X^{\sigma}(\Omega)\rightarrow
H^{\sigma}_{\overline{\Omega}}(\mathbb{R}^{n})$. This operator defines the
continuous imbedding $X^{\sigma}(\Omega)\hookrightarrow H^{\sigma,(0)}(\Omega)$.
Hence, by Theorem \ref{th11.2} a restriction of \eqref{eq11.8} establishes a
homeomorphism
\begin{equation}\label{eq12.6}
(L,B):D^{\sigma+2q,(2q)}_{L,X}(\Omega)\leftrightarrow
\mathbf{X}_{\sigma}(\Omega,\partial\Omega).
\end{equation}
Its domain is the Hilbert space
\begin{gather*}
D^{\sigma+2q,(2q)}_{L,X}(\Omega):=\{u\in H^{\sigma+2q,(2q)}(\Omega):Lu\in
X^{\sigma}(\Omega)\},\\
\|u\|_{D^{\sigma+2q,(2q)}_{L,X}(\Omega)}^{2}:=\|u\|_{H^{\sigma+2q,(2q)}(\Omega)}^{2}+
\|Lu\|_{X^{\sigma}(\Omega)}^{2}.
\end{gather*}
It follows from \eqref{eq12.6} that $D^{\infty}_{L,X}(\Omega)$ is dense in
$D^{\sigma+2q,(2q)}_{L,X}(\Omega)$.

According to Ya.A.~Roitberg \cite[Sec. 6.1, Theorem 6.1.1]{Roitberg96} we have the
equivalence of norms
\begin{equation}\label{eq12.7}
\|u\|_{H^{\sigma+2q,(2q)}(\Omega)}\asymp\bigl(\,\|u\|_{H^{\sigma+2q,(0)}(\Omega)}^{2}+
\|Lu\|_{H^{\sigma,(0)}(\Omega)}^{2}\,\bigr)^{1/2},\quad u\in
C^{\infty}(\,\overline{\Omega}\,).
\end{equation}
This result and the continuous imbedding $X^{\sigma}(\Omega)\hookrightarrow
H^{\sigma,(0)}(\Omega)$ imply
\begin{equation}\label{eq12.8}
\|u\|_{D^{\sigma+2q,(2q)}_{L,X}(\Omega)}\asymp
\bigl(\,\|u\|_{H^{\sigma+2q,(0)}(\Omega)}^{2}+
\|Lu\|_{X^{\sigma}(\Omega)}^{2}\,\bigr)^{1/2},\quad u\in
C^{\infty}(\,\overline{\Omega}\,).
\end{equation}
Thus, $D^{\sigma+2q,(2q)}_{L,X}(\Omega)$ is the completion of
$D^{\infty}_{L,X}(\Omega)$ with respect to the norm which is the right-hand side of
\eqref{eq12.8}.

Consider the mapping $u\mapsto u_{0}$ that takes each $u\in
D^{\sigma+2q,(2q)}_{L,X}(\Omega)$ to the initial component $u_{0}\in
H^{\sigma+2q,(0)}(\Omega)$ of the vector $T_{2q}u$. Here the operator $T_{2q}$ is
that in Theorem \ref{th11.1} for $r=2q$.

If $-2q-1/2\leq\sigma<0$, then
$H^{\sigma+2q,(0)}(\Omega)=\mathrm{H}^{\sigma+2q}(\Omega)$ by \eqref{eq12.2}. Now,
we may assert that the mapping $u\mapsto u_{0}$ establishes a homeomorphism of
$D^{\sigma+2q,(2q)}_{L,X}(\Omega)$ onto $D^{\sigma+2q}_{L,X}(\Omega)$. Hence,
\eqref{eq12.6} implies the required homeomorphism
\begin{equation}\label{eq12.9}
(L,B):D^{\sigma+2q}_{L,X}(\Omega)\leftrightarrow
\mathbf{X}_{\sigma}(\Omega,\partial\Omega).
\end{equation}

Further, if $\sigma<-2q-1/2$, then
$H^{\sigma+2q,(0)}(\Omega)=H^{\sigma+2q}_{\overline{\Omega}}(\mathbb{R}^{n})$. Then
using \eqref{eq12.1} and Roitberg's result \cite[Sec. 6.2, Theorem 6.2]{Roitberg96}
we can prove that the mapping $u\mapsto u_{0}\!\upharpoonright\!\Omega$ establishes
a homeomorphism of $D^{\sigma+2q,(2q)}_{L,X}(\Omega)$ onto
$D^{\sigma+2q}_{L,X}(\Omega)$. Hence, \eqref{eq12.6} implies \eqref{eq12.9} in this
case as well. See \cite[Sec. 4]{09MFAT2} for more details.

\begin{remark}\label{rem12.3}
A proposition similar to Theorem \ref{th12.1} was proved in Magenes's survey
\cite[Sec.~6.10]{Magenes65} for non half-integer $\sigma\leq-2q$ and the Dirichlet
problem, the space $X^{\sigma}(\Omega)$ obeying some different conditions depending
on the problem. Our Condition \ref{cond12.1} (I$_{\sigma}$) does not depend on it.
\end{remark}

\begin{remark}\label{rem12.4}
Ya.A.~Roitberg \cite[Sec.~2.4]{Roitberg71} considered a condition on the space
$X^{\sigma}(\Omega)$, which was somewhat stronger than Condition \ref{cond12.1}
(I$_{\sigma}$). He required additionally that
$C^{\infty}(\,\overline{\Omega}\,)\subset X^{\sigma}(\Omega)$. Under this stronger
condition, Roitberg \cite[Sec.~2.4]{Roitberg71}, \cite[Sec. 6.2, p. 190]{Roitberg96}
proved the boundedness of the operator \eqref{eq12.3} for all $\sigma<0$.
Homeomorphism Theorem for this operator was formulated in the survey \cite[Sec. 7.9,
p.~85]{Agranovich97} provided that $-2q\leq s\leq0$ and
$\mathcal{N}=\mathcal{N}^{+}=\{0\}$. We also mention the analogs of Theorem
\ref{th12.1} proved by Yu.V.~Kostarchuk and Ya.A.~Roitberg
\cite[Theorem~4]{KostarchukRoitberg73}, \cite[Sec.~1.3.8]{Roitberg99}. In these
analogs, Roitberg's condition is used, but solutions of an elliptic boundary-value
problem are considered in $H^{\sigma+2q,(2q)}(\Omega)$. Note that Roitberg's
condition does not include the important case where $X^{\sigma}(\Omega)=\{0\}$ and
does not cover some weighted spaces
$X^{\sigma}(\Omega)=\varrho\mathrm{H}^{\sigma}(\Omega)$, which we consider.
\end{remark}

Let us consider some applications of Theorem \ref{th12.1} caused by a particular
choice of the space $X^{\sigma}(\Omega)$. Apparently, the space
$X^{\sigma}(\Omega):=\{0\}$ satisfies Condition \ref{cond12.1} (I$_{\sigma}$). In
this case, Theorem \ref{th12.1} coincides with Theorem \ref{th10.1} for
$s:=\sigma+2q<2q$. It is remarkable that, despite $\mathrm{H}^{s}(\Omega)\neq
H^{s}(\Omega)$ for half-integer $s<0$, we have
\begin{equation}\label{eq12.10}
\{u\in\mathrm{H}^{s}(\Omega):\,Lu=0\;\;\mbox{in}\;\;\Omega\}=\{u\in
H^{s}(\Omega):\,Lu=0\;\;\mbox{in}\;\;\Omega\},
\end{equation}
the norms in $\mathrm{H}^{s}(\Omega)$ and $H^{s}(\Omega)$ being equivalent on the
distributions $u$ appearing in \eqref{eq12.10}.

It is also evident that the space $X^{\sigma}(\Omega):=L_{2}(\Omega)$ satisfies
Condition \ref{cond12.1} (I$_{\sigma}$) for every $\sigma<0$. In this important
case, Theorem \ref{th12.1} coincides with Theorem \ref{thLM1}.

We can describe all the Sobolev inner product spaces satisfying Condition
\ref{cond12.1}.

\begin{lemma}\label{lem12.1}
Let $\sigma<0$ and $\lambda\in\mathbb{R}$. The space
$X^{\sigma}(\Omega):=\mathrm{H}^{\lambda}(\Omega)$ satisfies Condition
$\ref{cond12.1}$ ($\mathrm{I}_{\sigma}$) if and only if
\,$\lambda\geq\max\,\{\sigma,-1/2\}$.
\end{lemma}

Indeed, we can restrict ourselves to the $\lambda<0$ case. Then the space
$X^{\sigma}(\Omega):=\mathrm{H}^{\lambda}(\Omega)$ satisfies Condition
\ref{cond12.1} (I$_{\sigma}$) if and only if the mapping $\mathcal{O}$ establishes
the dense continuous embedding $\mathrm{H}^{\lambda}(\Omega)\hookrightarrow
H^{\sigma}_{\overline{\Omega}}(\mathbb{R}^{n})$. By the duality, this embedding is
equivalent to the dense continuous embedding $H^{-\sigma}(\Omega)\hookrightarrow
H^{-\lambda}_{0}(\Omega)$, which is valid if and only if $-\sigma\geq-\lambda$ and
$H^{-\lambda}_{0}(\Omega)=H^{-\lambda}(\Omega)$. Since the latter equality
$\Leftrightarrow-\lambda\leq1/2$, the lemma is proved.

The next individual theorem results from Theorem \ref{th12.1} and Lemma
\ref{lem12.1}.

\begin{theorem}\label{th12.2}
Let $\sigma<0$ and $\lambda\geq\max\,\{\sigma,-1/2\}$. Then the mapping
$u\mapsto(Lu,Bu)$, with $u\in C^{\infty}(\,\overline{\Omega}\,)$, extends uniquely
to a continuous linear operator
\begin{equation}\label{eq12.11}
(L,B):\,\{u\in\mathrm{H}^{\sigma+2q}(\Omega):Lu\in
\mathrm{H}^{\lambda}(\Omega)\}\rightarrow
\mathrm{H}^{\lambda}(\Omega)\oplus\bigoplus_{j=1}^{q}\,H^{\sigma+2q-m_{j}-1/2}(\partial\Omega)
\end{equation}
provided that its domain is endowed with the norm
$$
\bigl(\,\|u\|_{\mathrm{H}^{\sigma+2q}(\Omega)}^{2}+
\|Lu\|_{\mathrm{H}^{\lambda}(\Omega)}^{2}\bigr)^{1/2}.
$$
The domain is a Hilbert space with respect to this norm. Moreover, the operator
\eqref{eq12.11} is Fredholm, and its index is $\dim\mathcal{N}-\dim\mathcal{N}^{+}$.
\end{theorem}

Here, it is useful to discuss the special case where $\lambda=\sigma$. If
$-1/2<\lambda=\sigma<0$, then the domain of \eqref{eq12.11} coincides with
$H^{\sigma+2q}(\Omega)$ and we arrive at Theorem \ref{th9.1} for $s=\sigma+2q$ and
$\varphi\equiv1$. If $\lambda=\sigma=-1/2$, then the domain is narrower than
$\mathrm{H}^{2q-1/2}(\Omega)$ and is equal to $H^{2q-1/2,(2q)}(\Omega)$ in view of
\eqref{eq12.8} and \eqref{eq12.2} so that we get Theorem \ref{th11.2} for $s=2q-1/2$
and $\varphi\equiv1$.

In Theorem \ref{th12.2}, we always have
$X^{\sigma}(\Omega)\subseteq\mathrm{H}^{-1/2}(\Omega)$. But we can get a space
$X^{\sigma}(\Omega)$ containing an extensive class of distributions
$f\notin\mathrm{H}^{-1/2}(\Omega)$ and satisfying Condition \ref{cond12.1}
(I$_{\sigma}$) if we use certain weighted spaces
$\varrho\mathrm{H}^{\sigma}(\Omega)$.

In this connection, recall the following.

\begin{definition}\label{def12.1}
Let $X(\Omega)$ be a Banach space lying in $\mathcal{D}'(\Omega)$. A function
$\varrho$ given in $\Omega$ is called a multiplier in $X(\Omega)$ if the operator of
multiplication by $\varrho$ is defined and bounded on $X(\Omega)$.
\end{definition}

Let $\sigma<-1/2$ and consider the next condition.

\begin{condition}[we name it as II$_{\sigma}$]\label{cond12.2}
The function $\varrho$ is a multiplier in $H^{-\sigma}(\Omega)$, and
\begin{equation}\label{eq12.12}
D_{\nu}^{j}\,\varrho=0\;\;\mbox{on}\;\;\partial\Omega\;\;\mbox{for every}\;\;
j\in\mathbb{Z}\;\;\mbox{such that}\;\;0\leq j<-\sigma-1/2.
\end{equation}
\end{condition}

Note if $\varrho$ is a multiplier in $H^{-\sigma}(\Omega)$, then evidently
$\varrho\in H^{-\sigma}(\Omega)$ so that, by Theorem \ref{th8.4}, the trace of
$D_{\nu}^{j}\varrho$ on $\partial\Omega$ is well defined in \eqref{eq12.12}. A
description of the set of all multipliers in $H^{-\sigma}(\Omega)$ is given in
\cite[Sec. 9.3.3]{MazyaShaposhnikova09}.

Using Condition \ref{cond12.2} (II$_{\sigma}$), we can describe the class of all
weighted Sobolev inner product spaces of order $\sigma$ that satisfies Condition
\ref{cond12.1} (I$_{\sigma}$).

\begin{lemma}\label{lem12.2}
Let $\sigma<-1/2$, and let a function $\varrho\in C^{\infty}(\Omega)$ be positive.
The space $X^{\sigma}(\Omega):=\varrho\mathrm{H}^{\sigma}(\Omega)$ satisfies
Condition $\ref{cond12.1}$ ($\mathrm{I}_{\sigma}$) if and only if $\varrho$ meets
Condition $\ref{cond12.2}$ ($\mathrm{II}_{\sigma}$).
\end{lemma}

Indeed, using the intrinsic description of $H^{-\sigma}_{0}(\Omega)$ mentioned in
Subsection \ref{sec8.4}, we can prove that $\varrho$ satisfies Condition
\ref{cond12.2} ($\mathrm{II}_{\sigma}$) if and only if the multiplication by
$\varrho$ is a bounded operator $M_{\varrho}:H^{-\sigma}(\Omega)\rightarrow
H^{-\sigma}_{0}(\Omega)$. The latter is equivalent, by the duality, to the
boundedness of the operator $M_{\varrho}:\mathrm{H}^{\sigma}(\Omega)\rightarrow
H^{\sigma}_{\overline{\Omega}}(\mathbb{R}^{n})$. Note that the mapping
$f\mapsto\varrho^{-1}f$ establishes the homeomorphism $M_{\varrho^{-1}}:
\varrho\mathrm{H}^{\sigma}(\Omega)\leftrightarrow\mathrm{H}^{\sigma}(\Omega)$.
Therefore, we conclude that $\varrho$ satisfies Condition \ref{cond12.2}
($\mathrm{II}_{\sigma}$) if and only if the identity operator
$M_{\varrho}\,M_{\varrho^{-1}}$ establishes a continuous embedding
$\mathcal{O}:\varrho\mathrm{H}^{\sigma}(\Omega)\rightarrow
H^{\sigma}_{\overline{\Omega}}(\mathbb{R}^{n})$. The embedding means that the space
$X^{\sigma}(\Omega)=\varrho\mathrm{H}^{\sigma}(\Omega)$ satisfies Condition
\ref{cond12.1} ($\mathrm{I}_{\sigma}$).

The next individual theorem results from Theorem \ref{th12.1} and Lemma
\ref{lem12.2}.

\begin{theorem}\label{th12.3}
Let $\sigma<-1/2$, and let a positive function $\varrho\in C^{\infty}(\Omega)$
satisfy Condition $\ref{cond12.2}$ ($\mathrm{II}_{\sigma}$). Then the mapping
$u\rightarrow(Lu,Bu)$, with $u\in C^{\infty}(\,\overline{\Omega}\,)$,
$Lu\in\varrho\mathrm{H}^{\sigma}(\Omega)$, extends uniquely to a continuous linear
operator
\begin{equation}\label{eq12.13}
(L,B):\,\bigl\{u\in\mathrm{H}^{\sigma+2q}(\Omega):
Lu\in\varrho\mathrm{H}^{\sigma}(\Omega)\bigr\}\rightarrow
\varrho\mathrm{H}^{\sigma}(\Omega)
\oplus\bigoplus_{j=1}^{q}\,H^{\sigma+2q-m_{j}-1/2}(\partial\Omega)
\end{equation}
provided that its domain is endowed with the norm
$$
\bigl(\,\|u\|_{\mathrm{H}^{\sigma+2q}(\Omega)}^{2}+
\|\varrho^{-1}Lu\|_{\mathrm{H}^{\sigma}(\Omega)}^{2}\bigr)^{1/2}.
$$
The domain is a Hilbert space with respect to this norm. Moreover, the operator
\eqref{eq12.13} is Fredholm, and its index is $\dim\mathcal{N}-\dim\mathcal{N}^{+}$.
\end{theorem}

We give an important example of a function $\varrho$ satisfying Condition
\ref{cond12.2} (II$_{\sigma}$) for fixed $\sigma<-1/2$ if we set
$\varrho:=\varrho_{1}^{\delta}$ provided that $\varrho_{1}$ meets \eqref{eq12.4} and
that $\delta\geq-\sigma-1/2\in\mathbb{Z}$ or $\delta>-\sigma-1/2\notin\mathbb{Z}$.

It is useful to compare Theorem \ref{thLM2} (the second Lions--Magenes theorem) with
Theorems \ref{th12.2} and \ref{th12.3}. For non half-integer $\sigma<-1/2$, Theorem
\ref{thLM2} is the special case of Theorem \ref{th12.3}, where
$\varrho:=\varrho_{1}^{-\sigma}$. For the half-integer values of $\sigma<-1/2$,
Theorem \ref{thLM2} follows from this case by the interpolation with the power
parameter $t^{1/2}$. Finally, if $-1/2\leq\sigma<0$, then Theorem \ref{thLM2} is a
consequence of Theorem \ref{th12.2}, in which we can take the space
$X^{\sigma}(\Omega):=\mathrm{H}^{\sigma}(\Omega)$ containing
$\varrho_{1}^{-\sigma}\mathrm{H}^{\sigma}(\Omega)$.

\subsection{Individual theorems on classes of H\"ormander spaces}\label{sec12.3}

Here we give analogs of Theorems \ref{th12.1}, \ref{th12.2}, and \ref{th12.3} for
some classes of H\"ormander spaces. The proofs of the analogs are similar to those
outlined in the previous subsection.

First, we state the key theorem, an analog of Theorems \ref{th12.1}. Let $\sigma<0$
and $\varphi\in\mathcal{M}$. Suppose that a Hilbert space
$X^{\sigma,\varphi}(\Omega)$ is embedded continuously in $\mathcal{D}'(\Omega)$.
Consider the following analog of Condition \ref{cond12.1} (I$_{\sigma}$).

\begin{condition}[we name it as I$_{\sigma,\varphi}$]\label{cond12.3}
The set $X^{\infty}(\Omega):=X^{\sigma,\varphi}(\Omega)\cap
C^{\infty}(\,\overline{\Omega}\,)$ is dense in $X^{\sigma,\varphi}(\Omega)$, and
there exists a number $c>0$ such that
$\|\mathcal{O}f\|_{H^{\sigma,\varphi}(\mathbb{R}^{n})}\leq
c\,\|f\|_{X^{\sigma,\varphi}(\Omega)}$ for all $f\in X^{\infty}(\Omega)$, where
$\mathcal{O}f$ is defined by \eqref{eq8.15}.
\end{condition}

The domain of $(L,B)$ is defined by the formula
$$
D^{\sigma+2q,\varphi}_{L,X}(\Omega):=\{u\in H^{\sigma+2q,\varphi}(\Omega):\,Lu\in
X^{\sigma,\varphi}(\Omega)\}
$$
and endowed with the graph inner product
$$
(u_{1},u_{2})_{D^{\sigma+2q}_{L,X}(\Omega)}:=
(u_{1},u_{2})_{H^{\sigma+2q}(\Omega)}+(Lu_{1},Lu_{2})_{X^{\sigma}(\Omega)}.
$$
The space $D^{\sigma+2q,\varphi}_{L,X}(\Omega)$ is Hilbert.

Our key theorem on classes of H\"ormander spaces is the following.

\begin{theorem}\label{th12.4}
Let $\varphi\in\mathcal{M}$, and let a number $\sigma<0$ be such that
\begin{equation}\label{eq12.14}
\sigma+2q\neq1/2-k\quad\mbox{for every integer}\quad k\geq1.
\end{equation}
Suppose that $X^{\sigma,\varphi}(\Omega)$ is an arbitrary Hilbert space imbedded
continuously in $\mathcal{D}'(\Omega)$ and satisfying Condition $\ref{cond12.3}$
($\mathrm{I}_{\sigma,\varphi}$). Then:
\begin{enumerate}
\item[i)] The set $D^{\infty}_{L,X}(\Omega):=\{u\in
C^{\infty}(\,\overline{\Omega}\,):Lu\in X^{\sigma,\varphi}(\Omega)\}$ is dense in
$D^{\sigma+2q,\varphi}_{L,X}(\Omega)$.

\item[ii)] The mapping $u\rightarrow(Lu,Bu)$, with $u\in
D^{\infty}_{L,X}(\Omega)$, extends uniquely to a continuous linear operator
\begin{equation}\label{eq12.15}
(L,B):\,D^{\sigma+2q,\varphi}_{L,X}(\Omega)\rightarrow
X^{\sigma,\varphi}(\Omega)\oplus
\bigoplus_{j=1}^{q}\,H^{\sigma+2q-m_{j}-1/2,\varphi}(\partial\Omega)
=:\mathbf{X}_{\sigma,\varphi}(\Omega,\partial\Omega),
\end{equation}
\item[iii)] The operator $\eqref{eq12.15}$ is Fredholm. Its kernel is $\mathcal{N}$, and
its range consists of all the vectors
$(f,g_{1},\ldots,g_{q})\in\mathbf{X}_{\sigma,\varphi}(\Omega,\partial\Omega)$ that
satisfy \eqref{eq9.8}.
\item[iv)] If $\mathcal{O}(X^{\infty}(\Omega))$ is dense in
$H^{\sigma,\varphi}_{\overline{\Omega}}(\mathbb{R}^{n})$, then the index of
\eqref{eq12.15} is $\dim\mathcal{N}-\dim\mathcal{N}^{+}$.
\end{enumerate}
\end{theorem}

Note that the condition \eqref{eq12.14} is stipulated by that, in definition of
$D^{\sigma+2q,\varphi}_{L,X}(\Omega)$, we use the space
$H^{\sigma+2q,\varphi}(\Omega)$, rather than an appropriate analog of
$\mathrm{H}^{\sigma+2q}(\Omega)$, which is different from
$H^{\sigma+2q,\varphi}(\Omega)$ if $\sigma+2q$ is negative and half-integer.

The following two individual theorems result from the key theorem. The first of them
is for nonweighted H\"ormander spaces
$X^{\sigma,\varphi}(\Omega):=H^{\lambda,\eta}(\Omega)$. In view of Theorem
\ref{th9.1}, we can confine ourselves to the $\sigma<-1/2$ case.

\begin{theorem}\label{th12.5}
Let $\sigma<-1/2$, the condition \eqref{eq12.14} be fulfilled, $\lambda>-1/2$, and
$\varphi,\eta\in\mathcal{M}$. Then the mapping $u\mapsto(Lu,Bu)$, with $u\in
C^{\infty}(\,\overline{\Omega}\,)$, extends uniquely to a continuous linear operator
\begin{equation}\label{eq12.16}
(L,B):\,\{u\in H^{\sigma+2q,\varphi}(\Omega):Lu\in
H^{\lambda,\eta}(\Omega)\}\rightarrow H^{\lambda,\eta}(\Omega)\oplus
\bigoplus_{j=1}^{q}\,H^{\sigma+2q-m_{j}-1/2,\varphi}(\partial\Omega)
\end{equation}
provided that its domain is endowed with the norm
$$
\bigl(\,\|u\|_{H^{\sigma+2q,\varphi}(\Omega)}^{2}+
\|Lu\|_{H^{\lambda,\eta}(\Omega)}^{2}\bigr)^{1/2}.
$$
The domain is a Hilbert space with respect to this norm. Moreover, the operator
\eqref{eq12.16} is Fredholm, and its index is $\dim\mathcal{N}-\dim\mathcal{N}^{+}$.
\end{theorem}

It is remarkable that, in this individual theorem, the solution and right-hand side
of the elliptic equation $Lu=f$ can be of different supplementary smoothness,
$\varphi$ and~$\eta$.

The second individual theorem is for weighted H\"ormander spaces
$X^{\sigma,\varphi}(\Omega):=\varrho H^{\sigma,\varphi}(\Omega)$, namely
\begin{gather*}
\varrho H^{\sigma,\varphi}(\Omega):=\{f=\varrho v:\,v\in
H^{\sigma,\varphi}(\Omega)\,\},\\
(f_{1},f_{2})_{\varrho H^{\sigma,\varphi}(\Omega)}:=
(\varrho^{-1}f_{1},\varrho^{-1}f_{2})_{H^{\sigma,\varphi}(\Omega)}.
\end{gather*}
Here $\sigma<-1/2$, $\varphi\in\mathcal{M}$, and the function $\varrho\in
C^{\infty}(\Omega)$ is positive. The space $\varrho H^{\sigma,\varphi}(\Omega)$ is
Hilbert.

\begin{theorem}\label{th12.6}
Let $\sigma<-1/2$, the condition \eqref{eq12.14} be valid, and
$\varphi\in\mathcal{M}$. Suppose that a positive function $\varrho\in
C^{\infty}(\Omega)$ is a multiplier in $H^{-\sigma,1/\varphi}(\Omega)$ and satisfies
\eqref{eq12.12}. Then the mapping $u\rightarrow(Lu,Bu)$, with $u\in
C^{\infty}(\,\overline{\Omega}\,)$, $Lu\in\varrho H^{\sigma,\varphi}(\Omega)$,
extends uniquely to a continuous linear operator
\begin{gather}\label{eq12.17}
(L,B):\bigl\{u\in H^{\sigma+2q,\varphi}(\Omega): Lu\in\varrho
H^{\sigma,\varphi}(\Omega)\bigr\}\rightarrow\notag \\
\varrho H^{\sigma,\varphi}(\Omega)
\oplus\bigoplus_{j=1}^{q}H^{\sigma+2q-m_{j}-1/2,\varphi}(\partial\Omega)
\end{gather}
provided that its domain is endowed with the norm
$$
\bigl(\,\|u\|_{H^{\sigma+2q,\varphi}(\Omega)}^{2}+
\|\varrho^{-1}Lu\|_{H^{\sigma,\varphi}(\Omega)}^{2}\bigr)^{1/2}.
$$
The domain is a Hilbert space with respect to this norm. Moreover, the operator
\eqref{eq12.17} is Fredholm, and its index is $\dim\mathcal{N}-\dim\mathcal{N}^{+}$.
\end{theorem}

We get a wide enough class of weight functions $\varrho$ satisfying the condition of
this theorem if we set $\varrho:=\varrho_{1}^{\delta}$, where $\varrho_{1}$ is
subject to \eqref{eq12.4} and $\delta>-\sigma-1/2$.

\section{Other results}\label{sec13}

In this section, we outline applications of H\"ormander spaces to other classes of
elliptic problems, namely to nonregular elliptic boundary-value problems,
par\-a\-me\-ter-elliptic problems, mixed elliptic problems, and elliptic systems. We
recall the statements of these problems and formulate theorems about properties of
the correspondent operators. As for Sobolev spaces, the Fredholm property and its
implications will be preserved for some classes of H\"ormander spaces. The theorems
stated below are deduced from the Sobolev case with the help of the interpolation
with an appropriate function parameter. We will not sketch the proofs and only will
refer to the authors' relevant papers.

\subsection{Nonregular elliptic boundary-value problems}\label{sec13.1}

Here we suppose that the boundary-value problem \eqref{eq9.1}, \eqref{eq9.2} is
elliptic in $\Omega$ but can be nonregular. This means that it satisfies conditions
i) and ii) of Definition \ref{def9.1} but need not meet condition iii). Theorems
\ref{th9.1}--\ref{th9.4} remain valid for this boundary-value problem except for the
description of the operator range and the index formula given in Theorem
\ref{th9.1}. The exception is caused by that the boundary-value problem need not
have a formally adjoint boundary-value problem in the class of differential
equations. A~version of Theorem \ref{th9.1} in this situation is the following.

\begin{theorem}\label{th13.1}
Let $s>m+1/2$ and $\varphi\in\mathcal{M}$. Then the bounded linear operator
\eqref{eq9.4} is Fredholm. Its kernel coincides with $\mathcal{N}$, whereas its
range consists of all the vectors
$(f,g_{1},\ldots,g_{q})\in\mathcal{H}_{s,\varphi}(\Omega,\partial\Omega)$ such that
the equality in \eqref{eq9.8} is fulfilled for each $v\in W$. Here $W$ is a certain
finite-dimensional space that lies in
$C^{\infty}(\,\overline{\Omega}\,)\times(C^{\infty}(\partial\Omega))^{q}$ and does
not depend on $s$ and $\varphi$. The index of \eqref{eq9.4} is $\dim\mathcal{N}-\dim
W$ and is also independent of $s$, $\varphi$.
\end{theorem}

The proof is given in \cite[Sec. 4]{06UMJ3}. Recall, if the boundary-value problem
\eqref{eq9.1}, \eqref{eq9.2} is regular elliptic, then $W=\mathcal{N}^{+}$

\begin{example}\label{ex13.1}
The oblique derivative problem for the Laplace equation:
\begin{equation}\label{eq13.1}
\Delta u=f\;\;\mbox{in}\;\;\Omega,\quad\quad\frac{\partial
u}{\partial\eta}=g\;\;\mbox{on}\;\;\partial\Omega.
\end{equation}
Here $\eta$ is an infinitely smooth field of unit vectors $\eta(x)$,
$x\in\partial\Omega$. Suppose that $\dim\Omega=2$, then the boundary-value problem
\eqref{eq13.1} is elliptic in $\Omega$, but it is nonregular provided
$\partial\Omega_{\eta}\neq\varnothing$. Here $\partial\Omega_{\eta}$ denotes the set
of all $x\in\partial\Omega$ such that $\eta(x)$ is tangent to $\partial\Omega$. If
$\overline{\Omega}$ is a disk, then the correspondent operator index equals
$2-\delta(\eta)/\pi$, where $\delta(\eta)$ is the increment of the angle between
$i:=(1,\,0)$ and $\eta(x)$ when $x$ goes counterclockwise around $\partial\Omega$;
see, e.g., \cite[Ch. 19, \S~4]{Mikhlin68}. Note if $\dim\Omega\geq3$ and
$\partial\Omega_{\eta}\neq\varnothing$, then the boundary-value problem
\eqref{eq13.1} is not elliptic at all.
\end{example}

Other examples of nonregular elliptic boundary-value problems are given in
\cite[Sec.~4]{Roitberg69}.

At the end of this subsection, we recall the following important result concerning
an arbitrary boundary-value problem \eqref{eq9.1}, \eqref{eq9.1} (see, e.g.,
\cite[Sec. 2.4]{Agranovich97}). If the corresponding operator \eqref{eq9.4} is
Fredholm for certain $s\geq2q$ with $\varphi\equiv1$, then this problem is elliptic
in $\Omega$, i.e., the above-mentioned conditions i) and ii) are satisfied.

\subsection{Parameter-elliptic problems}\label{sec13.2}

Such problems were distinguished by S. Agmon and L.~Nirenberg \cite{Agmon62,
AgmonNirenberg63}, M.S.~Agranovich and M.I.~Vishik \cite{AgranovichVishik64} as a
class of elliptic boundary-value problems that depend on a complex-valued parameter,
say $\lambda$, and possess the following remarkable property. Providing
$|\lambda|\gg1$, the operator correspondent to the problem establishes a
homeomorphism on appropriate pairs of Sobolev spaces, and moreover the operator norm
admits a two-sided a~priory estimate with constants independent of $\lambda$.
Parameter-elliptic problems were applied to the spectral theory of elliptic
operators and to parabolic equations. Some wider classes of parameter-elliptic
operators and boundary-value problems were investigated by M.S.~Agranovich
\cite{Agranovich90, Agranovich92}, R.~Denk, R.~Mennicken and L.R.~Vol\-e\-vich
\cite{DenkMennickenVolevich98, DenkMennickenVolevich01}, G.~Grubb
\cite[Ch.~2]{Grubb96}, A.N.~Kozhevnikov \cite{Kozhevnikov73, Kozhevnikov96,
Kozhevnikov97} (see also the surveys \cite{Agranovich94, Agranovich97}).

In this subsection, we give an application of H\"ormander spaces to
parameter-elliptic boundary-value problems considered by Agmon, Nirenberg, and
Agranovich, Vishik. Namely, we state a homeomorphism theorem on a class of
H\"ormander spaces and give a correspondent two-sided a~priory estimate for the
operator norm.

Recall the definition of the parameter-elliptic boundary-value problem. We consider
the nonhomogeneous boundary-value problem
\begin{equation}\label{eq13.2}
L(\lambda)\,u=f\quad\mbox{in}\quad\Omega,\quad\quad
B_{j}(\lambda)\,u=g_{j}\quad\mbox{on}\quad\partial\Omega,\quad j=1,\ldots,q,
\end{equation}
that depends on the parameter $\lambda\in\mathbb{C}$ as follows:
\begin{equation}\label{eq13.3}
L(\lambda):=\sum_{r=0}^{2q}\,\lambda^{2q-r}L_{r},\quad\quad
B_{j}(\lambda):=\sum_{r=0}^{m_{j}}\,\lambda^{m_{j}-r}B_{j,r}.
\end{equation}
Here $L_{r}=L_{r}(x,D)$, $x\in\overline{\Omega}$, and $B_{j,r}=B_{j,r}(x,D)$,
$x\in\partial\Omega$, are linear partial differential expressions of order $\leq r$
and with complex-valued infinitely smooth coefficients. As above, the integers $q$
and $m_{j}$ satisfy the equalities $q\geq1$ and $0\leq m_{j}\leq 2q-1$. Note that
$L(0)=L_{2q}$ and $B_{j}(0)=B_{j,m_{j}}$.

We associate certain homogeneous polynomials in $(\xi,\lambda)\in\mathbb{C}^{n+1}$
with partial differential expressions \eqref{eq13.3}. Namely, we set
$$
L^{(0)}(x;\xi,\lambda):=\sum_{r=0}^{2q}\,\lambda^{2q-r}L^{(0)}_{r}(x,\xi),
\quad\mbox{with}\quad
x\in\overline{\Omega},\;\xi\in\mathbb{C}^{n},\;\lambda\in\mathbb{C}.
$$
Here $L^{(0)}_{r}(x,\xi)$ is the principal symbol of $L_{r}(x,D)$ provided
$\mathrm{ord}\,L_{r}=r$, or $L^{(0)}_{r}(x,\xi)\equiv0$ if $\mathrm{ord}\,L_{r}<r$.
Similarly, for $j=1,\ldots,q$, we put
$$
B^{(0)}_{j}(x;\xi,\lambda):=\sum_{r=0}^{m_{j}}\,\lambda^{m_{j}-r}B^{(0)}_{j,r}(x,\xi),
\quad\mbox{with}\quad x\in\partial\Omega,\;\xi\in \mathbb{C}^{n},\;\lambda\in\mathbb{C}.
$$
Here $B^{(0)}_{j,r}(x,\xi)$ is the principal symbol of $B_{j,r}(x,D)$ provided
$\mathrm{ord}\,B_{j,r}=r$, or $B^{(0)}_{j,r}(x,\xi)\equiv0$ if $\mathrm{ord}\,B_{j,r}<r$.
Note that $L^{(0)}(x;\xi,\lambda)$ and $B^{(0)}_{j}(x;\xi,\lambda)$ are homogeneous
polynomials in $(\xi,\lambda)$ of the orders $2q$ and $m_{j}$ respectively.

Let $K$ be a fixed closed angle on the complex plain with vertex at the origin; here
we admits the case where $K$ degenerates into a ray.

\begin{definition}\label{def13.1}
The boundary-value problem \eqref{eq13.2} is called parameter-elliptic in the angle
$K$ if the following conditions are satisfied:
\begin{enumerate}
\item[i)] $L^{(0)}(x;\xi,\lambda)\neq0$ for each $x\in\overline{\Omega}$,
$\xi\in\mathbb{R}^{n}$, and $\lambda\in K$ whenever $|\xi|+|\lambda|\neq0$.
\item[ii)] Let $x\in\partial\Omega$, $\xi\in\mathbb{R}^{n}$, and $\lambda\in K$
be such that $\xi$ is tangent to $\partial\Omega$ at $x$ and that
$|\xi|+|\lambda|\neq0$. Then the polynomials $B^{(0)}_{j}(x;\xi+\tau\nu(x),\lambda)$
in $\tau$, $j=1,\ldots,q$,  are linearly independent modulo
$\prod_{j=1}^{q}(\tau-\tau^{+}_{j}(x;\xi,\lambda))$. Here
$\tau^{+}_{1}(x;\xi,\lambda),\ldots,\tau^{+}_{q}(x;\xi,\lambda)$ are all the
$\tau$-roots of $L^{(0)}(x;\xi+\tau\nu(x),\lambda)$ with $\mathrm{Im}\,\tau>0$, each
root being taken the number of times equal to its multiplicity.
\end{enumerate}
\end{definition}

\begin{remark}\label{rem13.1}
Condition ii) of Definition \ref{def13.1} is well stated in the sense that, for the
polynomial $L^{(0)}(x;\xi+\tau\nu(x),\lambda)$, the numbers of the $\tau$-roots with
$\mathrm{Im}\,\tau>0$ and of those with $\mathrm{Im}\,\tau<0$ are the same and equal
to $q$ if we take into account the roots multiplicity. Indeed, it follows from
condition i) that the partial differential expression
$$
L(x;D,D_{t}):=\sum_{r=0}^{2q}\,D_{t}^{2q-r}L_{r}(x,D),\quad x\in\overline{\Omega},
$$
is elliptic. Since the expression includes the derivation with respect to $n+1\geq3$
real arguments $x_{1},\ldots,x_{n},t$, its ellipticity is equivalent to the proper
ellipticity condition (see Remark \ref{rem9.1}). So, the $\tau$-roots of
$L^{(0)}(x;\xi+\tau\nu(x),\lambda)$ satisfy the indicated property.
\end{remark}

Let us give some instances of parameter-elliptic boundary-value problems \cite[Sec.
3.1 b)]{Agranovich97}.

\begin{example}\label{ex13.2}
Let differential expression $L(\lambda)$ satisfy condition i) of Definition
\ref{def13.1}. Then the Dirichlet boundary-value problem for the equation
$L(\lambda)=f$ is parameter-elliptic in the angle $K$. Here the boundary conditions
do not depend on the parameter $\lambda$.
\end{example}

\begin{example}\label{ex13.3}
The boundary-value problem
$$
\Delta u+\lambda^{2}u=f\;\;\mbox{in}\;\;\Omega,\quad\quad\frac{\partial
u}{\partial\nu}-\lambda u=g\;\;\mbox{on}\;\;\partial\Omega
$$
is parameter-elliptic in each angle $K_{\varepsilon}:=
\{\lambda\in\mathbb{C}:\,\varepsilon\leq|\mathrm{Im\,\lambda}|\leq\pi-\varepsilon\}$,
with $0<\varepsilon<\pi/2$, if the complex plane is slitted along the negative
semiaxis.
\end{example}

Further in this subsection the boundary-value problem \eqref{eq13.2} is supposed to
be parameter-elliptic in the angle $K$.

It follows from Definition \ref{def13.1} in view of Remark \ref{rem13.1} that the
boundary-value problem \eqref{eq13.2} is elliptic in $\Omega$ (and need not be
regular) provided $\lambda=0$. Since $\lambda$ is contained only in the lover order
terms of differential expressions $L(\lambda)$ and $B_{j}(\lambda)$, the problem is
elliptic in $\Omega$ for every $\lambda\in\mathbb{C}$. So, by Theorem \ref{th9.1},
we have the Fredholm bounded operator
\begin{equation}\label{eq13.4}
(L(\lambda),B(\lambda)):\,H^{s,\varphi}(\Omega)\rightarrow
\mathcal{H}_{s,\varphi}(\Omega,\partial\Omega)
\end{equation}
for each $s>m+1/2$, $\varphi\in\mathcal{M}$, and $\lambda\in\mathbb{C}$. The
operator index does not depend on $s$, $\varphi$, and on $\lambda$ because $\lambda$
influences only the lover order terms; see, e.g., \cite[Sec. 20.1, Theorem
20.1.8]{Hermander85}.  Moreover, since the boundary-value problem \eqref{eq13.2} is
parameter-elliptic in $K$, the operator \eqref{eq13.4} possesses the following
additional properties.

\begin{theorem}\label{th13.2}
\begin{enumerate}
\item[i)] There exists a number $\lambda_{0}>0$ such that for each $\lambda\in K$
with $|\lambda|\geq\nobreak\lambda_{0}$ and for any $s>m+1/2$,
$\varphi\in\mathcal{M}$, the operator \eqref{eq13.2} is a homeomorphism of
$H^{s,\varphi}(\Omega)$ onto $\mathcal{H}_{s,\varphi}(\Omega,\partial\Omega)$.
\item[ii)] Suppose that $s>2q$ and $\varphi\in\mathcal{M}$, then there is a
number $c=c(s,\varphi)\geq\nobreak1$ such that, for each $\lambda\in K$, with
$|\lambda|\geq\max\{\lambda_{0},1\}$, and for every $u\in H^{s,\varphi}(\Omega)$, we
have the following two-sided estimate
\begin{eqnarray}\label{eq13.5}
&&c^{-1}\bigl(\,\|u\|_{H^{s,\varphi}(\Omega)}+
|\lambda|^{s}\varphi(|\lambda|)\,\|u\|_{L_{2}(\Omega)}\,\bigr)\notag\\
&\leq&\|L(\lambda)u\|_{H^{s-2q,\varphi}(\Omega)}+
|\lambda|^{s-2q}\varphi(|\lambda|)\,\|L(\lambda)u\|_{L_{2}(\Omega)}\notag\\
&&+\,\sum_{j=1}^{q}\,
\bigl(\,\|B_{j}(\lambda)u\|_{H^{s-m_{j}-1/2,\varphi}(\partial\Omega)}\notag\\
&&+\,|\lambda|^{s-m_{j}-1/2}\varphi(|\lambda|)\,
\|B_{j}(\lambda)u\|_{L_{2}(\partial\Omega)}\,\bigr)\notag\\
&\leq& c\,\bigl(\,\|u\|_{H^{s,\varphi}(\Omega)}+
|\lambda|^{s}\varphi(|\lambda|)\,\|u\|_{L_{2}(\Omega)}\,\bigr).
\end{eqnarray}
Here $c$ does not depend on $u$ and $\lambda$.
\end{enumerate}
\end{theorem}

We should comment on assertion ii) of this theorem. For fixed $\lambda$, the
estimate \eqref{eq13.5} is written for the norms, non-Hilbert, that are equivalent
to $\|u\|_{H^{s,\varphi}(\Omega)}$ and
$\|(L(\lambda),B(\lambda))u\|_{\mathcal{H}_{s,\varphi}(\Omega,\partial\Omega)}$
respectively. The non-Hilbert norms are used to avoid cumbersome expressions. To
have the finite norm $\|L(\lambda)u\|_{L_{2}(\Omega)}$ in \eqref{eq13.5}, we suppose
that $s>2q$ is fulfilled instead of the condition $s>m+1/2$ used in assertion i).
Finally, the supplement condition $|\lambda|\geq1$ is caused by that the function
$\varphi(t)$ is defined for $t\geq1$. Note the estimate \eqref{eq13.5} is of
interest for $|\lambda|\gg1$ only.

In the Sobolev case where $s\geq2q$ and $\varphi\equiv1$, Theorem \ref{th13.2} was
proved by M.S.~Agranovich and M.I.~Vishik \cite[\S~4 and 5]{AgranovichVishik64}; see
also \cite[Sec. 3.2]{Agranovich97}. In general, the theorem is proved in
\cite[Sec.~7]{07UMJ5}. Note that the right-hand side of the estimate \eqref{eq13.5}
is valid without the assumption about the parameter-ellipticity of \eqref{eq13.2}.
Analogs of Theorem \ref{th13.2} for parameter-elliptic operators, scalar or matrix,
are proved in \cite{07Dop5, 07UMJ6, 08MFAT2}.

We note an important consequence of Theorem \ref{th13.2} i). Suppose that the
boun\-dary-value problem \eqref{eq13.2} is parameter-elliptic on a certain ray
$K:=\{\lambda\in\mathbb{C}:\arg\lambda =\mathrm{const}\}$. Then the operator
\eqref{eq13.4} is of zero index for each $s>m+1/2$, $\varphi\in\mathcal{M}$, and
$\lambda\in\mathbb{C}$.

\subsection{Mixed elliptic problems}

Here we consider a certain class of elliptic boun\-dary-value problems in multiply
connected bonded domains. As distinguished from the above, we allow the orders of
the boundary differential expressions to be distinct on different connected
components of the boundary. For instance, studying the Laplace equation in a ring,
one may set the Dirichlet condition on a chosen connected component of the ring
boundary and the Neumann condition on the other component. The problems under
consideration relate to the mixed elliptic boundary-value problems \cite{Peetre61,
Schechter60, Simanca87, VishikEskin69}. They have not investigated so completely as
the unmixed elliptic problems. This is concerned with some difficulties, that appear
when one reduces the mixed problem to a pseudodifferential operator on the boundary;
see, e.g., \cite{Simanca87}. In the problems we consider, the portions of boundary
on which the boundary expression has distinct orders do not adjoin to each other.
These problems are called formally mixed. They can be reduced locally to a model
elliptic problem in the half-space \cite{07Dop4}.

In this subsection, we suppose that the boundary of $\Omega$ consists of $r\geq2$
nonempty connected components $\Gamma_{1},\ldots,\Gamma_{r}$. Fix an integer
$q\geq1$ and consider a formally mixed boundary-value problem
\begin{equation}\label{eq13.6}
L\,u=f\quad\text{in}\;\;\Omega,\quad
B^{(k)}_{j}u=g_{k,j}\;\;\text{on}\;\;\Gamma_{k},\;\;j=1,\ldots,q,\;\;k=1,\ldots,r.
\end{equation}
Here the partial differential expression $L=L(x,D)$, $x\in\overline{\Omega}$, of
order $2q$, is the same as in Section \ref{sec9}, whereas
$B^{(k)}:=\{B^{(k)}_{j}:j=1,\ldots,q\}$ is a system of boundary linear partial
differential expressions given on the component $\Gamma_{k}$. Suppose that the
coefficients of the expressions $B^{(k)}_{j}=B^{(k)}_{j}(x,D)$, $x\in\Gamma_{k}$,
are infinitely smooth complex-valued functions and that all
$m^{(k)}_{j}:=\mathrm{ord}\,B^{(k)}_{j}\leq2q-1$. We denote
\begin{gather*}
\Lambda:=(L,B^{(1)}_{1},\ldots,B^{(1)}_{q},\ldots,B^{(r)}_{1},\ldots,B^{(r)}_{q}), \\
\mathcal{N}_{\Lambda}:=\{u\in C^{\infty}(\,\overline{\Omega}\,):\,\Lambda u=0\}, \\
m:=\max\,\{\mathrm{ord}\,B^{(k)}_{j}:\,j=1,\ldots,q,\;\;k=1,\ldots,r\}.
\end{gather*}

The mapping $u\mapsto\Lambda u$, $u\in C^{\infty}(\,\overline{\Omega}\,)$, extends
uniquely to a bounded linear operator
\begin{gather}\label{eq13.7}
\Lambda:\,H^{s,\varphi}(\Omega)\rightarrow
H^{s-2q,\;\varphi}(\Omega)\oplus\bigoplus_{k=1}^{r}\bigoplus_{j=1}^{q}
H^{s-m^{(k)}_{j}-1/2,\varphi}(\Gamma_{k})\\
=:\mathcal{H}_{s,\varphi}(\Omega,\Gamma_{1},\ldots,\Gamma_{r}) \notag
\end{gather}
for each $s>m+1/2$ and $\varphi\in\mathcal{M}$.

\begin{definition}\label{def13.2}
The formally mixed boundary-value problem \eqref{eq13.6} is called elliptic in the
multiply connected domain $\Omega$ if $L$ is proper elliptic on $\overline{\Omega}$
and if, for each $k=1,\ldots,r$, the system $B^{(k)}$ satisfies the Lopatinsky
condition with respect to $L$ on $\Gamma_{k}$.
\end{definition}

Suppose the mixed boundary-value problem \eqref{eq13.6} is elliptic in $\Omega$.
Then it has the following properties \cite{07Dop4}.

\begin{theorem}\label{th13.3}
Let $s>m+1/2$ and $\varphi\in\mathcal{M}$. Then the bounded linear operator
\eqref{eq13.7} is Fredholm. Its kernel coincides with $\mathcal{N}_{\Lambda}$,
whereas its range consists of all the vectors
$$
(f,g_{1,1},\ldots,g_{1,q},\ldots,g_{r,1},\ldots,g_{r,q})\in
\mathcal{H}_{s,\varphi}(\Omega,\Gamma_{1},\ldots,\Gamma_{r})
$$
such that
\begin{equation}\label{eq13.8}
(f,w_{0})_{\Omega}+\sum_{k=1}^{r}\,\sum_{j=1}^{q}\; (g_{k,j},w_{k,j})_{\Gamma_{k}}=0
\end{equation}
for each vector-valued function
$$
(w_{0},w_{1,1},\ldots,w_{1,q},\ldots,w_{r,1},\ldots,w_{r,q})\in W_{\Lambda}.
$$
Here $W_{\Lambda}$ is a certain finite-dimensional space that lies in
$$
C^{\infty}(\,\overline{\Omega}\,)\times\prod_{j=1}^{r}\,(C^{\infty}(\Gamma_{j}))^{q}
$$
and does not depend on $s$ and $\varphi$. The index of \eqref{eq13.7} is
$\dim\mathcal{N}-\dim W_{\Lambda}$ and is also independent of $s$,~$\varphi$.
\end{theorem}

It is self-clear that, in \eqref{eq13.8}, the notation $(\cdot,\cdot)_{\Gamma_{k}}$
stands for the inner product in $L_{2}(\Gamma_{k})$.

\subsection{Elliptic systems}\label{sec13.4}

Extensive classes of elliptic systems of linear partial differential equations were
introduced and investigated by I.G.~Petrovskii \cite{Petrovskii39} and A.~Douglis,
L.~Nirenberg \cite{DouglisNirenberg55}. For pseudodifferential equations, general
elliptic systems were studied by L.~H\"{o}rmander \cite[Sec. 1.0]{Hermander67}. He
proved a priori estimates for solutions of these systems in appropriate couples of
Sobolev inner product spaces of arbitrary real orders. If the system is given on a
closed smooth manifold, then the estimate is equivalent to the Fredholm property of
the correspondent elliptic matrix PsDO; see, e.g., the monograph \cite[Ch.
19]{Hermander85}, and the survey \cite[Sec. 3.2]{Agranovich94}. This fact is of
great importance in the theory of elliptic boundary-value problems because each of
these problems can be reduced to an elliptic system of pseudodifferential equations
on the boundary of the domain; see, e.g., \cite[Ch.~20]{Hermander85} and
\cite[Part~IV]{WlokaRowleyLawruk95}.

In this subsection, we examine the Petrovskii elliptic systems on the refined
Sobolev scale over a closed smooth manifold $\Gamma$ and generalize the results of
Subsection \ref{sec6.1} to these systems.

Let us consider a system of $p\geq2$ linear equations
\begin{equation}\label{eq13.9}
\sum_{k=1}^{p}\:A_{j,k}\,u_{k}=f_{j}\quad\mbox{on}\quad\Gamma,\quad j=1,\ldots,p.
\end{equation}
Here $A_{j,k}$, $j,k=1,\ldots,p$, are scalar classical pseudodifferential operators
of arbitrary real orders defined on $\Gamma$. We consider equations \eqref{eq13.9}
in the sense of the distribution theory so that
$u_{k},\,f_{j}\in\mathcal{D}'(\Gamma)$. Put
$m_{k}:=\max\{\mathrm{ord}\,A_{1,k},\ldots,\mathrm{ord}\,A_{p,k}\}$.

Let us rewrite the system \eqref{eq13.9} in the matrix form: $Au=f$ on $\Gamma$,
where $A:=(A_{j,k})$\; is a square matrix of order $p$, and
$u=\mathrm{col}\,(u_{1},\ldots,u_{p})$, $f=\mathrm{col}\,(f_{1},\ldots,f_{p})$ are
functional columns. The mapping $u\mapsto Au$ is a linear continuous operator on the
space $(\mathcal{D}'(\Gamma))^{p}$. By lemma \ref{lem6.1}, a restriction of the
mapping sets a bounded linear operator
\begin{equation}\label{eq13.10}
A:\,\bigoplus_{k=1}^{p}\,H^{s+m_{k},\,\varphi}(\Gamma)\rightarrow
(H^{s,\varphi}(\Gamma))^{p}
\end{equation}
for each $s\in\mathbb{R}$ and $\varphi\in\mathcal{M}$.

\begin{definition}\label{def13.3}
The system \eqref{eq13.9} and the matrix PsDO $A$ are called Petrovskii elliptic on
$\Gamma$ if $\det\bigl(a^{(0)}_{j,k}(x,\xi)\bigr)_{j,k=1}^{p}\neq0$ for each point
$x\in\Gamma$ and covector $\xi\in T^{\ast}_{x}\Gamma\setminus\{0\}$. Here
$a_{j,k}^{(0)}(x,\xi)$ is the principal symbol of $A_{j,k}$ provided
$\mathrm{ord}\,A_{j,k}=m_{k}$; otherwise $a_{j,k}^{(0)}(x,\xi)\equiv0$.
\end{definition}

We suppose that the system $Au=f$ is elliptic on $\Gamma$. Then both the spaces
\begin{gather*}
N:=\bigl\{\,u\in(C^{\infty}(\Gamma))^{p}: \,Au=0\;\;\mbox{on}\;\;\Gamma\,\bigr\},\\
N^{+}:=\bigl\{v\in(C^{\infty}(\Gamma))^{p}:\,A^{+}v=0
\;\;\mbox{on}\;\;\Gamma\,\bigr\}
\end{gather*}
are finite-dimensional \cite[Sec. 3.2]{Agranovich94}. Here $A^{+}$ is the matrix
pseudodifferential operator formally adjoint to $A$ with respect to the inner
product in $(L_{2}(\Gamma))^{p}$.

\begin{theorem}\label{th13.4}
The operator \eqref{eq13.10} corresponding to the elliptic system is Fredholm for
each $s\in\mathbb{R}$ and $\varphi\in\mathcal{M}$. Its kernel coincides with $N$,
whereas its range consists of all the vectors $f\in(H^{s,\varphi}(\Gamma))^{p}$ such
that $\sum_{j=1}^{p}\,(f_{j},v_{j})_{\Gamma}=0$ for each $(v_{1},\ldots,v_{p})\in
N^{+}$. The index of \eqref{eq13.10} is equal to $\dim N-\dim N^{+}$ and independent
of $s$ and $\varphi$.
\end{theorem}

This theorem is proved in \cite{08BPAS3} together with other properties of the
system \eqref{eq13.9}. They are similar to that given in Subsection \ref{sec6.1}, in
which the scalar case is treated. We also refer to the second author's papers
\cite{07Dop5, 08MFAT2, 08UMB3, 09UMJ3} devoted to various classes of elliptic
systems in H\"ormander spaces.

\subsection{Boundary-value problems for elliptic systems}\label{sec13.5}

Boundary-value problems for various classes of elliptic systems of linear partial
differential equations were investigated by S.~Agmon, A.~Douglis, and L.~Nirenberg,
M.S.~Agranovich and A.S.~Dynin, L.~H\"ormander, L.N.~Slobodetskii, V.A.~Solonnikov,
L.R.~Volevich; see the until now unique monograph \cite{WlokaRowleyLawruk95} devoted
especially to these problems, the survey \cite[\S~6]{Agranovich97} and the
references given therein. It was proved that the operator correspondent to the
problem is Fredholm on appropriate pairs of the positive order Sobolev spaces.
Regarding the boundary-value problems for Petrovskii elliptic systems, we extend
this result over the one-sided refined Sobolev scale.

Let us consider a system of $p\geq2$ partial differential equations
\begin{equation}\label{eq13.11}
\sum_{k=1}^{p}\,L_{j,k}\,u_{k}=f_{j}\quad\mbox{in}\quad\Omega,\quad j=1,\ldots,p.
\end{equation}
Here $L_{j,k}=L_{j,k}(x,D)$, $x\in\overline{\Omega}$, $j,k=1,\ldots,p$, are scalar
linear partial differential expressions given on $\overline{\Omega}$. The expression
$L_{j,k}$ is of an arbitrary finite order, the coefficients of $L_{j,k}$ are
supposed to be complex-valued and infinitely smooth on $\overline{\Omega}$. Put
$m_{k}:=\max\{\mathrm{ord}\,L_{1,k},\ldots,\mathrm{ord}\,L_{p,k}\}$ so that $m_{k}$
is the maximal order of derivative of the unknown function $u_{k}$. Suppose that all
$m_{k}\geq1$ and that $\sum_{k=1}^{p}\,m_{k}$ is even, say $2q$.

We consider the solutions of \eqref{eq13.11} that satisfy the boundary conditions
\begin{equation}\label{eq13.12}
\sum_{k=1}^{p}\,B_{j,k}\,u_{k}=g_{j}\quad\mbox{on}\quad\partial\Omega,\quad
j=1,\ldots,q.
\end{equation}
Here $B_{j,k}=B_{j,k}(x,D)$, with $x\in\partial\Omega$, $j=1,\ldots,q$, and
$k=1,\ldots,p$, are boundary linear partial differential expressions with infinitely
smooth coefficients. We suppose $\mathrm{ord}\,B_{j,k}\leq m_{k}-1$ and set
$r_{j}:=\min\,\{m_{k}-\mathrm{ord}\,B_{j,k}:\,k=1,\ldots,p\}$ admitting
$\mathrm{ord}\,B_{j,k}:=-\infty$ for $B_{j,k}\equiv0$; thus
$\mathrm{ord}\,B_{j,k}\leq m_{k}-r_{j}$.

Let us write the boundary-value problem \eqref{eq13.11}, \eqref{eq13.12} in the
matrix form
$$
Lu=f\;\;\mbox{in}\;\;\Omega,\quad Bu=g\;\;\mbox{on}\;\;\partial\Omega.
$$
Here $L:=(L_{j,k})_{j,k=1}^{p}$ and $B:=(B_{j,k})_{\substack{j=1,\ldots,q
\\ k=1,\ldots,p}}$ are matrix differential expressions, whereas
$u:=\mathrm{col}\,(u_{1},\ldots,u_{p})$, $f:=\mathrm{col}\,(f_{1},\ldots,f_{p})$, and
$g:=\mathrm{col}\,(g_{1},\ldots,g_{q})$ are function columns.

It follows from Lemma \ref{eq9.1} that the mapping $u\mapsto(Lu,Bu)$,
$u\in(C^{\infty}(\,\overline{\Omega}\,))^{p}$, extends uniquely to a continuous
linear operator
\begin{gather}\label{eq13.13}
(L,B):\,\bigoplus_{k=1}^{p}H^{s+m_{k},\varphi}(\Omega)\rightarrow
(H^{s,\varphi}(\Omega))^{p}\oplus\bigoplus_{j=1}^{q}
H^{s+r_{j}-1/2,\varphi}(\partial\Omega)\\
=:\mathbf{H}_{s,\varphi}(\Omega,\partial\Omega) \notag
\end{gather}
for each $s>-r+1/2$ and $\varphi\in\mathcal{M}$, with
$r:=\min\{r_{1},\ldots,r_{q}\}\geq1$. We are interested in properties of this
operator provided the boundary-value problem is elliptic in the Petrovskii sense.
Recall the ellipticity definition.

With $L$ and $B$ we associate the matrixes of homogeneous polynomials
$$
L^{(0)}(x,\xi):=\bigl(L_{j,k}^{(0)}(x,\xi)\bigr)_{j,k=1}^{p},\quad B^{(0)}(x,\xi):=
\bigl(B_{j,k}^{(0)}(x,\xi)\bigr)_{\substack{j=1,\ldots,q
\\ k=1,\ldots,p}}.
$$
Here $L_{j,k}^{(0)}(x,\xi)$, $x\in\overline{\Omega}$, $\xi\in\mathbb{C}^{n}$, is the
principal symbol of $L_{j,k}(x,D)$ provided $\mathrm{ord}\,L_{j,k}=m_{k}$; otherwise
$L_{j,k}^{(0)}(x,\xi)\equiv0$. Similarly, $B_{j,k}^{(0)}(x,\xi)$, $x\in\partial\Omega$,
$\xi\in\mathbb{C}^{n}$, is the principal symbol of $B_{j,k}(x,D)$ provided
$\mathrm{ord}\,B_{j,k}=m_{k}-r_{j}$; otherwise $B_{j,k}^{(0)}(x,\xi)\equiv0$.

\begin{definition}\label{def13.4}
The boundary-value problem \eqref{eq13.11}, \eqref{eq13.12} is called Petrovskii
elliptic in $\Omega$ if the following conditions are satisfied:
\begin{enumerate}
\item[i)] System \eqref{eq13.11} is proper elliptic on $\overline{\Omega}$; i.e.,
condition i) of Definition \ref{def9.1} is fulfilled, with the notation $\det
L^{(0)}(x,\xi'+\tau\xi'')$ being placed instead of $L^{(0)}(x,\xi'+\tau\xi'')$.
\item[ii)] Relations \eqref{eq13.12} satisfies the Lopatinsky condition with
respect to \eqref{eq13.11} on $\partial\Omega$; i.e., for an arbitrary point
$x\in\partial\Omega$ and for each vector $\xi\neq0$ tangent to $\partial\Omega$ at
$x$, the rows of the matrix $B^{(0)}(x,\xi+\tau\nu(x))\times
L^{(0)}_{\mathrm{c}}(x,\xi+\tau\nu(x))$ are linearly independent polynomials, in
$\tau\in\mathbb{R}$, modulo
$\prod_{j=1}^{q}\bigl(\tau-\tau^{+}_{j}(x;\xi,\nu(x))\bigr)$. Here
$L^{(0)}_{\mathrm{c}}(x,\xi)$ is the transpose of the matrix composed by the
cofactors of the matrix $L^{(0)}(x,\xi)$ elements.
\end{enumerate}
\end{definition}

Note, if condition i) is satisfied, then the system \eqref{eq13.11} is Petrovskii
elliptic on $\overline{\Omega}$, i.e. $\det L^{(0)}(x,\xi)\neq0$ for each
$x\in\overline{\Omega}$ and $\xi\in\mathbb{R}^{n}\setminus\{0\}$. The converse is
true provided that $\dim\Omega\geq3$; see \cite[Sec. 6.1~a)]{Agranovich97}.

\begin{example}\label{ex13.4}
The elliptic boundary-value problem for the Cauchy-Riemann system:
\begin{gather*}
\frac{\partial u_{1}}{\partial x_{1}}-\frac{\partial u_{2}}{\partial
x_{2}}=f_{1},\quad \frac{\partial u_{1}}{\partial x_{2}}+\frac{\partial
u_{2}}{\partial x_{1}}=f_{2}\quad\mbox{in}\quad\Omega,\\
u_{1}+u_{2}=g\quad\mbox{on}\quad\partial\Omega.
\end{gather*}
Here $n=p=2$ and $m_{1}=m_{2}=1$, so that $q=1$. The Cauchy-Riemann system is an
instance of homogeneous elliptic systems, which satisfy Definition \ref{def13.4}
with $m_{1}=\ldots=m_{p}$.
\end{example}

\begin{example}\label{ex13.5}
The Petrovskii elliptic boundary-value problem
\begin{gather*}
\frac{\partial u_{1}}{\partial x_{1}}-\frac{\partial^{3}u_{2}}{\partial
x_{2}^{3}}=f_{1},\quad \frac{\partial u_{1}}{\partial x_{2}}+\frac{\partial^{3}
u_{2}}{\partial x_{1}^{3}}=f_{2}\quad\mbox{in}\quad\Omega, \\
u_{1}=g_{1},\quad u_{2}\,\biggl(\mbox{or}\;\frac{\partial
u_{2}}{\partial\nu},\;\mbox{or}\;\frac{\partial^{2}
u_{2}}{\partial\nu^{2}}\biggr)=g_{2}\quad\mbox{on}\quad\partial\Omega.
\end{gather*}
Here $n=p=2$, $m_{1}=1$, and $m_{2}=3$ so that  $q=2$. This system is not
homo\-ge\-ne\-ous elliptic.
\end{example}

Other examples of elliptic systems, of various kinds, are given in
\cite[\S~6.2]{Agranovich97}.

Suppose the boundary-value problem \eqref{eq13.11}, \eqref{eq13.12} is Petrovskii
elliptic in $\Omega$. Then it has the following properties \cite{07Dop6}.

\begin{theorem}\label{th13.5}
Let $s>-r+1/2$ and $\varphi\in\mathcal{M}$. Then the bounded linear operator
\eqref{eq13.13} is Fredholm. The kernel $\mathcal{N}$ of \eqref{eq13.13} lies in
$(C^{\infty}(\,\overline{\Omega}\,))^{p}$ and does not depend on $s$ and $\varphi$.
The range of \eqref{eq13.13} consists of all the vectors
$(f_{1},\ldots,f_{p};g_{1},\ldots,g_{q})\in
\mathbf{H}_{s,\varphi}(\Omega,\partial\Omega)$ such that
$$
\sum_{j=1}^{p}\,(f_{j},w_{j})_{\Omega}+\sum_{j=1}^{q}\,
(g_{j},h_{j})_{\partial\Omega}=0
$$
for each vector-valued function $(w_{1},\ldots,w_{p};\,h_{1},\ldots,h_{q})\in W$.
Here $W$ is a certain finite-dimensional space that lies in
$(C^{\infty}(\,\overline{\Omega}\,))^{p}\times (C^{\infty}(\Gamma))^{q}$. The index
of the operator \eqref{eq13.13} is $\dim\mathcal{N}-\dim W$ and independent of
$s$,~$\varphi$.
\end{theorem}

\medskip

\bibliographystyle{amsplain}

\end{document}